\documentclass[reqno,fleqn]{amsart}
\usepackage{amsfonts,amsthm,amssymb}
\usepackage{lipsum}
\usepackage[utf8]{inputenc}
\usepackage[english]{babel}

\usepackage{lmodern,textcomp}
\usepackage{amsmath}
\usepackage{mathtools}
\usepackage{latexsym}
\usepackage{tikz}
\usepackage{fancyvrb}
\usepackage{epsfig}
\usepackage{pstricks,slashbox,multirow}
\usepackage{graphicx}
\usepackage{rotating}
\usetikzlibrary{positioning}
\usepackage{subcaption}
\usepackage{datetime}
\newtheorem{theorem}{Theorem}

\newtheorem{cor}[theorem]{Corollary}
\newtheorem{definition}[theorem]{Definition}

\newtheorem{prm}{Problem}

\newtheorem{rem}[theorem]{Remark}

\title[A study on Type-2 isomorphic $C_n(R)$: Part 2: Type-2 isomorphic $C_{n}(R)$ for $n$ = 16,24,27]{A study on Type-2 isomorphic circulant graphs. \\ Part 2: Type-2 isomorphic circulant graphs of orders 16, 24, 27}

\author{\sc Vilfred Kamalappan} 
\address{Department of Mathematics, Central University of Kerala, Periye, Kasaragod, Kerala, India - 671 316.}
\email{vilfredkamal@gmail.com}

\subjclass[2010]{05C60, 05C25, 05C75.}

\keywords{Circulant graph, Cayley Isomorphism (CI) property, Type-1 isomorphism, Type-2 isomorphism, Type-1 group of $C_{n}(R)$, Type-2 group of $C_{n}(R)$ w.r.t. $m$, $(T2_{n,m}(C_n(R)), ~\circ)$, $(V_{n,m}(C_n(R)), ~\circ)$.}

\date{}

\begin{document}

\begin{abstract} This article is the $2^{nd}$ part of a detailed study on Type-2 isomorphic circulant graphs having ten parts \cite{v2-1}-\cite{v2-10}. Definition of Type-2 isomorphism of circulant graphs $C_n(R)$ w.r.t. $m$ modified in \cite{v2-2-arX} was further modified to definition \ref{d6} in \cite{v2-1} by the author by considering $m > 1$ divides $\gcd(n, r)$, $m^3$ divides $n$ and $r\in R$ and studied Type-2 isomorphic circulant graphs w.r.t. $m$ = 2. In this paper, using the modified definition \ref{d6}, we obtain all 8 and 32 number of pairs of Type-2 isomorphic circulant graphs of orders 16 and 24, respectively and all 27 number of triples of Type-2 isomorphic circulant graphs of order 27.
\end{abstract}

\maketitle

	
\section{Introduction}

Circulant graphs $C_n(R)$ and $C_n(S)$ are said to be \emph{Adam's isomorphic} if there exist some $a\in \mathbb{Z}_n^*$ such that $S = a R$ under arithmetic reflexive modulo $n$ \cite{ad67}. We call Adam's isomorphism as {\em Type-1 isomorphism}. In 1996, Vilfred \cite{v96} defined Type-2 isomorphism of $C_n(R)$ w.r.t. $m$ $\ni$ $m$ = $\gcd(n, r) > 1$, $r\in R$ and $r,n\in\mathbb{N}$ and proved that circulant graphs $C_{16}(1,2,7)$ and $C_{16}(2,3,5)$ are Type-2 isomorphic w.r.t. $m$ = 2. He studied such graphs for $m$ = 2 in \cite{v13,v20} and with Wilson \cite{v24} obtained families of isomorphic circulant graphs of Type-2 w.r.t. $m$ = 3,5,7. In \cite{v2-2-arX}, Vilfred modified the definition of Type-2 isomorphism of circulant graphs $C_n(R)$ w.r.t. $m$ by considering $m > 1$ as a divisor of $\gcd(n, r)$ and $r\in R$ and in \cite{v2-1}, he further modified the definition to $r\in R$ and $m > 1$ and $m^3$ are divisors of $\gcd(n, r)$ and $n$, respectively. This modification is based on Theorems \ref{t10}  and \ref{t12} and it simplifies our calculations while finding isomorphic circulant graphs of Type-2. In this paper, we use the modified definition \ref{d6} to study Type-2 isomorphic circulant graphs of orders 16, 24 and 27 and show that the total number of pairs of Type-2 isomorphic circulant graphs of orders 16 and 24 are 8 and 32, respectively and the total number of triples of Type-2 isomorphic circulant graphs of order 27 are 12. To establish Type-2 isomorphism w.r.t. $m$ among $C_n(R)$ and $C_n(S)$, we follow remark \ref{r7} and to obtain more Type-2 isomorphic circulant graphs from known Type-2 isomorphic graph(s), we follow remark \ref{r8}. For basic definitions and results on isomorphic circulant graphs, refer \cite{v2-1}.

\section{Preliminaries }  

We present here a few definitions and results that are required in the subsequent sections.

\begin{theorem}{\rm\cite{v20}  \quad \label{t1}  In $C_n(R)$, the length of a cycle of period $r$ is $n/\gcd(n,r)$ and the number of disjoint periodic cycles of period $r$ is $\gcd(n,r)$, $r \in R$. \hfill $\Box$}
\end{theorem}

\begin{cor}{\rm \cite{v20} \quad \label{c2}  In $C_n(R)$, length of a cycle of period $r$ is $n$ if and only if $\gcd(n,r)$ = $1$, $r \in R$. \hfill $\Box$}
\end{cor}

\begin{theorem}{\rm \cite{v20} \quad \label{t3} If $C_n(R)$ and $C_n(S)$ are isomorphic, then there exists a bijection $f$ from $R \to S$ such that $\gcd(n, r)$ = $\gcd(n, f(r))$ for all $r\in R$.  \hfill $\Box$}
\end{theorem}

\begin{definition}{\rm\cite{ad67}} \quad \label{d4} For $R =$ $\{r_1$, $r_2$, $\ldots$, $r_k\}$ and $S$ = $\{s_1$, $s_2$, $\ldots$, $s_k\}$, circulant graphs $C_n(R)$ and $C_n(S)$ are {\it Adam's isomorphic} if there exists a positive integer $x$ $\ni$ $\gcd(n, x)$ = 1 and $S$ = $\{xr_1$, $xr_2$, $\ldots$, $xr_k\}_n^*$ where $<r_i>_n^*$, the {\it reflexive modular reduction} of a sequence $< r_i >$, is the sequence obtained by reducing each $r_i$ under modulo $n$ to yield $r_i'$ and then replacing all resulting terms $r_i'$ which are larger than $\frac{n}{2}$ by $n-r_i'.$  
\end{definition}

A circulant graph $C_n(R)$ is said to have {\em Cayley Isomorphism (CI) property} if whenever $C_n(S)$ is isomorphic to $C_n(R)$, they are Adam’s isomorphic \cite{v24}.

\begin{theorem} \cite{v24} \label{t5} Let $Ad_n(C_n(R))$ = $\{\varphi_{n,x}(C_n(R)) = C_n(xR): x\in\varphi_n \}$. Then, $C_n(S)\in Ad_n(C_n(R))$ if and only if $Ad_n(C_n(R))$ = $Ad_n(C_n(S))$ if and only if $C_n(R)\in Ad_n(C_n(S))$. \hfill $\Box$
\end{theorem}

In \cite{v2-1}, Vilfred modified the definition of Type-2 isomorphism of $C_n(R)$ w.r.t. $m$ as follows and hereafter we use the same definition. 

\begin{definition} \cite{v2-1} \quad  \label{d6} Let $V(K_n) = \{u_0,u_1,u_2,...,u_{n-1}\}$, $V(C_n(R))$ = $\{v_0, v_1, v_2, ... , v_{n-1}\}$, $|R| \geq 3$, $r\in R$ and $m > 1$ and $m^3$ be divisors of $\gcd(n, r)$ and $n$, respectively. Define 1-1 mapping $\theta_{n,m,t} :$ $V(C_n(R)) \rightarrow V(K_n)$ such that $\theta_{n,m,t}(v_x) = u_{x+jtm}$,  $\theta_{n,m,t}((v_x, v_{x+s}))$ = $(\theta_{n,m,t}(v_x),$ $\theta_{n,m,t}(v_{x+s}))$ under subscript arithmetic modulo $n$ and $\theta_{n,m,t}(C_n(R))$ = $C_n(\theta_{n,m,t}(R))$ where $\theta_{n,m,t}(R)$ in $C_n(\theta_{n,m,t}(R))$ is calculated under the reflexive modulo $n$, $\forall$ $x \in \mathbb{Z}_n$, $x = qm+j,$ $0 \leq j \leq m-1$, $s\in R$ and $0 \leq q,t \leq \frac{n}{m} -1$. And for a particular value of $t,$ if  $\theta_{n,m,t}(C_n(R))$ = $C_n(S)$ for some $S$  and  $S \neq yR$ for all $y\in \varphi_n$ under reflexive modulo $n,$ then $C_n(R)$ and $C_n(S)$ are called {\em isomorphic circulant graphs of Type-2 w.r.t. $m$.} 
	
When $C_n(R)$ and $C_n(S)$ are Type-2 isomorphic w.r.t. $m$, then we also say that $C_{kn}(kR)$  and $C_{kn}(kS)$ are Type-2 isomorphic w.r.t. $m$, $k\in\mathbb{N}$. Here, $k.C_n(T)$ = $C_{kn}(kT)$, $k\in\mathbb{N}$. 	 
\end{definition}

\begin{rem} \cite{v2-1}  \label{r7} Following steps are used to establish Type-2 isomorphism w.r.t. $m$ between circulant graphs $C_n(R)$ and $C_n(S)$. (i) $R$ $\neq$ $S$ and $|R| = |S| \geq 3$; (ii) $\exists$ $r\in R,S$ and $m > 1$ $\ni$ $m$ is a divisor of $\gcd(n, r)$, $m^3$ is a divisor of $n$ and for some $t$ $\ni$ $1 \leq t \leq \frac{n}{m} -1$, $\theta_{n,m,t}(C_n(R))$ = $C_n(S)$ and (iii) $S$ $\neq$ $xR$ for all $x\in\varphi_n$ under arithmetic reflexive modulo $n$. 

Thus, if $C_n(R)$ and $C_n(S)$ are Type-2 isomorphic circulant graphs w.r.t. $m$, then there exist $r\in R,S$, $m > 1$ and some $t$ $\ni$ $m$ is a divisor of $\gcd(n, r)$, $m^3$ is a divisor of $n$, $1 \leq t \leq \frac{n}{m} -1$, $\theta_{n,m,t}(C_n(R))$ = $C_n(S)$ and $S$ $\neq$ $xR$ for all $x\in\varphi_n$ under arithmetic reflexive modulo $n$.
\end{rem} 

\begin{rem}  \cite{v2-1} \label{r8} \quad The calculation on jump sizes $r_i$s which are integer multiples of $m$ need not be done under the transformation $\theta_{n,m,t}$, while searching for possible value(s) of $t$ for which the transformed graph $\theta_{n,m,t}(C_n(R))$ is circulant of the form $C_n(S)$ for some $S \subseteq [1, \frac{n}{2}]$, as there is no change in these $r_i$s where $r\in R$ and $m > 1$ and $m^3$ are divisors of $\gcd(n, r)$ and $n$, respectively. 

Thus, if $\theta_{n,m,t}(C_n(R))$ = $C_n(S)$ for some $S$ and thereby $C_n(R)$ $\cong$ $C_n(S)$, then $\theta_{n,m,t}(C_n(R \cup mT))$ = $C_n(S \cup mT)$ for any $T$ and thereby $C_n(R \cup mT)$ $\cong$ $C_n(S \cup mT)$.

Also, for a given $C_n(R)$, w.r.t. different values of $m$, we may get different Type-2 isomorphic circulant graphs.
\end{rem}

\begin{rem} \cite{v2-1} \label{r9} \quad {\rm For given $C_n(R)$ and $C_n(S)$ when either $\theta_{n,m,t}(C_n(R))$ = $C_n(S)$ for some $t$ or $C_n(xR)$ = $C_n(S)$ for some $x$, then $C_n(R)$ and $C_n(S)$ are isomorphic, $0 \leq t \leq \frac{n}{m} -1$ and $x\in\varphi_n$. }
\end{rem}

\begin{theorem}{\rm \cite{v20}}\quad \label{t10} {\rm For $n \geq 2$, $1 \leq 2s-1 \leq 2n-1$, $n \neq 2s-1$, $R$ = $\{2,2s-1, 4n-(2s-1)\}$ and $S$ = $\{ 2,$ $2n-(2s-1)$, $2n+2s-1 \}$, $\theta_{8n,2,n}(C_{8n}(R))$ = $C_{8n}(S)$ = $\theta_{8n,2,3n}(C_{8n}(R)),$ $\theta_{8n,2,n}(C_{8n}(S))$ = $C_{8n}(R)$ = $\theta_{8n,2,3n}(C_{8n}(S))$ and circulant graphs $C_{8n}(R)$ and $C_{8n}(S)$ are Type-2 isomorphic  w.r.t. $m$ = 2. When $n$ = $2s-1$, the two circulant graphs are the same. \hfill $\Box$}
\end{theorem}

\begin{theorem}{\rm \cite{v20}}\quad \label{t11} {\rm For $n \geq 2$, $1 \leq 2s-1 < 2s'-1 \leq [\frac{n}{2}]$, $0 \leq t \leq [\frac{n}{2}]$, $R$ = $\{2,2s-1, 2s'-1\}$ and $n,s,s'\in \mathbb{N}$, if $\theta_{n,2,t}(C_n(R))$ and $C_n(R)$ are  isomorphic circulant graphs of Type-2 w.r.t. $m$ = 2 for some $t$, then $n \equiv 0~(mod ~ 2^3)$, $2s-1+2s'-1$ = $\frac{n}{2}$, $2s-1 \neq \frac{n}{8}$, $t$ = $\frac{n}{8}$ or $\frac{3n}{8}$, $1 \leq 2s-1 \leq \frac{n}{4}$ and $n \geq 16$. In particular, $\theta_{8n,2,n}(C_{8n}(R))$ = $C_{8n}(S)$ = $\theta_{8n,2,3n}(C_{8n}(R))$ and $C_{8n}(R)$ and $C_{8n}(S)$ are Type-2 isomorphic w.r.t. $m$ = 2 when $R$ = $\{2, 2s-1, 4n-(2s-1)\}$, $S$ = $\{2, 2n-(2s-1), 2n+2s-1\}$, $n\geq 2$ and $n,s\in \mathbb{N}$. \hfill $\Box$}
\end{theorem}

\begin{theorem} \cite{v2-6} \label{t12} {\rm Let $p$ be an odd prime number, $1 \leq i \leq p$, $1 \leq x \leq p-1$, $y\in\mathbb{N}_0$, $0 \leq y \leq np-1$, $1 \leq x+yp \leq np^2-1$, $d^{np^3, x+yp}_i = (i-1)xpn+x+yp$,  $R^{np^3, x+yp}_i$ $=$ $\{p$, $d^{np^3, x+yp}_i$, $np^2-d^{np^3, x+yp}_i$, $np^2+d^{np^3, x+yp}_i$, $2np^2-d^{np^3, x+yp}_i$, $2np^2+$ $d^{np^3, x+yp}_i,$ $3np^2-d^{np^3, x+yp}_i$, $3np^2+d^{np^3, x+yp}_i$, . . . , $(p-1)np^2$ - $d^{np^3, x+yp}_i$, $(p-1)np^2+d^{np^3, x+yp}_i$, $np^3-d^{np^3, x+yp}_i$, $np^3-p\}$ and $i,j,n,x\in\mathbb{N}$. Then, for a given set of values of $n$, $p$, $x$ and $y$, $\theta_{np^3,p,jn} (C_{np^3}(R^{np^3, x+yp}_i))$ = $C_{np^3}(R^{np^3, x+yp}_{i+j})$ and the $p$ circulant graphs $C_{np^3}(R^{np^3, x+yp}_i)$ are isomorphic of Type-2 w.r.t.  $p$, $1 \leq i,j \leq p$ where $i+j$ in $R^{np^3, x+yp}_{i+j}$ is calculated under addition modulo $p$ and $C_{np^3}(R^{np^3, x+yp}_0)$ = $C_{np^3}(R^{np^3, x+yp}_p)$. \hfill $\Box$}
\end{theorem}

\section{Main results}

Given a circulant graph $C_n(R)$ having isomorphic circulant graphs of Type-2 w.r.t. $m$, remark \ref{r8} helps us to obtain more  isomorphic graphs which covers Type-2 w.r.t. $m$ as well as some Adam's isomorphic graphs of $C_n(R)$. In this section, we obtain all pairs of Type-2 isomorphic circulant graphs of orders 16 and 24 and all triples of Type-2 isomorphic circulant graphs of order 27. We show that the total number of pairs of Type-2 isomorphic circulant graphs of orders 16 and 24 are 8 and 32, respectively and the total number of triples of Type-2 isomorphic circulant graphs of order 27 are 12.

\subsection{Finding Type-2 isomorphic circulant graphs of order 16}

In this section, we find all Type-2 isomorphic circulant graphs of order 16. In \cite{v96} and \cite{v2-1}, it is shown that $C_{16}(1,2,7)$ and $C_{16}(2,3,5)$ are Type-2 isomorphic w.r.t. $m$ = 2. In the next problem, we find all pairs of Type-2 isomorphic circulant graphs of order 16 and show that the total number of pairs of Type-2 isomorphic circulant graphs of order 16 are 8 and each of these pairs is Type-2 isomorphic w.r.t. $m$ = 2.

\begin{prm}\quad \label{p1} {\rm Find all Type-2 isomorphic circulant graphs of  order 16. }
\end{prm}
\noindent
{\bf Solution.}\quad Here, we consider circulant graphs of the form $C_{n}(R)$ with $n$ = 16 = $2^4$ and so the possible values of $m > 1$ for existence of isomorphic circulant graphs of Type-2 w.r.t. $m$ $\ni$ $r\in R$, $m$ is a divisor of $\gcd(n, r)$ and $m^3$ is a divisor of $n$ is $m$ = 2. We start with finding isomorphic circulant graphs of Type-2 w.r.t. $m$ = 2. In \cite{v96}, it is proved that $\theta_{16,2,2}(C_{16}(1,2,7))$ = $C_{16}(2,3,5)$ and circulant graphs 

(1) $C_{16}(1,2,7)$ and $C_{16}(2,3,5)$ are isomorphic of Type-2 w.r.t. $m$ = 2.

Moreover, jump sizes 1,3,5,7 cover all the jump sizes other than multiples of 2 in these two circulant graphs and so each pair of Type-2 isomorphic circulant graphs of order 16 must contain jump sizes 1,3,5,7 other than multiples of 2.

Using remark \ref{r8} on this pair of isomorphic circulant graphs of Type-2 w.r.t. $m$ = 2 and by 2 = $\gcd(16, 2)$ = $\gcd(16, 6)$ = $\gcd(16, 2,4)$ = $\gcd(16, 2,6)$ = $\gcd(16, 2,8)$ = $\gcd(16, 4,6)$ = $\gcd(16, 6,8)$ = $\gcd(16, 2,4,6)$ = $\gcd(16, 2,4,8)$ = $\gcd(16, 2,6,8)$ = $\gcd(16, 4,6,8)$ = $\gcd(16, 2,4,6,8)$, $\gcd(16, 4)$ = 4 = $\gcd(16, 4,8)$ and $\gcd(16, 8)$ = 8, we consider the following pairs of isomorphic circulant graphs for possible isomorphic circulant graphs of Type-2 w.r.t. $m$ = 2. 		
\begin{enumerate}
\item [\rm (2)] $C_{16}(1,4,7)$, $C_{16}(3,4,5) = \theta_{16,2,2}(C_{16}(1,4,7))$; 

\item [\rm (3)] $C_{16}(1,6,7)$, $C_{16}(3,5,6) = \theta_{16,2,2}(C_{16}(1,6,7))$; 
	
\item [\rm (4)] $C_{16}(1,7,8)$, $C_{16}(3,5,8) = \theta_{16,2,2}(C_{16}(1,7,8))$; 

\item [\rm (5)] $C_{16}(1,2,4,7)$, $C_{16}(2,3,4,5) = \theta_{16,2,2}(C_{16}(1,2,4,7))$; 
	
\item [\rm (6)] $C_{16}(1,2,6,7)$, $C_{16}(2,3,5,6) = \theta_{16,2,2}(C_{16}(1,2,6,7))$; 
	
\item [\rm (7)] $C_{16}(1,2,7,8)$, $C_{16}(2,3,5,8) = \theta_{16,2,2}(C_{16}(1,2,7,8))$;
	
\item [\rm (8)] $C_{16}(1,4,6,7)$, $C_{16}(3,4,5,6) = \theta_{16,2,2}(C_{16}(1,4,6,7))$; 

\item [\rm (9)] $C_{16}(1,4,7,8)$, $C_{16}(3,4,5,8) = \theta_{16,2,2}(C_{16}(1,4,7,8))$; 
	
\item [\rm (10)] $C_{16}(1,6,7,8)$, $C_{16}(3,5,6,8) = \theta_{16,2,2}(C_{16}(1,6,7,8))$; 
	
\item [\rm (11)] $C_{16}(1,2,4,6,7)$, $C_{16}(2,3,4,5,6) = \theta_{16,2,2}(C_{16}(1,2,4,6,7))$; 
	
\item [\rm (12)] $C_{16}(1,2,4,7,8)$, $C_{16}(2,3,4,5,8) = \theta_{16,2,2}(C_{16}(1,2,4,7,8))$;
	
\item [\rm (13)] $C_{16}(1,2,6,7,8)$, $C_{16}(2,3,5,6,8) = \theta_{16,2,2}(C_{16}(1,2,6,7,8))$; 
	
\item [\rm (14)] $C_{16}(1,4,6,7,8)$, $C_{16}(3,4,5,6,8) = \theta_{16,2,2}(C_{16}(1,4,6,7,8))$ and
	
\item [\rm (15)] $C_{16}(1,2,4,6,7,8)$, $C_{16}(2,3,4,5,6,8) = \theta_{16,2,2}(C_{16}(1,2,4,6,7,8))$. 
\end{enumerate}
Among the 15 pairs of isomorphic circulant graphs, the following 7 pairs are Type-1 isomorphic since

\begin{enumerate}
\item [\rm (2)]	$C_{16}(1,4,7)$, $C_{16}(3,4,5)$ = $C_{16}(3(1,4,7))$. Also, $m$ = 2 is a divisor of $\gcd(16, 4)$ = 4; 
	
\item [\rm (4)] $C_{16}(1,7,8)$, $C_{16}(3,5,8)$ = $C_{16}(3(1,7,8))$. Here, $m$ = 2 is a divisor of $\gcd(16, 8)$ = 8; 	
	
\item [\rm (6)]	$C_{16}(1,2,6,7)$, $C_{16}(2,3,5,6)$ = $C_{16}(3(1,2,6,7))$. Here, $m$ = 2;
	
\item [\rm (9)]	$C_{16}(1,4,7,8)$, $C_{16}(3,4,5,8)$ = $C_{16}(3(1,4,7,8))$. Here, $m$ = 2 is a divisor of $\gcd(16, 4)$ = 4;

\item [\rm (11)]	$C_{16}(1,2,4,6,7)$, $C_{16}(2,3,4,5,6)$ = $C_{16}(3(1,2,4,6,7))$ and  $m$ = 2 is a divisor of $\gcd(16, 4)$ = 4; 
	
\item [\rm (13)]	$C_{16}(1,2,6,7,8)$, $C_{16}(2,3,5,6,8)$ = $C_{16}(3(1,2,6,7,8))$ and $m$ = 2 is a divisor of $\gcd(16, 6)$ = 2; 
	
\item [\rm (15)]	$C_{16}(1,2,4,6,7,8)$, $C_{16}(2,3,4,5,6,8)$ = $C_{16}(3(1,2,4,6,7,8))$ and $m$ = 2 is a divisor of $\gcd(16, 4)$.
\end{enumerate}
And all other pairs of isomorphic circulant graphs are of Type-2 w.r.t. $m$ = 2 since 
\begin{enumerate}
\item [\rm (3)]	$\theta_{16,2,2}(C_{16}(1,6,7))$ = $C_{16}(3,5,6)$. $\Rightarrow$ $C_{16}(1,6,7) \cong C_{16}(3,5,6)$. Also, 
	
$Ad_{16}(C_{16}(1,6,7))$ = $\{C_{16}(1,6,7), C_{16}(2,3,5) \}$. $\Rightarrow$ $C_{16}(3,5,6) \notin Ad_{16}(C_{16}(1,6,7))$.
	
$\Rightarrow$ $C_{16}(1,6,7)$ and $C_{16}(3,5,6)$ are Type-2 isomorphic w.r.t. $m$ = 2; 
	
\item [\rm (5)]	$\theta_{16,2,2}(C_{16}(1,2,4,7))$ = $C_{16}(2,3,4,5)$. $\Rightarrow$ $C_{16}(1,2,4,7) \cong C_{16}(2,3,4,5)$. Also, 
	
$Ad_{16}(C_{16}(1,2,4,7))$ = $\{C_{16}(1,2,4,7), C_{16}(3,4,5,6) \}$. $\Rightarrow$ $C_{16}(2,3,4,5) \notin Ad_{16}(C_{16}(1,2,4,7))$.
	
$\Rightarrow$ $C_{16}(1,2,4,7)$ and $C_{16}(2,3,4,5)$ are Type-2 isomorphic w.r.t. $m$ = 2;
	
\item [\rm (7)]	$\theta_{16,2,2}(C_{16}(1,2,7,8))$ = $C_{16}(2,3,5,8)$. $\Rightarrow$ $C_{16}(1,2,7,8) \cong C_{16}(2,3,5,8)$. Also, 
	
	$Ad_{16}(C_{16}(1,2,7,8))$ = $\{C_{16}(1,2,7,8), C_{16}(3,5,6,8) \}$. $\Rightarrow$ $C_{16}(2,3,5,8) \notin Ad_{16}(C_{16}(1,2,7,8))$.
	
$\Rightarrow$ $C_{16}(1,2,7,8)$ and $C_{16}(2,3,5,8)$ are Type-2 isomorphic w.r.t. $m$ = 2;
	
\item [\rm (8)]	$\theta_{16,2,2}(C_{16}(1,4,6,7))$ = $C_{16}(3,4,5,6)$. $\Rightarrow$ $C_{16}(1,4,6,7) \cong C_{16}(3,4,5,6)$. Also, 
	
	$Ad_{16}(C_{16}(1,4,6,7))$ = $\{C_{16}(1,4,6,7), C_{16}(2,3,4,5) \}$. $\Rightarrow$ $C_{16}(3,4,5,6) \notin Ad_{16}(C_{16}(1,4,6,7))$.
	
$\Rightarrow$ $C_{16}(1,4,6,7)$ and $C_{16}(3,4,5,6)$ are Type-2 isomorphic w.r.t. $m$ = 2;
	
\item [\rm (10)]	$\theta_{16,2,2}(C_{16}(1,6,7,8))$ = $C_{16}(3,5,6,8)$. $\Rightarrow$ $C_{16}(1,6,7,8) \cong C_{16}(3,5,6,8)$. Also, 
	
	$Ad_{16}(C_{16}(1,6,7,8))$ = $\{C_{16}(1,6,7,8), C_{16}(2,3,5,8) \}$. $\Rightarrow$ $C_{16}(3,5,6,8) \notin Ad_{16}(C_{16}(1,6,7,8))$.
	
$\Rightarrow$ $C_{16}(1,6,7,8)$ and $C_{16}(3,5,6,8)$ are Type-2 isomorphic w.r.t. $m$ = 2;
	
\item [\rm (12)]	$\theta_{16,2,2}(C_{16}(1,2,4,7,8))$ = $C_{16}(2,3,4,5,8)$. $\Rightarrow$ $C_{16}(1,2,4,7,8) \cong C_{16}(2,3,4,5,8)$. Also, 	
	
	$Ad_{16}(C_{16}(1,2,4,7,8))$ = $\{C_{16}(1,2,4,7,8), C_{16}(3,4,5,6,8) \}$. 
	
	$\Rightarrow$ $C_{16}(2,3,4,5,8) \notin Ad_{16}(C_{16}(1,2,4,7,8))$.
	
$\Rightarrow$ $C_{16}(1,2,4,7,8)$ and $C_{16}(2,3,4,5,8)$ are Type-2 isomorphic w.r.t. $m$ = 2 and
	
\item [\rm (14)] $\theta_{16,2,2}(C_{16}(1,4,6,7,8))$ = $C_{16}(3,4,5,6,8)$. $\Rightarrow$ $C_{16}(1,4,6,7,8) \cong C_{16}(3,4,5,6,8)$. Also, 
	
$Ad_{16}(C_{16}(1,4,6,7,8))$ = $\{C_{16}(1,4,6,7,8), C_{16}(2,3,4,5,8) \}$. 
	
$\Rightarrow$ $C_{16}(3,4,5,6,8) \notin Ad_{16}(C_{16}(1,4,6,7,8))$.
	
$\Rightarrow$ $C_{16}(1,4,6,7,8)$ and $C_{16}(3,4,5,6,8)$ are Type-2 isomorphic w.r.t. $m$ = 2.  
\end{enumerate}

\begin{figure}[ht]
	\centerline{\includegraphics[width=6in]{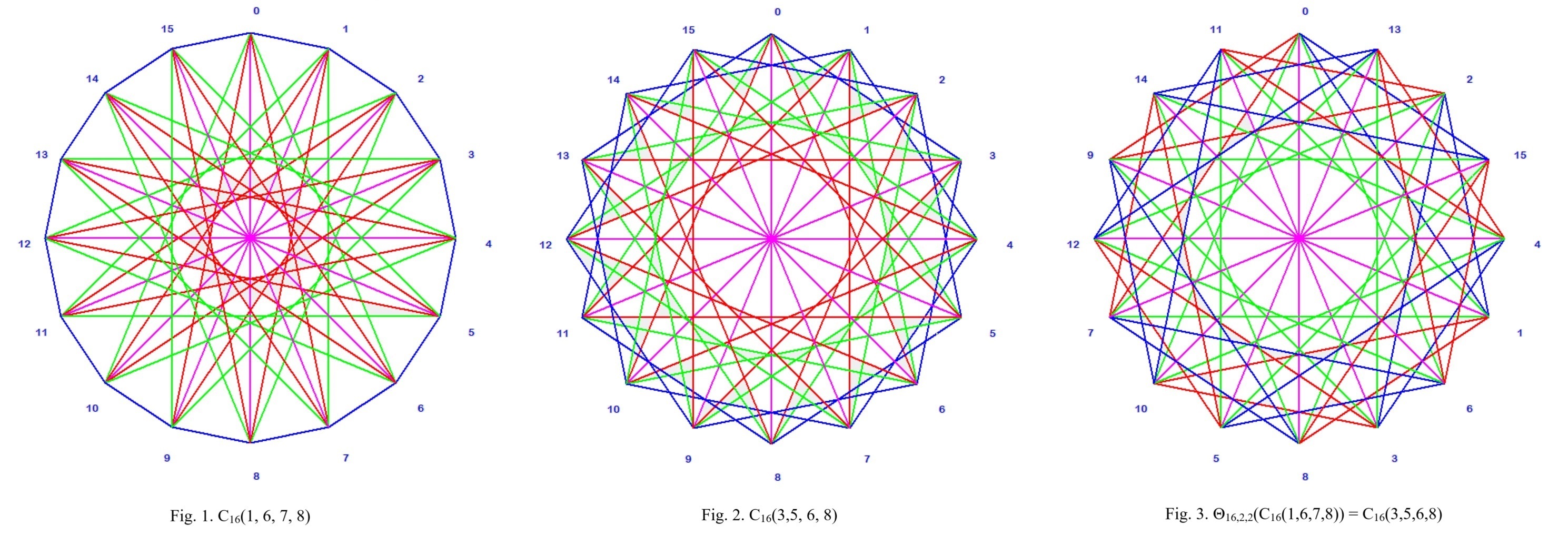}}
\end{figure}

Thus, there are 8 pairs of Type-2 isomorphic circulant graphs of order 16 and each pair is isomorphic of Type-2 w.r.t. $m$ = 2 only. Type-2 isomorphic circulant graphs $C_{16}(1, 6, 7, 8)$ and $C_{16}(3, 5, 6, 8)$ are given in Figures 1 and 2 and $\theta_{16,2,2}(C_{16}(1,6,7,8))$ = $C_{16}(3,5,6,8)$ in Figure 3. \hfill $\Box$   

\subsection{Finding Type-2 isomorphic circulant graphs of order 24}

In this section, we find all Type-2 isomorphic circulant graphs of order 24. In \cite{v2-1}, we proved that $C_{24}(1,2,8,11)$ and $C_{24}(2,5,7,8)$ are Type-2 isomorphic w.r.t. $m$ = 2 whereas $C_{24}(1,2,10,11)$ and $C_{24}(2,5,7,10)$ are Type-1 isomorphic. Here, we show, at first, that many pairs of isomorphic circulant graphs of the form $C_{24}(2r_1+1,x,2r_2+1)$ and $C_{24}(2s_1+1,x,2s_2+1)$ are of Type-1 where $x\in\mathbb{Z}_{24}$ $\ni$ $\gcd(24, x) > 1$, $r_1,r_2\in\mathbb{N}_0$ and $r_1 < r_2$ and then, we show that the pairs of circulant graphs $C_{24}(1,2,11)$, $C_{24}(2,5,7)$; and $C_{24}(1,10,11)$, $C_{24}(5,7,10)$ are Type-2 isomorphic w.r.t. $m$ = 2. And finally, we obtain all pairs of Type-2 isomorphic circulant graphs of order of 24 from the two pairs of Type-2 isomorphic circulant graphs. 

While searching for isomorphic circulant graphs of order 24 = $3\times 2^3$ and of Type-2 w.r.t. $m$, $m$ = 2 is the only possible value of $m$ for the existence of isomorphic circulant graphs of Type-2 w.r.t. $m$. Correspondingly, we consider graphs $C_{24}(2r_1+1,x,2r_2+1)$ and $\theta_{24,2,3t}(C_{24}(2r_1+1,x,2r_2+1))$ for $t$ = 1 to 3 and $x$ = 2, 10 for checking whether $\theta_{24,2,3t}(C_{24}(2r_1+1,x,2r_2+1))$ = $C_{24}(S)$ for some $S$ or not and also their type of isomorphism that exists among each pair of circulant graphs $\theta_{24,2,3t}(C_{24}(2r_1+1,x,2r_2+1)$ and $C_{24}(S)$, $r_1,r_2\in\mathbb{N}_0$, $r_1 < r_2$ and $x\in S$. In the next problem, we obtain some isomorphic circulant graphs of Type-1, each having CI-property.

\begin{prm}\quad \label{p2} {\rm Show that the following is true.
\begin{enumerate}
\item [\rm (a1)]  $C_{24}(1,2,3)$, $C_{24}(5,9,10)$, $C_{24}(3,7,10)$ and $C_{24}(2,9,11)$ are Type-1 isomorphic and each has CI-property.

\item [\rm (a2)]   $C_{24}(1,2,5)$, $C_{24}(1,5,10)$, $C_{24}(7,10,11)$ and $C_{24}(2,7,11)$ are Type-1 isomorphic and each  has CI-property.

\item [\rm (a3)] $C_{24}(1,2,7)$, $C_{24}(5,10,11)$, $C_{24}(1,7,10)$ and $C_{24}(2,5,11)$ are Type-1 isomorphic and each  has CI-property.

\item [\rm (a4)] $C_{24}(1,2,9)$, $C_{24}(3,5,10)$, $C_{24}(7,9,10)$ and $C_{24}(2,3,11)$ are Type-1 isomorphic and each  has CI-property.

\item [\rm (a5)]  $C_{24}(2,3,5)$, $C_{24}(1,9,10)$, $C_{24}(3,10,11)$ and $C_{24}(2,7,9)$ are Type-1 isomorphic and each  has CI-property.

\item [\rm (a6)]  $C_{24}(2,3,7)$, $C_{24}(9,10,11)$, $C_{24}(1,3,10)$ and $C_{24}(2,5,9)$ are Type-1 isomorphic and each  has CI-property.

\item [\rm (a7)]  $C_{24}(2,3,9)$ and $C_{24}(3,9,10)$ are Type-1 isomorphic and each  has CI-property.

\item [\rm (a8)]  $C_{24}(2,3,11)$, $C_{24}(7,9,10)$, $C_{24}(3,5,10)$ and $C_{24}(1,2,9)$ are Type-1 isomorphic and each  has CI-property. (Same as (a4)).

\item [\rm (a9)]  $C_{24}(2,5,9)$, $C_{24}(1,3,10)$, $C_{24}(9,10,11)$ and $C_{24}(2,3,7)$ are Type-1 isomorphic and each  has CI-property.  (Same as (a6)).

\item [\rm (a10)]  $C_{24}(2,5,11)$, $C_{24}(1,7,10)$, $C_{24}(5,10,11)$ and $C_{24}(1,2,7)$ are Type-1 isomorphic and each  has CI-property.  (Same as (a3)).

\item [\rm (a11)]  $C_{24}(2,7,9)$, $C_{24}(3,10,11)$, $C_{24}(1,9,10)$ and $C_{24}(2,3,5)$ are Type-1 isomorphic and each  has CI-property. (Same as (a5)).

\item [\rm (a12)]  $C_{24}(2,7,11)$, $C_{24}(7,10,11)$, $C_{24}(1,5,10)$ and $C_{24}(1,2,5)$ are Type-1 isomorphic and each  has CI-property. (Same as (a2)).

\item [\rm (a13)]  $C_{24}(2,9,11)$, $C_{24}(3,7,10)$, $C_{24}(5,9,10)$ and $C_{24}(1,2,3)$ are Type-1 isomorphic and each  has CI-property.  (Same as (a1)).

\item [\rm (b1)]  $C_{24}(1,3,10)$, $C_{24}(2,5,9)$, $C_{24}(2,3,7)$ and $C_{24}(9,10,11)$ are Type-1 isomorphic and each  has CI-property. (Same as (a6)).

\item [\rm (b2)]  $C_{24}(1,5,10)$, $C_{24}(1,2,5)$, $C_{24}(2,7,11)$ and $C_{24}(7,10,11)$ are Type-1 isomorphic and each  has CI-property.  (Same as (a2)).

\item [\rm (b3)]  $C_{24}(1,7,10)$, $C_{24}(2,5,11)$, $C_{24}(1,2,7)$ and $C_{24}(5,10,11)$ are Type-1 isomorphic and each  has CI-property.  (Same as (a3)).

\item [\rm (b4)]  $C_{24}(1,9,10)$, $C_{24}(2,3,5)$, $C_{24}(2,7,9)$ and $C_{24}(3,10,11)$ are Type-1 isomorphic and each  has CI-property. (Same as (a5)).

\item [\rm (b5)]  $C_{24}(3,5,10)$, $C_{24}(1,2,9)$, $C_{24}(2,3,11)$ and $C_{24}(7,9,10)$ are Type-1 isomorphic and each  has CI-property.  (Same as (a4)).

\item [\rm (b6)]  $C_{24}(3,7,10)$, $C_{24}(2,9,11)$, $C_{24}(1,2,3)$ and $C_{24}(5,9,10)$ are Type-1 isomorphic and each  has CI-property.  (Same as (a1)).

\item [\rm (b7)]  $C_{24}(3,9,10)$ and $C_{24}(2,3,9)$ are Type-1 isomorphic and each  has CI-property.  (Same as (a7)).

\item [\rm (b8)]  $C_{24}(3,10,11)$, $C_{24}(2,7,9)$, $C_{24}(2,3,5)$ and $C_{24}(1,9,10)$ are Type-1 isomorphic and each  has CI-property.  (Same as (a5)).

\item [\rm (b9)]  $C_{24}(5,9,10)$, $C_{24}(1,2,3)$, $C_{24}(2,9,11)$ and $C_{24}(3,7,10)$ are Type-1 isomorphic and each  has CI-property.  (Same as (a1)).

\item [\rm (b10)]  $C_{24}(5,10,11)$, $C_{24}(1,2,7)$, $C_{24}(2,5,11)$ and $C_{24}(1,7,10)$ are Type-1 isomorphic and each  has CI-property.  (Same as (a3)).

\item [\rm (b11)]  $C_{24}(7,9,10)$, $C_{24}(2,3,11)$, $C_{24}(1,2,9)$ and $C_{24}(3,5,10)$ are Type-1 isomorphic and each  has CI-property.  (Same as (a4)).

\item [\rm (b12)]  $C_{24}(7,10,11)$, $C_{24}(2,7,11)$, $C_{24}(1,2,5)$ and $C_{24}(1,5,10)$ are Type-1 isomorphic and each  has CI-property.  (Same as (a2)).

\item [\rm (b13)]  $C_{24}(9,10,11)$, $C_{24}(2,3,7)$, $C_{24}(2,5,9)$ and $C_{24}(1,3,10)$ are Type-1 isomorphic and each  has CI-property.  (Same as (a6)).

\item [\rm (c1)]  $C_{24}(1,3,4)$, $C_{24}(4,5,9)$, $C_{24}(3,4,7)$ and $C_{24}(4,9,11)$ are Type-1 isomorphic and each  has CI-property. 

\item [\rm (c2)]  $C_{24}(1,4,5)$ and $C_{24}(4,7,11)$ are Type-1 isomorphic and each  has CI-property.  

\item [\rm (c3)]  $C_{24}(1,4,7)$ and $C_{24}(4,5,11)$ are Type-1 isomorphic and each  has CI-property. 

\item [\rm (c4)]  $C_{24}(1,4,9)$, $C_{24}(3,4,5)$, $C_{24}(4,7,9)$ and $C_{24}(3,4,11)$ are Type-1 isomorphic and each  has CI-property.  

\item [\rm (c5)]  $C_{24}(1,4,11)$ and $C_{24}(4,5,7)$ are Type-1 isomorphic and each  has CI-property.

\item [\rm (c6)]  $C_{24}(3,4,5)$, $C_{24}(1,4,9)$, $C_{24}(3,4,11)$ and $C_{24}(4,7,9)$ are Type-1 isomorphic and each  has CI-property.    (Same as (c4)).

\item [\rm (c7)]  $C_{24}(3,4,7)$, $C_{24}(4,9,11)$, $C_{24}(1,3,4)$ and $C_{24}(4,5,9)$ are Type-1 isomorphic and each  has CI-property.    (Same as (c1)).

\item [\rm (c8)]  $C_{24}(3,4,9)$ has CI-property.  

\item [\rm (c9)]  $C_{24}(3,4,11)$, $C_{24}(4,7,9)$, $C_{24}(3,4,5)$ and $C_{24}(1,4,9)$ are Type-1 isomorphic and each  has CI-property.    (Same as (c4)).

\item [\rm (c10)]  $C_{24}(4,5,7)$ and $C_{24}(1,4,11)$ are Type-1 isomorphic and each  has CI-property.    (Same as (c5)).

\item [\rm (c11)]  $C_{24}(4,5,9)$, $C_{24}(1,3,4)$, $C_{24}(4,9,11)$ and $C_{24}(3,4,7)$ are Type-1 isomorphic and each  has CI-property.   (Same as (c1)).

\item [\rm (c12)]  $C_{24}(4,5,11)$ and $C_{24}(1,4,7)$ are Type-1 isomorphic and each  has CI-property.   (Same as (c3)). 

\item [\rm (c13)]  $C_{24}(4,7,9)$, $C_{24}(3,4,11)$, $C_{24}(1,4,9)$ and $C_{24}(3,4,5)$ are Type-1 isomorphic and each  has CI-property.    (Same as (c4)).

\item [\rm (c14)]  $C_{24}(4,7,11)$ and $C_{24}(1,4,5)$ are Type-1 isomorphic and each  has CI-property.   (Same as (c2)). 

\item [\rm (c15)]  $C_{24}(4,9,11)$, $C_{24}(3,4,7)$, $C_{24}(4,5,9)$ and $C_{24}(1,3,4)$  are Type-1 isomorphic and each  has CI-property.    (Same as (c1)). 

\item [\rm (d1)]  $C_{24}(1,3,6)$, $C_{24}(5,6,9)$, $C_{24}(3,6,7)$ and $C_{24}(6,9,11)$ are Type-1 isomorphic and each  has CI-property.

\item [\rm (d2)]  $C_{24}(1,5,6)$ and $C_{24}(6,7,11)$ are Type-1 isomorphic and each  has CI-property.

\item [\rm (d3)]  $C_{24}(1,6,7)$ and $C_{24}(5,6,11)$ are Type-1 isomorphic and each  has CI-property.

\item [\rm (d4)]  $C_{24}(1,6,9)$, $C_{24}(3,5,6)$, $C_{24}(6,7,9)$ and $C_{24}(3,6,11)$ are Type-1 isomorphic and each  has CI-property.

\item [\rm (d5)]  $C_{24}(1,6,11)$ and $C_{24}(5,6,7)$ are Type-1 isomorphic and each  has CI-property.

\item [\rm (d6)]  $C_{24}(3,5,6)$, $C_{24}(1,6,9)$, $C_{24}(3,6,11)$ and $C_{24}(6,7,9)$ are Type-1 isomorphic and each has CI-property.   (Same as (d4)).

\item [\rm (d7)]  $C_{24}(3,6,7)$, $C_{24}(6,9,11)$, $C_{24}(1,3,6)$ and $C_{24}(5,6,9)$ are Type-1 isomorphic and each  has CI-property.   (Same as (d1)).

\item [\rm (d8)]  $C_{24}(3,6,9)$ has CI-property.

\item [\rm (d9)]  $C_{24}(3,6,11)$, $C_{24}(6,7,9)$, $C_{24}(3,5,6)$ and $C_{24}(1,6,9)$ are Type-1 isomorphic and each  has CI-property.   (Same as (d4)).

\item [\rm (d10)]  $C_{24}(5,6,7)$ and $C_{24}(1,6,11)$ are Type-1 isomorphic and each  has CI-property.   (Same as (d5)).

\item [\rm (d11)]  $C_{24}(5,6,9)$, $C_{24}(1,3,6)$, $C_{24}(6,9,11)$ and $C_{24}(3,6,7)$ are Type-1 isomorphic and each  has CI-property.   (Same as (d1)).

\item [\rm (d12)]  $C_{24}(5,6,11)$ and $C_{24}(1,6,7)$ are Type-1 isomorphic and each  has CI-property.   (Same as (d3)).

\item [\rm (d13)]  $C_{24}(6,7,9)$, $C_{24}(3,6,11)$, $C_{24}(1,6,9)$ and $C_{24}(3,5,6)$ are Type-1 isomorphic and each  has CI-property.   (Same as (d4)).

\item [\rm (d14)]  $C_{24}(6,7,11)$ and $C_{24}(1,5,6)$ are Type-1 isomorphic and each  has CI-property.   (Same as (d2)).

\item [\rm (d15)]  $C_{24}(6,9,11)$, $C_{24}(3,6,7)$, $C_{24}(5,6,9)$ and $C_{24}(1,3,6)$ are Type-1 isomorphic and each  has CI-property.   (Same as (d1)).

\item [\rm (e1)]  $C_{24}(1,3,8)$,  $C_{24}(5,8,9)$, $C_{24}(3,7,8)$ and $C_{24}(8,9,11)$ are Type-1 isomorphic and each  has CI-property.

\item [\rm (e2)]  $C_{24}(1,5,8)$ and $C_{24}(7,8,11)$ are Type-1 isomorphic and each  has CI-property.

\item [\rm (e3)]  $C_{24}(1,7,8)$ and $C_{24}(5,8,11)$ are Type-1 isomorphic and each  has CI-property.

\item [\rm (e4)]  $C_{24}(1,8,9)$,  $C_{24}(3,5,8)$, $C_{24}(7,8,9)$ and $C_{24}(3,8,11)$ are Type-1 isomorphic and each  has CI-property.

\item [\rm (e5)]  $C_{24}(1,8,11)$ and $C_{24}(5,7,8)$ are Type-1 isomorphic and each  has CI-property.

\item [\rm (e6)]  $C_{24}(3,5,8)$,  $C_{24}(1,8,9)$, $C_{24}(3,8,11)$ and $C_{24}(7,8,9)$ are Type-1 isomorphic and each  has CI-property.   (Same as (e4)).

\item [\rm (e7)]  $C_{24}(3,7,8)$,  $C_{24}(8,9,11)$, $C_{24}(1,3,8)$ and $C_{24}(5,8,9)$ are Type-1 isomorphic and each  has CI-property.  (Same as (e1)).

\item [\rm (e8)]  $C_{24}(3,8,9)$ has CI-property.

\item [\rm (e9)]  $C_{24}(3,8,11)$,  $C_{24}(7,8,9)$, $C_{24}(3,5,8)$ and $C_{24}(1,8,9)$ are Type-1 isomorphic and each  has CI-property.  (Same as (e4)).

\item [\rm (e10)]  $C_{24}(5,7,8)$ and $C_{24}(1,8,11)$ are Type-1 isomorphic and each  has CI-property.  (Same as (e5)).

\item [\rm (e11)]  $C_{24}(5,8,9)$,  $C_{24}(1,3,8)$, $C_{24}(8,9,11)$ and $C_{24}(3,7,8)$ are Type-1 isomorphic and each  has CI-property.   (Same as (e1)).

\item [\rm (e12)]  $C_{24}(5,8,11)$ and $C_{24}(1,7,8)$ are Type-1 isomorphic and each  has CI-property.  (Same as (e3)).

\item [\rm (e13)]  $C_{24}(7,8,9)$,  $C_{24}(3,8,11)$, $C_{24}(1,8,9)$ and $C_{24}(3,5,8)$ are Type-1 isomorphic and each  has CI-property.  (Same as (e4)).

\item [\rm (e14)]  $C_{24}(7,8,11)$ and $C_{24}(1,5,8)$ are Type-1 isomorphic and each  has CI-property.  (Same as (e2)).

\item [\rm (e15)]  $C_{24}(8,9,11)$,  $C_{24}(3,7,8)$, $C_{24}(5,8,9)$ and $C_{24}(1,3,8)$ are Type-1 isomorphic and each  has CI-property.  (Same as (e1)).

\item [\rm (f1)]  $C_{24}(1,3,12)$, $C_{24}(5,9,12)$, $C_{24}(3,7,12)$ and $C_{24}(9,11,12)$ are Type-1 isomorphic and each  has CI-property.

\item [\rm (f2)]  $C_{24}(1,5,12)$ and $C_{24}(7,11,12)$ are Type-1 isomorphic and each  has CI-property.

\item [\rm (f3)]  $C_{24}(1,7,12)$ and $C_{24}(5,11,12)$ are Type-1 isomorphic and each  has CI-property.

\item [\rm (f4)]  $C_{24}(1,9,12)$, $C_{24}(3,5,12)$, $C_{24}(7,9,12)$ and $C_{24}(3,11,12)$ are Type-1 isomorphic and each  has CI-property.

\item [\rm (f5)]  $C_{24}(1,11,12)$ and $C_{24}(5,7,12)$ are Type-1 isomorphic and each  has CI-property.

\item [\rm (f6)]  $C_{24}(3,5,12)$, $C_{24}(1,9,12)$, $C_{24}(3,11,12)$ and $C_{24}(7,9,12)$ are Type-1 isomorphic and each  has CI-property.  (Same as (f4)).

\item [\rm (f7)]  $C_{24}(3,7,12)$, $C_{24}(9,11,12)$, $C_{24}(1,3,12)$ and $C_{24}(5,9,12)$ are Type-1 isomorphic and each  has CI-property.  (Same as (f1)).

\item [\rm (f8)]  $C_{24}(3,9,12)$  has CI-property.

\item [\rm (f9)]  $C_{24}(3,11,12)$, $C_{24}(7,9,12)$, $C_{24}(3,5,12)$ and $C_{24}(1,9,12)$ are Type-1 isomorphic and each  has CI-property.   (Same as (f4)).

\item [\rm (f10)]  $C_{24}(5,7,12)$ and $C_{24}(1,11,12)$ are Type-1 isomorphic and each  has CI-property.  (Same as (f5)).

\item [\rm (f11)]  $C_{24}(5,9,12)$, $C_{24}(1,3,12)$, $C_{24}(9,11,12)$ and $C_{24}(3,7,12)$ are Type-1 isomorphic and each  has CI-property.  (Same as (f1)).

\item [\rm (f12)]  $C_{24}(5,11,12)$ and $C_{24}(1,7,12)$ are Type-1 isomorphic and each  has CI-property.  (Same as (f3)).

\item [\rm (f13)]  $C_{24}(7,9,12)$, $C_{24}(3,11,12)$, $C_{24}(1,9,12)$ and $C_{24}(3,5,12)$ are Type-1 isomorphic and each  has CI-property.  (Same as (f4)).

\item [\rm (f14)]  $C_{24}(7,11,12)$ and $C_{24}(1,5,12)$ are Type-1 isomorphic and each  has CI-property.  (Same as (f2)).

\item [\rm (f15)]  $C_{24}(9,11,12)$, $C_{24}(3,7,12)$, $C_{24}(5,9,12)$ and $C_{24}(1,3,12)$ are Type-1 isomorphic and each  has CI-property.   (Same as (f1)).

\item [\rm (g1)]  $C_{24}(1,3,5)$, $C_{24}(1,5,9)$, $C_{24}(3,7,11)$ and $C_{24}(7,9,11)$ are Type-1 isomorphic and each  has CI-property.

\item [\rm (g2)]  $C_{24}(1,3,7)$ and $C_{24}(5,9,11)$ are Type-1 isomorphic and each  has CI-property.

\item [\rm (g3)]  $C_{24}(1,3,9)$, $C_{24}(3,5,9)$, $C_{24}(3,7,9)$ and $C_{24}(3,9,11)$ are Type-1 isomorphic and each  has CI-property.

\item [\rm (g4)]  $C_{24}(1,3,11)$, $C_{24}(5,7,9)$, $C_{24}(3,5,7)$ and $C_{24}(1,9,11)$ are Type-1 isomorphic and each  has CI-property.

\item [\rm (g5)]  $C_{24}(3,5,7)$, $C_{24}(1,9,11)$, $C_{24}(1,3,11)$ and $C_{24}(5,7,9)$ are Type-1 isomorphic and each  has CI-property.  (Same as (g4)).

\item [\rm (g6)]  $C_{24}(3,5,9)$, $C_{24}(1,3,9)$, $C_{24}(3,9,11)$ and $C_{24}(3,7,9)$ are Type-1 isomorphic and each  has CI-property.  (Same as (g3)).

\item [\rm (g7)]  $C_{24}(3,5,11)$ and $C_{24}(1,7,9)$ are Type-1 isomorphic and each  has CI-property. 

\item [\rm (g8)]  $C_{24}(3,7,9)$, $C_{24}(3,9,11)$, $C_{24}(1,3,9)$ and $C_{24}(3,5,9)$ are Type-1 isomorphic and each  has CI-property.  (Same as (g3)).

\item [\rm (g9)]  $C_{24}(3,7,11)$, $C_{24}(7,9,11)$, $C_{24}(1,3,5)$ and $C_{24}(1,5,9)$ are Type-1 isomorphic and each  has CI-property.  (Same as (g1)).

\item [\rm (g10)]  $C_{24}(3,9,11)$, $C_{24}(3,7,9)$, $C_{24}(3,5,9)$ and $C_{24}(1,3,9)$ are Type-1 isomorphic and each  has CI-property.  (Same as (g3)).

\item [\rm (h1)]  $C_{24}(1,3,9)$, $C_{24}(3,5,9)$, $C_{24}(3,7,9)$ and $C_{24}(3,9,11)$ are Type-1 isomorphic and each  has CI-property.  (Same as (g3)).

\item [\rm (h2)]  $C_{24}(1,5,9)$, $C_{24}(1,3,5)$, $C_{24}(7,9,11)$ and $C_{24}(3,7,11)$ are Type-1 isomorphic and each  has CI-property.   (Same as (g1)).

\item [\rm (h3)]  $C_{24}(1,7,9)$ and $C_{24}(3,5,11)$ are Type-1 isomorphic and each  has CI-property.

\item [\rm (h4)]  $C_{24}(1,9,11)$, $C_{24}(3,5,7)$, $C_{24}(5,7,9)$ and $C_{24}(1,3,11)$ are Type-1 isomorphic and each  has CI-property.  (Same as (g4)).

\item [\rm (h5)]  $C_{24}(3,5,9)$, $C_{24}(1,3,9)$, $C_{24}(3,9,11)$ and $C_{24}(3,7,9)$ are Type-1 isomorphic and each  has CI-property.  (Same as (g3)).

\item [\rm (h6)]  $C_{24}(3,7,9)$, $C_{24}(3,9,11)$, $C_{24}(1,3,9)$ and $C_{24}(3,5,9)$ are Type-1 isomorphic and each  has CI-property.  (Same as (g3)).

\item [\rm (h7)]  $C_{24}(3,9,11)$, $C_{24}(3,7,9)$, $C_{24}(3,5,9)$ and $C_{24}(1,3,9)$ are Type-1 isomorphic and each  has CI-property.  (Same as (g3)).

\item [\rm (h8)]  $C_{24}(5,7,9)$, $C_{24}(1,3,11)$, $C_{24}(1,9,11)$ and $C_{24}(3,5,7)$ are Type-1 isomorphic and each  has CI-property.  (Same as (g4)).

\item [\rm (h9)]  $C_{24}(5,9,11)$ and $C_{24}(1,3,7)$ are Type-1 isomorphic and each  has CI-property.  (Same as (g2)).

\item [\rm (h10)]  $C_{24}(7,9,11)$, $C_{24}(3,7,11)$, $C_{24}(1,5,9)$ and $C_{24}(1,3,5)$ are Type-1 isomorphic and each  has CI-property.   (Same as (g1)).
\end{enumerate} }
\end{prm}
 \noindent
{\bf Solution.}  We use remarks \ref{r7} and \ref{r9} to prove the results. Also, for a given circulant graph $C_n(R)$, if all $C_n(S)$ $\ni$ $C_n(S)$ = $\theta_{n,m,t}(C_n(R))$ for some $t$, $1 \leq t \leq \frac{n}{m}-1$ and $C_n(S)\in T1_n(C_n(R))$, then $C_n(R)$ has no isomorphic circulant graph of Type-2  w.r.t. $m$ where $r\in R,S$ and $m > 1$ and $m^3$ are divisors of $\gcd(n, r)$ and $n$, respectively. 

We consider graphs $C_{24}(2r_1+1,x,2r_2+1)$ and $\theta_{24,2,3t}(C_{24}(2r_1+1,x,2r_2+1))$ for $t$ = 1 to 3 and $x$ = 2,3,4,6,8, 9,10,12 for checking whether $\theta_{24,2,3t}(C_{24}(2r_1+1,x,2r_2+1))$ = $C_{24}(S)$ for some $S$ or not and also their type of isomorphism that exists among the pairs of circulant graphs $\theta_{24,2,3t}(C_{24}(2r_1+1,x,2r_2+1))$ and $C_{24}(S)$, $r_1,r_2\in\mathbb{N}_0$, $r_1 < r_2$ and $x\in S$. Solutions to all cases in this problem are similar, except cases $(g7)$ and $(h3)$ and to simplify our work, we present full calculations only to a few cases but important values related to all cases are presented in Tables 1 to 4. 
\begin{enumerate} 
\item [\rm (a1)] We have, 

$Ad_{24}(C_{24}(1,2,3))$ = $\{C_{24}(1,2,3)$, $C_{24}(5,9,10)$ = $C_{24}(5(1,2,3))$, 

\hfill $C_{24}(3,7,10)$ = $C_{24}(7(1,2,3))$, $C_{24}(2,9,11)$ = $C_{24}(11(1,2,3))\}$ 

\hfill = $Ad_{24}(C_{24}(5,9,10))$ = $Ad_{24}(C_{24}(3,7,10))$ = $Ad_{24}(C_{24}(2,9,11))$; 

$\theta_{24,2,6}(C_{24}(1,2,3))$ = $\theta_{24,2,6}(C_{24}(1,2,3, 21,22,23))$  = $C_{24}(\theta_{24,2,6}(1,2,3,  21,22,23))$

\hspace{2.85cm}  = $C_{24}(13,2,15, 9,22,11)$ = $C_{24}(2,9,11, 13,15,22)$ 

\hspace{2.85cm} = $C_{24}(2,9,11)\in Ad_{24}(C_{24}(1,2,3))$; 

$\theta_{24,2,6}(C_{24}(5,9,10))$ = $\theta_{24,2,6}(C_{24}(5,9,10, 14,15,19))$  = $C_{24}(\theta_{24,2,6}(5,9,10, 14,15,19))$

\hspace{3cm}  = $C_{24}(17,21,10, 14,3,7)$ = $C_{24}(3,7,10, 14,17,21)$ 

\hspace{3cm} = $C_{24}(3,7,10)\in Ad_{24}(C_{24}(5,9,10))$; 

$\theta_{24,2,6}(C_{24}(3,7,10))$ = $\theta_{24,2,6}(C_{24}(3,7,10, 14,17,21))$  = $C_{24}(\theta_{24,2,6}(3,7,10, 14,17,21))$

\hspace{3cm}  = $C_{24}(15,19,10, 14,5,9)$ = $C_{24}(5,9,10, 14,15,19)$ 

\hspace{3cm} = $C_{24}(5,9,10)\in Ad_{24}(C_{24}(3,7,10))$; 

$\theta_{24,2,6}(C_{24}(2,9,11))$ = $\theta_{24,2,6}(C_{24}(2,9,11, 13,15,22))$  = $C_{24}(\theta_{24,2,6}(2,9,11, 13,15,22))$

\hspace{3cm}  = $C_{24}(2,21,23, 1,3,22)$ = $C_{24}(1,2,3, 21,22,23)$ 

\hspace{3cm} = $C_{24}(1,2,3)\in Ad_{24}(C_{24}(2,9,11))$. 
\\
For $1 \leq i \leq 5$ and $i\in\mathbb{N}$,

$\theta_{24,2,6+i}(C_{24}(1,2,3))$ = $\theta_{24,2,6+i}(C_{24}(1,2,3, 21,22,23))$ 

\hspace{3.15cm} = $C_{24}(13+2i,2,15+2i, 9+2i,22,11+2i)$; 

$\theta_{24,2,6+i}(C_{24}(5,9,10))$ = $\theta_{24,2,6+i}(C_{24}(5,9,10, 14,15,19))$ 

\hspace{3.3cm} = $C_{24}(17+2i,21+2i,10, 14,3+2i,7+2i)$; 
 
$\theta_{24,2,6+i}(C_{24}(3,7,10))$ = $\theta_{24,2,6+i}(C_{24}(3,7,10, 14,17,21))$ 

\hspace{3.3cm} = $C_{24}(15+2i,19+2i,10, 14,5+2i,9+2i)$; 

$\theta_{24,2,6+i}(C_{24}(2,9,11))$ = $\theta_{24,2,6+i}(C_{24}(2,9,11, 13,15,22))$ 

\hspace{3.3cm} = $C_{24}(2,21+2i,23+2i, 1+2i,3+2i,22)$. 
\\
It is easy to see that for $i$ = 1 to 5,

 $\theta_{24,2,6+i}(C_{24}(1,2,3)), \theta_{24,2,6+i}(C_{24}(5,9,10)), \theta_{24,2,6+i}(C_{24}(3,7,10))$, $\theta_{24,2,6+i}(C_{24}(2,9,11))$ $\neq$ $C_{24}(S)$ for any $S$. Also, see Table 1 to Table 4.

$\Rightarrow$ $C_{24}(1,2,3)$, $C_{24}(5,9,10)$, $C_{24}(3,7,10)$ and $C_{24}(2,9,11)$ have no isomorphic circulant graph of Type-2 and these 4 circulant graphs are isomorphic of Type-1. 

This implies,  circulant graphs $C_{24}(1,2,3)$, $C_{24}(5,9,10)$, $C_{24}(3,7,10)$ and $C_{24}(2,9,11)$ are Type-1 isomorphic and each has CI-property. 

\begin{table}
	\caption{ Calculation of $\theta_{24,2,t}(\{1,2,3, 21,22,23\})$ for $t$ = 1 to 12 w.r.t. $C_{24}(1,2,3)$.}
\begin{center}
\scalebox{.9}{
 \begin{tabular}{||c||*{6}{c|}|c||c||}\hline \hline 
	Jump size $x$ &  1 & ~2~ & 3 & 21 & ~22~ & 23 & Pairwise Equidistant  \\ \cline{1-7} 
$\theta_{24,2,t}(x)$ & $\theta_{24,2,t}(1)$ & 2 & $\theta_{24,2,t}(3)$ & $\theta_{24,2,t}(21)$ & 22 & $\theta_{24,2,t}(23)$ & from $0$ or not in $\mathbb{Z}_{24}$ 	\\\hline \hline
	$\theta_{24,2,1}(x)$ & 3 & 2 & 5 & 23 & 22 & 1 & Not\\
	$\theta_{24,2,2}(x)$ & 5 & 2 & 7 & 1 & 22 & 3 & Not \\
	$\theta_{24,2,3}(x)$ & 7 & 2 & 9 & 3 & 22 & 5 & Not \\		
	$\theta_{24,2,4}(x)$ & 9 & 2 & 11 & 5 & 22 & 7 & Not \\
	$\theta_{24,2,5}(x)$ & 11 & 2 & 13 & 7 & 22 & 9 & Not \\   
	$\theta_{24,2,6}(x)$ & 13 & 2 & 15 & 9 & 22 & 11 & Yes (Type-1 with $11\in\varphi_{24}$) \\ 
	$\theta_{24,2,7}(x)$ & 15 & 2 & 17 & 11 & 22 & 13 & Not\\
	$\theta_{24,2,8}(x)$ & 17 & 2 & 19 & 13 & 22 & 15 & Not \\
	$\theta_{24,2,9}(x)$ & 19 & 2 & 21 & 15 & 22 & 17 & Not \\		
	$\theta_{24,2,10}(x)$ & 21 & 2 & 23 & 17 & 22 & 19 & Not \\
	$\theta_{24,2,11}(x)$ & 23 & 2 & 1 & 19 & 22 & 21 & Not \\   \hline \hline 
	$\theta_{24,2,12}(x)$ & 1 & 2 & 3 & 21 & 22 & 23 & Yes (same) \\  \hline \hline
\end{tabular}}
\end{center}
\end{table} 

\begin{table}
	\caption{ Calculation of $\theta_{24,2,t}(\{5,9,10, 14,15,19\})$ for $t$ = 1 to 11 w.r.t. $C_{24}(5,9,10)$.}
\begin{center}
\scalebox{.9}{
 \begin{tabular}{||c||*{6}{c|}|c||c||}\hline \hline 
	Jump size $x$ &  5 & 9 & ~10~ & ~14~ & 15 & 19 & Pairwise Equidistant  \\ \cline{1-7} 
$\theta_{24,2,t}(x)$ & $\theta_{24,2,t}(5)$ & $\theta_{24,2,t}(9)$ & 10 & 14 & $\theta_{24,2,t}(15)$ & $\theta_{24,2,t}(19)$ & from $0$ or not in $\mathbb{Z}_{24}$ 	\\\hline \hline
	$\theta_{24,2,1}(x)$ & 7 & 11 & 10 & 14 & 17 & 21 & Not\\
	$\theta_{24,2,2}(x)$ & 9 & 13 & 10 & 14 & 19 & 23 & Not \\
	$\theta_{24,2,3}(x)$ & 11 & 15 & 10 & 14 & 21 & 1 & Not \\		
	$\theta_{24,2,4}(x)$ & 13 & 17 & 10 & 14 & 23 & 3 & Not \\
	$\theta_{24,2,5}(x)$ & 15 & 19 & 10 & 14 & 1 & 5 & Not \\   
	$\theta_{24,2,6}(x)$ & 17 & 21 & 10 & 14 & 3 & 7 & Yes (Type-1 with $11\in\varphi_{24}$) \\ 
	$\theta_{24,2,7}(x)$ & 19 & 23 & 10 & 14 & 5 & 9 & Not\\
	$\theta_{24,2,8}(x)$ & 21 & 1 & 10 & 14 & 7 & 11 & Not \\
	$\theta_{24,2,9}(x)$ & 23 & 3 & 10 & 14 & 9 & 13 & Not \\		
	$\theta_{24,2,10}(x)$ & 1 & 5 & 10 & 14 & 11 & 15 & Not \\
	$\theta_{24,2,11}(x)$ & 3 & 7 & 10 & 14 & 13 & 17 & Not \\   \hline \hline 
\end{tabular}}
\end{center}
\end{table} 

\begin{table}
	\caption{ Calculation of $\theta_{24,2,t}(\{3,7,10, 14,17,21\})$ for $t$ = 1 to 11 w.r.t. $C_{24}(3,7,10)$.}
\begin{center}
\scalebox{.9}{
 \begin{tabular}{||c||*{6}{c|}|c||c||}\hline \hline 
	Jump size $x$ &  3 & 7 & ~10~ & ~14~ & 17 & 21 & Pairwise Equidistant  \\ \cline{1-7} 
$\theta_{24,2,t}(x)$ & $\theta_{24,2,t}(3)$ & $\theta_{24,2,t}(7)$ & 10 & 14 & $\theta_{24,2,t}(17)$ & $\theta_{24,2,t}(21)$ & from $0$ or not in $\mathbb{Z}_{24}$ 	\\\hline \hline
	$\theta_{24,2,1}(x)$ & 5 & 9 & 10 & 14 & 19 & 23 & Not\\
	$\theta_{24,2,2}(x)$ & 7 & 11 & 10 & 14 & 21 & 1 & Not \\
	$\theta_{24,2,3}(x)$ & 9 & 13 & 10 & 14 & 23 & 3 & Not \\		
	$\theta_{24,2,4}(x)$ & 11 & 15 & 10 & 14 & 1 & 5 & Not \\
	$\theta_{24,2,5}(x)$ & 13 & 17 & 10 & 14 & 3 & 7 & Not \\   
	$\theta_{24,2,6}(x)$ & 15 & 19 & 10 & 14 & 5 & 9 & Yes (Type-1 with $11\in\varphi_{24}$) \\ 
	$\theta_{24,2,7}(x)$ & 17 & 21 & 10 & 14 & 7 & 11 & Not\\
	$\theta_{24,2,8}(x)$ & 19 & 23 & 10 & 14 & 9 & 13 & Not \\
	$\theta_{24,2,9}(x)$ & 21 & 1 & 10 & 14 & 11 & 15 & Not \\		
	$\theta_{24,2,10}(x)$ & 23 & 3 & 10 & 14 & 13 & 17 & Not \\
	$\theta_{24,2,11}(x)$ & 1 & 5 & 10 & 14 & 15 & 19 & Not \\   \hline \hline 
\end{tabular}}
\end{center}
\end{table} 
\begin{table}
	\caption{ Calculation of $\theta_{24,2,t}(\{2,7,11, 13,17,22\})$ for $t$ = 1 to 11 w.r.t. $C_{24}(2,7,11)$.}
\begin{center}
\scalebox{.9}{
 \begin{tabular}{||c||*{6}{c|}|c||c||}\hline \hline 
	Jump size $x$ &  ~2~ & $\theta_{24,2,t}(7)$ & $\theta_{24,2,t}(11)$ & $\theta_{24,2,t}(13)$ & $\theta_{24,2,t}(17)$ & ~22~ & Pairwise Equidistant  \\ \cline{1-7} 
$\theta_{24,2,t}(x)$ & 2 & $\theta_{24,2,t}(7)$ & $\theta_{24,2,t}(11)$ & $\theta_{24,2,t}(13)$ & $\theta_{24,2,t}(17)$ & 22 & from $0$ or not in $\mathbb{Z}_{24}$ 	\\\hline \hline
	$\theta_{24,2,1}(x)$ & 2 & 9 & 13 & 15 & 19 & 22 & Not\\
	$\theta_{24,2,2}(x)$ & 2 & 11 & 15 & 17 & 21 & 22 & Not \\
	$\theta_{24,2,3}(x)$ & 2 & 13 & 17 & 19 & 23 & 22 & Not \\		
	$\theta_{24,2,4}(x)$ & 2 & 15 & 19 & 21 & 1 & 22 & Not \\
	$\theta_{24,2,5}(x)$ & 2 & 17 & 21 & 23 & 3 & 22 & Not \\   
	$\theta_{24,2,6}(x)$ & 2 & 19 & 23 & 1 & 5 & 22 & Yes (Type-1 with $11\in\varphi_{24}$) \\ 
	$\theta_{24,2,7}(x)$ & 2 & 21 & 1 & 3 & 7 & 22 & Not\\
	$\theta_{24,2,8}(x)$ & 2 & 23 & 3 & 5 & 9 & 22 & Not \\
	$\theta_{24,2,9}(x)$ & 2 & 1 & 5 & 7 & 11 & 22 & Not \\		
	$\theta_{24,2,10}(x)$ & 2 & 3 & 7 & 9 & 13 & 22 & Not \\
	$\theta_{24,2,11}(x)$ & 2 & 5 & 9 & 11 & 15 & 22 & Not \\   \hline \hline 
\end{tabular}}
\end{center}
\end{table} 

\item [\rm (a2)]  $Ad_{24}(C_{24}(1,2,5))$ = $\{C_{24}(1,2,5)$, $C_{24}(1,5,10)$ = $C_{24}(5(1,2,5))$, 

\hfill $C_{24}(7,10,11)$ = $C_{24}(7(1,2,5))$, $C_{24}(2,7,11)$ = $C_{24}(11(1,2,5))\}$ 

\hfill = $Ad_{24}(C_{24}(1,5,10))$ = $Ad_{24}(C_{24}(7,10,11))$ = $Ad_{24}(C_{24}(2,7,11))$; 

$\theta_{24,2,6}(C_{24}(1,2,5))$ = $\theta_{24,2,6}(C_{24}(1,2,5, 19,22,23))$  = $C_{24}(\theta_{24,2,6}(1,2,5, 19,22,23))$

\hspace{2.85cm}  = $C_{24}(13,2,17, 7,22,11)$ = $C_{24}(2,7,11, 13,17,22)$ 

\hspace{2.85cm} = $C_{24}(2,7,11)\in Ad_{24}(C_{24}(1,2,5))$; 

$\theta_{24,2,6}(C_{24}(1,5,10))$ = $\theta_{24,2,6}(C_{24}(1,5,10, 14,19,23))$  = $C_{24}(\theta_{24,2,6}(1,5,10, 14,19,23))$

\hspace{3cm}  = $C_{24}(13,17,10, 14,7,11)$ = $C_{24}(7,10,11, 13,14,17)$ 

\hspace{3cm} = $C_{24}(7,10,11)\in Ad_{24}(C_{24}(1,5,10))$; 

$\theta_{24,2,6}(C_{24}(7,10,11))$ = $\theta_{24,2,6}(C_{24}(7,10,11, 13,14,17))$  = $C_{24}(\theta_{24,2,6}(7,10,11, 13,14,17))$

\hspace{3.2cm}  = $C_{24}(19,10,23, 1,14,5)$ = $C_{24}(1,5,10, 14,19,23)$ 

\hspace{3.2cm} = $C_{24}(1,5,10)\in Ad_{24}(C_{24}(7,10,11))$; 

$\theta_{24,2,6}(C_{24}(2,7,11))$ = $\theta_{24,2,6}(C_{24}(2,7,11, 13,17,22))$  = $C_{24}(\theta_{24,2,6}(2,7,11, 13,17,22))$

\hspace{3cm}  = $C_{24}(2,19,23, 1,5,22)$ = $C_{24}(1,2,5, 19,22,23)$ 

\hspace{3cm} = $C_{24}(1,2,5)\in Ad_{24}(C_{24}(2,7,11))$. 
\\
For $1 \leq i \leq 5$ and $i\in\mathbb{N}$,

$\theta_{24,2,6+i}(C_{24}(1,2,5))$ = $\theta_{24,2,6+i}(C_{24}(1,2,5, 19,22,23))$ 

\hspace{3.1cm} = $C_{24}(13+2i,2,17+2i, 7+2i,22,11+2i)$;  

$\theta_{24,2,6+i}(C_{24}(1,5,10))$ = $\theta_{24,2,6+i}(C_{24}(1,5,10, 14,19,23))$ 

\hspace{3.35cm} = $C_{24}(13+2i,17+2i,10, 14,7+2i,11+2i)$;

$\theta_{24,2,6+i}(C_{24}(7,10,11))$ = $\theta_{24,2,6+i}(C_{24}(7,10,11, 13,14,17))$ 

\hspace{3.5cm} = $C_{24}(19+2i,10,23+2i, 1+2i,14,5+2i)$;

$\theta_{24,2,6+i}(C_{24}(2,7,11))$ = $\theta_{24,2,6+i}(C_{24}(2,7,11, 13,17,22))$ 

\hspace{3.35cm} = $C_{24}(2,19+2i,23+2i, 1+2i,5+2i,22)$.
\\
From the above relations, we get, for $i$ = 1 to 5,

 $\theta_{24,2,6+i}(C_{24}(1,2,5)), \theta_{24,2,6+i}(C_{24}(1,5,10)), \theta_{24,2,6+i}(C_{24}(7,10,11))$, $\theta_{24,2,6+i}(C_{24}(2,7,11))$ $\neq$ $C_{24}(S)$ for any $S$.

$\Rightarrow$ $C_{24}(1,2,5)$, $C_{24}(1,5,10)$, $C_{24}(7,10,11)$ and $C_{24}(2,7,11)$ have no isomorphic circulant graph of Type-2 and these 4 circulant graphs are Type-1 isomorphic. 

This implies,  circulant graphs $C_{24}(1,2,5)$, $C_{24}(1,5,10)$, $C_{24}(7,10,11)$ and $C_{24}(2,7,11)$ are Type-1 isomorphic and each has CI-property. 

\item [\rm (a3)] $Ad_{24}(C_{24}(1,2,7))$ = $\{C_{24}(1,2,7)$, $C_{24}(5,10,11)$ = $C_{24}(5(1,2,7))$, 

\hfill $C_{24}(1,7,10)$ = $C_{24}(7(1,2,7))$, $C_{24}(2,5,11)$ = $C_{24}(11(1,2,7))\}$

\hfill = $Ad_{24}(C_{24}(1,5,10))$ = $Ad_{24}(C_{24}(7,10,11))$ = $Ad_{24}(C_{24}(2,7,11))$;

$\theta_{24,2,6}(C_{24}(1,2,7))$ = $\theta_{24,2,6}(C_{24}(1,2,7, 17,22,23))$  = $C_{24}(\theta_{24,2,6}(1,2,7, 17,22,23))$

\hspace{2.85cm}  = $C_{24}(13,2,19, 5,22,11)$ = $C_{24}(2,5,11, 13,19,22)$ 

\hspace{2.85cm} = $C_{24}(2,5,11)\in Ad_{24}(C_{24}(1,2,7))$; 

$\theta_{24,2,6}(C_{24}(5,10,11))$ = $\theta_{24,2,6}(C_{24}(5,10,11, 13,14,19))$  = $C_{24}(\theta_{24,2,6}(5,10,11, 13,14,19))$

\hspace{3.2cm}  = $C_{24}(17,10,23, 1,14,7)$ = $C_{24}(1,7,10, 14,17,23)$ 

\hspace{3.2cm} = $C_{24}(1,7,10)\in Ad_{24}(C_{24}(5,10,11))$;

$\theta_{24,2,6}(C_{24}(1,7,10))$ = $\theta_{24,2,6}(C_{24}(1,7,10, 14,17,23))$  = $C_{24}(\theta_{24,2,6}(1,7,10, 14,17,23))$

\hspace{3cm}  = $C_{24}(13,19,10, 14,5,11)$ = $C_{24}(5,10,11, 13,14,19)$ 

\hspace{3cm} = $C_{24}(5,10,11)\in Ad_{24}(C_{24}(1,7,10))$;

$\theta_{24,2,6}(C_{24}(2,5,11))$ = $\theta_{24,2,6}(C_{24}(2,5,11, 13,19,22))$  = $C_{24}(\theta_{24,2,6}(2,5,11, 13,19,22))$

\hspace{3cm}  = $C_{24}(2,17,23, 1,7,22)$ = $C_{24}(1,2,7, 17,22,23)$ 

\hspace{3cm} = $C_{24}(1,2,7)\in Ad_{24}(C_{24}(2,5,11))$.
\\
For $1 \leq i \leq 5$ and $i\in\mathbb{N}$,

$\theta_{24,2,6+i}(C_{24}(1,2,7))$ = $\theta_{24,2,6+i}(C_{24}(1,2,7, 17,22,23))$ 

\hspace{3.15cm} = $C_{24}(13+2i,2,19+2i, 5+2i,22,11+2i)$; 

$\theta_{24,2,6+i}(C_{24}(5,10,11))$ = $\theta_{24,2,6+i}(C_{24}(5,10,11, 13,14,19))$ 

\hspace{3.5cm} = $C_{24}(17+2i,10,23+2i, 1+2i,14,7+2i)$;

$\theta_{24,2,6+i}(C_{24}(1,7,10))$ = $\theta_{24,2,6+i}(C_{24}(1,7,10, 14,17,23))$ 

\hspace{3.35cm} = $C_{24}(13+2i,19+2i,10, 14,5+2i,11+2i)$

$\theta_{24,2,6+i}(C_{24}(2,5,11))$ = $\theta_{24,2,6+i}(C_{24}(2,5,11, 13,19,22))$ 

\hspace{3.35cm} = $C_{24}(2,17+2i,23+2i, 1+2i,7+2i,22)$.

From the above relations, we get, for $i$ = 1 to 5,

 $\theta_{24,2,6+i}(C_{24}(1,2,7)), \theta_{24,2,6+i}(C_{24}(5,10,11)), \theta_{24,2,6+i}(C_{24}(1,7,10))$, $\theta_{24,2,6+i}(C_{24}(2,5,11))$ $\neq$ $C_{24}(S)$ for any $S$.

$\Rightarrow$ $C_{24}(1,2,7)$, $C_{24}(5,10,11)$, $C_{24}(1,7,10)$ and $C_{24}(2,5,11)$ have no isomorphic circulant graph of Type-2 w.r.t. $m$ = 2 and these 4 circulant graphs are Type-1 isomorphic. 

This implies,  circulant graphs $C_{24}(1,2,7)$, $C_{24}(5,10,11)$, $C_{24}(1,7,10)$ and $C_{24}(2,5,11)$ are Type-1 isomorphic and each has CI-property.  

\item [\rm (a7)]   $Ad_{24}(C_{24}(2,3,9))$ = $\{C_{24}(2,3,9)$, $C_{24}(3,9,10)$ = $C_{24}(5(2,3,9))\}$ = $Ad_{24}(C_{24}(3,9,10))$;

$\theta_{24,2,3}(C_{24}(2,3,9))$ = $\theta_{24,2,3}(C_{24}(2,3,9,  15,21,22))$ = $C_{24}(\theta_{24,2,3}(2,3,9,  15,21,22))$

\hspace{2.85cm}   = $C_{24}(2,9,15,  21,3,22)$ = $C_{24}(2,3,9, 15,21,22)$

\hspace{2.85cm} = $C_{24}(2,3,9)$ = $\theta_{24,2,6}(C_{24}(2,3,9))$ = $\theta_{24,2,9}(C_{24}(2,3,9))$;

$\theta_{24,2,3}(C_{24}(3,9,10))$ = $\theta_{24,2,3}(C_{24}(3,9,10,  14,15,21))$ = $C_{24}(\theta_{24,2,3}(3,9,10,  14,15,21))$ 

\hspace{3cm}  = $C_{24}(9,15,10,  14,21,3)$  = $C_{24}(3,9,10, 14,15,21)$ 

\hspace{3cm} = $C_{24}(3,9,10)$  = $\theta_{24,2,6}(C_{24}(3,9,10))$ = $\theta_{24,2,9}(C_{24}(3,9,10))$.
\\
For $t$ = 1 to 3,

$\theta_{24,2,3t+1}(C_{24}(2,3,9))$ = $\theta_{24,2,3t+1}(C_{24}(2,3,9, 15,21,22))$ 

\hspace{3.3cm} = $C_{24}(2,5+6t,11+6t, 17+6t, 23+6t, 22)$; 

$\theta_{24,2,3t+2}(C_{24}(2,3,9))$ = $\theta_{24,2,3t+2}(C_{24}(2,3,9, 15,21,22))$ 

\hspace{3.3cm} = $C_{24}(2,7+6t,13+6t, 19+6t, 1+6t, 22)$;  

$\theta_{24,2,3t+1}(C_{24}(3,9,10))$ = $\theta_{24,2,3t+1}(C_{24}(3,9,10, 14,15,21))$ 

\hspace{3.5cm} = $C_{24}(5+6t,11+6t,10, 14,17+6t,23+6t)$; and 

$\theta_{24,2,3t+2}(C_{24}(3,9,10))$ = $\theta_{24,2,3t+2}(C_{24}(3,9,10, 14,15,21))$ 

\hspace{3.47cm} = $C_{24}(7+6t,13+6t,10, 14,19+6t,1+6t)$.
\\
$\Rightarrow$ $\theta_{24,2,3t+1}(C_{24}(2,3,9)), \theta_{24,2,3t+2}(C_{24}(2,3,9), \theta_{24,2,3t+1}(C_{24}(3,9,10))$, $\theta_{24,2,3t+2}(C_{24}(3,9,10)$ $\neq$ $C_{24}(S)$ for any $S$ and $t$ = 1,2,3 (i.e., when $1 \leq 3t+1, 3t+2 \leq 11$ and $t\in\mathbb{N}$).

$\Rightarrow$ $C_{24}(2,3,9)$ and $C_{24}(3,9,10)$ have no Type-2 isomorphic circulant graph. 

$\Rightarrow$  Circulant graphs $C_{24}(2,3,9)$ and $C_{24}(3,9,10)$ are Type-1 isomorphic and each has CI-property.

\item [\rm (c5)]   $Ad_{24}(C_{24}(1,4,11))$ = $\{C_{24}(1,4,11), C_{24}(4,5,7) = C_{24}(5(1,4,11))\}$ = $Ad_{24}(C_{24}(4,5,7))$;

$\theta_{24,2,3}(C_{24}(1,4,11))$ = $\theta_{24,2,3}(C_{24}(1,4,11, 13,20,23))$ = $C_{24}(\theta_{24,2,3}(1,4,11, 13,20,23))$ 

\hspace{3cm}  = $C_{24}(7,4,17, 19,20,5)$ = $C_{24}(4,5,7, 17,19,20)$

\hspace{3cm}  = $C_{24}(4,5,7)$ = $\theta_{24,2,9}(C_{24}(1,4,11))\in Ad_{24}(C_{24}(1,4,11))$,   

$\theta_{24,2,6}(C_{24}(1,4,11))$ = $C_{24}(1,4,11)$;   

$\theta_{24,2,3}(C_{24}(4,5,7))$ = $\theta_{24,2,3}(C_{24}(4,5,7, 17,19,20))$ = $C_{24}(\theta_{24,2,3}(4,5,7, 17,19,20))$ 

\hspace{2.85cm}  = $C_{24}(4,11,13, 23,1,20)$ = $C_{24}(1,4,11, 13,20,23)$

\hspace{2.85cm}  = $C_{24}(1,4,11)$ = $\theta_{24,2,9}(C_{24}(4,5,7))\in Ad_{24}(C_{24}(4,5,7))$,   

$\theta_{24,2,6}(C_{24}(4,5,7))$ = $C_{24}(4,5,7)$.   
\\
For $t$ = 1 to 3  (i.e., when $1 \leq 3t+1, 3t+2 \leq 11$ and $t\in\mathbb{N}$),

$\theta_{24,2,3t+1}(C_{24}(1,4,11))$ = $\theta_{24,2,3t+1}(C_{24}(1,4,11, 13,20,23))$ 

\hspace{3.5cm} = $C_{24}(3+6t,4,13+6t, 15+6t, 20,1+6t)$,  

$\theta_{24,2,3t+2}(C_{24}(1,4,11))$ = $\theta_{24,2,3t+2}(C_{24}(1,4,11, 13,20,23))$ 

\hspace{3.5cm} = $C_{24}(5+6t,4,15+6t, 17+6t, 20,3+6t)$; and  

$\theta_{24,2,3t+1}(C_{24}(4,5,7))$ = $\theta_{24,2,3t+1}(C_{24}(4,5,7, 17,19,20))$ 

\hspace{3.3cm} = $C_{24}(4,7+6t,9+6t, 19+6t,21+6t,20)$,  

$\theta_{24,2,3t+2}(C_{24}(4,5,7))$ = $\theta_{24,2,3t+2}(C_{24}(4,5,7, 17,19,20))$ 

\hspace{3.3cm} = $C_{24}(4,9+6t,11+6t, 21+6t,23+6t,20)$
\\
$\Rightarrow$ $\theta_{24,2,3t+1}(C_{24}(1,4,11)), \theta_{24,2,3t+2}(C_{24}(1,4,11), \theta_{24,2,3t+1}(C_{24}(4,5,7))$, $\theta_{24,2,3t+2}(C_{24}(4,5,7)$ $\neq$ $C_{24}(S)$ for any $S$ and $t$ = 1,2,3.

$\Rightarrow$ $C_{24}(1,4,11)$ and $C_{24}(4,5,7)$ have no Type-2 isomorphic circulant graph.  

$\Rightarrow$ Circulant graphs $C_{24}(1,4,11)$ and $C_{24}(4,5,7)$ are Type-1 isomorphic and each has CI-property.

\item [\rm (d5)]  $Ad_{24}(C_{24}(1,6,11))$ = $\{C_{24}(1,6,11), C_{24}(5,6,7) = C_{24}(5(1,6,11))\}$ = $Ad_{24}(C_{24}(5,6,7))$;

$\theta_{24,2,3}(C_{24}(1,6,11))$ = $\theta_{24,2,3}(C_{24}(1,6,11, 13,18,23))$ = $C_{24}(\theta_{24,2,3}(1,6,11, 13,18,23))$ 

\hspace{3cm}  = $C_{24}(7,6,17, 19,18,5)$ = $C_{24}(5,6,7, 17,18,19)$

\hspace{3cm}  = $C_{24}(5,6,7)$ = $\theta_{24,2,9}(C_{24}(1,6,11))\in Ad_{24}(C_{24}(1,6,11))$,   

$\theta_{24,2,6}(C_{24}(1,6,11))$ = $C_{24}(1,6,11)$;   

$\theta_{24,2,3}(C_{24}(5,6,7))$ = $\theta_{24,2,3}(C_{24}(5,6,7, 17,18,19))$ = $C_{24}(\theta_{24,2,3}(5,6,7, 17,18,19))$ 

\hspace{2.85cm}  = $C_{24}(11,6,13, 23,18,1)$ = $C_{24}(1,6,11, 13,18,23)$

\hspace{2.85cm}  = $C_{24}(1,6,11)$ = $\theta_{24,2,9}(C_{24}(5,6,7))\in Ad_{24}(C_{24}(5,6,7))$,   

$\theta_{24,2,6}(C_{24}(5,6,7))$ = $C_{24}(5,6,7)$.   
\\
For $t$ = 1 to 3  (i.e., when $1 \leq 3t+1, 3t+2 \leq 11$ and $t\in\mathbb{N}$),

$\theta_{24,2,3t+1}(C_{24}(1,6,11))$ = $\theta_{24,2,3t+1}(C_{24}(1,6,11, 13,18,23))$ 

\hspace{3.5cm} = $C_{24}(3+6t,6,13+6t, 15+6t, 18,1+6t)$,  

$\theta_{24,2,3t+2}(C_{24}(1,6,11))$ = $\theta_{24,2,3t+2}(C_{24}(1,6,11, 13,18,23))$ 

\hspace{3.5cm} = $C_{24}(5+6t,6,15+6t, 17+6t, 18,3+6t)$; and  

$\theta_{24,2,3t+1}(C_{24}(5,6,7))$ = $\theta_{24,2,3t+1}(C_{24}(5,6,7, 17,18,19))$ 

\hspace{3.5cm} = $C_{24}(7+6t,6,9+6t, 19+6t,18,21+6t)$,  

$\theta_{24,2,3t+2}(C_{24}(5,6,7))$ = $\theta_{24,2,3t+2}(C_{24}(5,6,7, 17,18,19))$ 

\hspace{3.5cm} = $C_{24}(9+6t,6,11+6t, 21+6t,18,23+6t)$
\\
$\Rightarrow$ $\theta_{24,2,3t+1}(C_{24}(1,6,11)), \theta_{24,2,3t+2}(C_{24}(1,6,11), \theta_{24,2,3t+1}(C_{24}(5,6,7))$, $\theta_{24,2,3t+2}(C_{24}(5,6,7)$ $\neq$ $C_{24}(S)$ for any $S$ and $t$ = 1,2,3.

$\Rightarrow$ $C_{24}(1,6,11)$ and $C_{24}(5,6,7)$ have no Type-2 isomorphic circulant graph.  

$\Rightarrow$ Circulant graphs $C_{24}(1,6,11)$ and $C_{24}(5,6,7)$ are Type-1 isomorphic and each has CI-property.

\item [\rm (e5)]  $Ad_{24}(C_{24}(1,8,11))$ = $\{C_{24}(1,8,11)$, $C_{24}(5,7,8)$ = $C_{24}(5(1,8,11))\}$ =$Ad_{24}(C_{24}(5,7,8))$;

$\theta_{24,2,3}(C_{24}(1,8,11))$ = $\theta_{24,2,3}(C_{24}(1,8,11, 13,16,23))$ = $C_{24}(\theta_{24,2,3}(1,8,11, 13,16,23))$ 

\hspace{3cm}  = $C_{24}(7,8,17, 19,16,5)$ = $C_{24}(5,7,8, 16,17,19)$

\hspace{3cm}  = $C_{24}(5,7,8)$ = $\theta_{24,2,9}(C_{24}(1,8,11))\in Ad_{24}(C_{24}(1,8,11))$,   

$\theta_{24,2,6}(C_{24}(1,8,11))$ = $C_{24}(1,8,11)$;   

$\theta_{24,2,3}(C_{24}(5,7,8))$ = $\theta_{24,2,3}(C_{24}(5,7,8, 16,17,19))$ = $C_{24}(\theta_{24,2,3}(5,7,8, 16,17,19))$ 

\hspace{2.8cm}  = $C_{24}(11,13,8, 16,23,1)$ = $C_{24}(1,8,11, 13,16,23)$

\hspace{2.8cm}  = $C_{24}(1,8,11)$ = $\theta_{24,2,9}(C_{24}(5,7,8))\in Ad_{24}(C_{24}(5,7,8))$,   

$\theta_{24,2,6}(C_{24}(5,7,8))$ = $C_{24}(5,7,8)$.   
\\
For $t$ = 1 to 3  (i.e., when $1 \leq 3t+1, 3t+2 \leq 11$ and $t\in\mathbb{N}$),

$\theta_{24,2,3t+1}(C_{24}(1,8,11))$ = $\theta_{24,2,3t+1}(C_{24}(1,8,11, 13,16,23))$ 

\hspace{3.5cm} = $C_{24}(3+6t,8,13+6t, 15+6t, 16,1+6t)$,  

$\theta_{24,2,3t+2}(C_{24}(1,8,11))$ = $\theta_{24,2,3t+2}(C_{24}(1,8,11, 13,16,23))$ 

\hspace{3.5cm} = $C_{24}(5+6t,8,15+6t, 17+6t, 16,3+6t)$; and  

$\theta_{24,2,3t+1}(C_{24}(5,7,8))$ = $\theta_{24,2,3t+1}(C_{24}(5,7,8, 16,17,19))$ 

\hspace{3.3cm} = $C_{24}(7+6t,9+6t,8, 16,19+6t,21+6t)$,  

$\theta_{24,2,3t+2}(C_{24}(5,7,8))$ = $\theta_{24,2,3t+2}(C_{24}(5,7,8, 16,17,19))$ 

\hspace{3.3cm} = $C_{24}(9+6t,11+6t,8, 16,21+6t,23+6t)$.
\\
$\Rightarrow$ $\theta_{24,2,3t+1}(C_{24}(1,8,11)), \theta_{24,2,3t+2}(C_{24}(1,8,11), \theta_{24,2,3t+1}(C_{24}(5,7,8))$, $\theta_{24,2,3t+2}(C_{24}(5,7,8)$ $\neq$ $C_{24}(S)$ for any $S$ and $t$ = 1,2,3.

$\Rightarrow$ $C_{24}(1,8,11)$ and $C_{24}(5,7,8)$ have no Type-2 isomorphic circulant graph.  

$\Rightarrow$ Circulant graphs $C_{24}(1,8,11)$ and $C_{24}(5,7,8)$ are Type-1 isomorphic and each has CI-property.

\item [\rm (f5)]  $Ad_{24}(C_{24}(1,11,12))$ = $\{C_{24}(1,11,12), C_{24}(5,7,12) = C_{24}(5(1,11,12))\}$ =$Ad_{24}(C_{24}(5,7,12))$;

$\theta_{24,2,3}(C_{24}(1,11,12))$ = $\theta_{24,2,3}(C_{24}(1,11, 12, 13,23))$ = $C_{24}(\theta_{24,2,3}(1,11, 12, 13,23))$ 

\hspace{3.2cm}  = $C_{24}(7,17, 12, 19,5)$ = $C_{24}(5,7, 12, 17,19)$

\hspace{3.2cm}  = $C_{24}(5,7,12)$ = $\theta_{24,2,9}(C_{24}(1,11,12))\in Ad_{24}(C_{24}(1,11,12))$,   

$\theta_{24,2,6}(C_{24}(1,11,12))$ = $C_{24}(1,11,12)$;   

$\theta_{24,2,3}(C_{24}(5,7,12))$ = $\theta_{24,2,3}(C_{24}(5,7, 12, 17,19))$ = $C_{24}(\theta_{24,2,3}(5,7, 12, 17,19))$ 

\hspace{3cm}  = $C_{24}(11,13, 12, 23,1)$ = $C_{24}(1,11, 12, 13,23)$

\hspace{3cm}  = $C_{24}(1,11,12)$ = $\theta_{24,2,9}(C_{24}(5,7,12))\in Ad_{24}(C_{24}(5,7,12))$,   

$\theta_{24,2,6}(C_{24}(5,7,12))$ = $C_{24}(5,7,12)$.   
\\
For $t$ = 1 to 3  (i.e., when $1 \leq 3t+1, 3t+2 \leq 11$ and $t\in\mathbb{N}$),

$\theta_{24,2,3t+1}(C_{24}(1,11,12))$ = $\theta_{24,2,3t+1}(C_{24}(1,11, 12, 13,23))$ 

\hspace{3.5cm} = $C_{24}(3+6t,13+6t, 12, 15+6t,1+6t)$,  

$\theta_{24,2,3t+2}(C_{24}(1,11,12))$ = $\theta_{24,2,3t+2}(C_{24}(1,11, 12, 13,23))$ 

\hspace{3.5cm} = $C_{24}(5+6t,15+6t, 12, 17+6t,3+6t)$; and  

$\theta_{24,2,3t+1}(C_{24}(5,7,12))$ = $\theta_{24,2,3t+1}(C_{24}(5,7, 12, 17,19))$ 

\hspace{3.5cm} = $C_{24}(7+6t,9+6t, 12, 19+6t,21+6t)$,  

$\theta_{24,2,3t+2}(C_{24}(5,7,12))$ = $\theta_{24,2,3t+2}(C_{24}(5,7, 12, 17,19))$ 

\hspace{3.5cm} = $C_{24}(9+6t,11+6t, 12, 21+6t,23+6t)$.

This implies that $\theta_{24,2,3t+1}(C_{24}(1,11,12)), \theta_{24,2,3t+2}(C_{24}(1,11,12), \theta_{24,2,3t+1}(C_{24}(5,7,12))$, 

$\theta_{24,2,3t+2}(C_{24}(5,7,12)$ $\neq$ $C_{24}(S)$ for any $S$ and $t$ = 1,2,3.

$\Rightarrow$ $C_{24}(1,11,12)$ and $C_{24}(5,7,12)$ have no Type-2 isomorphic circulant graph.  

$\Rightarrow$ Circulant graphs $C_{24}(1,11,12)$ and $C_{24}(5,7,12)$ are Type-1 isomorphic and each has CI-property.  

\item [\rm (g7)]  $Ad_{24}(C_{24}(3,5,11))$ = $\{C_{24}(3,5,11), C_{24}(1,7,9) = C_{24}(5(3,5,11))\}$ =$Ad_{24}(C_{24}(1,7,9))$;

$\theta_{24,2,2}(C_{24}(3,5,11))$ = $\theta_{24,2,2}(C_{24}(3,5,11, 13,19,21))$ = $C_{24}(\theta_{24,2,2}(3,5,11, 13,19,21))$ 

\hspace{3cm}  = $C_{24}(7,9,15, 17,23,1)$ = $C_{24}(1,7,9, 15,17,23)$

\hspace{3cm}  = $C_{24}(1,7,9)$ = $\theta_{24,2,6}(C_{24}(3,5,11))$ = $\theta_{24,2,10}(C_{24}(3,5,11))$,   

$\theta_{24,2,4}(C_{24}(3,5,11))$ = $C_{24}(3,5,11)$ = $\theta_{24,2,8}(C_{24}(3,5,11))$;   

$\theta_{24,2,2}(C_{24}(1,7,9))$ = $\theta_{24,2,2}(C_{24}(1,7,9, 15,17,23))$ = $C_{24}(\theta_{24,2,2}(1,7,9, 15,17,23))$ 

\hspace{2.85cm}  = $C_{24}(5,11,13, 19,21,3)$ = $C_{24}(3,5,11, 13,19,21)$

\hspace{2.85cm}  = $C_{24}(3,5,11)$ = $\theta_{24,2,6}(C_{24}(1,7,9))$ = $\theta_{24,2,10}(C_{24}(1,7,9))$,   

$\theta_{24,2,4}(C_{24}(1,7,9))$ = $C_{24}(1,7,9)$ = $\theta_{24,2,8}(C_{24}(1,7,9))$,   

$\theta_{24,2,1}(3,5,11, 13,19,21)$ = $\{5,7,13, 15,21,23\}$,  

$\theta_{24,2,3}(3,5,11, 13,19,21)$ = $\{9,11,17, 19,1,3\}$ = $\{1,3,9, 11,17,19\}$,  

$\theta_{24,2,5}(3,5,11, 13,19,21)$ = $\{13,15,21, 23,5,7\}$ = $\{5,7,13, 15,21,23\}$,  

$\theta_{24,2,7}(3,5,11, 13,19,21)$ = $\{17,19,1, 3,9,11\}$ = $\{1,3,9, 11,17,19\}$,  

$\theta_{24,2,9}(3,5,11, 13,19,21)$ = $\{21,23,5, 7,13,15\}$ = $\{5,7,13, 15,21,23\}$,

$\theta_{24,2,11}(3,5,11, 13,19,21)$ = $\{1,3,9, 11,17,19\}$,  

$\theta_{24,2,1}(1,7,9, 15,17,23)$ = $\{3,9,11, 17,19,1\}$ = $\{1,3,9, 11,17,19\}$,  

$\theta_{24,2,3}(1,7,9, 15,17,23)$ = $\{7,13,15, 21,23,5\}$ = $\{5,7,13, 15,21,23\}$,  

$\theta_{24,2,5}(1,7,9, 15,17,23)$ = $\{11,17,19, 1,3,9\}$ = $\{1,3,9, 11,17,19\}$,  

$\theta_{24,2,7}(1,7,9, 15,17,23)$ = $\{15,21,23, 5,7,13\}$ = $\{5,7,13, 15,21,23\}$,  

$\theta_{24,2,9}(1,7,9, 15,17,23)$ = $\{19,1,3, 9,11,17\}$ = $\{1,3,9, 11,17,19\}$,

$\theta_{24,2,11}(1,7,9, 15,17,23)$ = $\{23,5,7, 13,15,21\}$ = $\{5,7,13, 15,21,23\}$.  

$\Rightarrow$ For $t$ = 1,3,5,7,9,11, $\theta_{24,2,t}(C_{24}(3,5,11)), \theta_{24,2,t}(C_{24}(1,7,9))$  $\neq$ $C_{24}(S)$ for any $S$  and $C_{24}(3,5,11)$ and $C_{24}(1,7,9)$ have no Type-2 isomorphic circulant graph.  

$\Rightarrow$ Circulant graphs $C_{24}(3,5,11)$ and $C_{24}(1,7,9)$  are Type-1 isomorphic and each has CI-property.  

\item [\rm (h3)]  Using the calculations involved in $(g7)$, we get,

$Ad_{24}(C_{24}(1,7,9))$ = $\{C_{24}(1,7,9), C_{24}(3,5,11) = C_{24}(5(1,7,9))\}$ =$Ad_{24}(C_{24}(3,5,11))$,

$\theta_{24,2,2}(C_{24}(1,7,9))$ = $C_{24}(3,5,11)$ = $\theta_{24,2,6}(C_{24}(1,7,9))$ = $\theta_{24,2,10}(C_{24}(1,7,9))$,   

$\theta_{24,2,4}(C_{24}(1,7,9))$ = $C_{24}(1,7,9)$ = $\theta_{24,2,8}(C_{24}(1,7,9))$,   

$\theta_{24,2,2}(C_{24}(3,5,11))$ = $C_{24}(1,7,9)$ = $\theta_{24,2,6}(C_{24}(3,5,11))$ = $\theta_{24,2,10}(C_{24}(3,5,11))$,   

$\theta_{24,2,4}(C_{24}(3,5,11))$ = $C_{24}(3,5,11)$ = $\theta_{24,2,8}(C_{24}(3,5,11))$, and 
\\
for $t$ = 1,3,5,7,9,11, $\theta_{24,2,t}(C_{24}(3,5,11)), \theta_{24,2,t}(C_{24}(1,7,9))$  $\neq$ $C_{24}(S)$ for any $S$ and $C_{24}(3,5,11)$ and $C_{24}(1,7,9)$ have no Type-2 isomorphic circulant graph.  

$\Rightarrow$ Circulant graphs $C_{24}(1,7,9)$ and $C_{24}(3,5,11)$ are Type-1 isomorphic and each has CI-property.  

Proof to all the cases given in the problem are similar, except $(g7)$ and $(h3)$, and in each case, we calculate $\theta_{24,2,3}(C_{24}(R))$, $\theta_{24,2,6}(C_{24}(R))$, $\theta_{24,2,9}(C_{24}(R))$ and $T1_{24}(C_{24}(R))$ corresponding to circulant graph $C_{24}(R)$ and show that all circulant graphs $C_{24}(S)$ which are $\theta_{24,2,3t}(C_{24}(R))$ for $t$ = 1 or 2 or 3 and belong to $T1_{24}(C_{24}(R))$. To simplify our work, calculated values of $\theta_{24,2,3}(C_{24}(R))$, $\theta_{24,2,6}(C_{24}(R))$, $\theta_{24,2,9}(C_{24}(R))$ and $T1_{24}(C_{24}(R))$ corresponding to each $C_{24}(R)$ are presented in Tables 5 to 8. For the cases $(g7)$ and $(h3)$, solutions are already presented in details. See Tables 5 to 8. 
\end{enumerate}

Thus, we could establish that $C_{24}(R)$ is not having Type-2 isomorphic circulant graph in each case. From the values of these tables, we get the result. \hfill $\Box$

\begin{table}
	\caption{ Calculation of $\theta_{24,2,3t}(C_{24}(R))$ and $T1_{24}(C_{24}(R))$ for $t$ = 1,2,3.}
\begin{center}
\scalebox{.85}{
 \begin{tabular}{||c||*{6}{c|}|c||}\hline \hline 
  & & & & & \\
S. No. & $C_{24}(R)$ &  $\theta_{24,2,3}(C_{24}(R))$ & $\theta_{24,2,6}(C_{24}(R))$ & $\theta_{24,2,9}(C_{24}(R))$ & $C_{24}(S)\in T1_{24}(C_{24}(R))$   \\

 & & & & &   \\ \hline \hline

(a1) & $C_{24}(1,2,3)$ & $\neq$ $C_{24}(S)$ & = $C_{24}(2,9,11)$ & $\neq$ $C_{24}(S)$ & $C_{24}(1,2,3)$, $C_{24}(5,9,10)$, \\ 

  & & & & & $C_{24}(3,7,10)$, $C_{24}(2,9,11)$. \\

(a2) & $C_{24}(1,2,5)$ &  $\neq$ $C_{24}(S)$ & = $C_{24}(2,7,11)$ & $\neq$ $C_{24}(S)$  &   $C_{24}(1,2,5)$, $C_{24}(1,5,10)$, \\ 

 & & & & & $C_{24}(7,10,11)$, $C_{24}(2,7,11)$. \\

(a3) & $C_{24}(1,2,7)$ &  $\neq$ $C_{24}(S)$ & = $C_{24}(2,5,11)$ & $\neq$ $C_{24}(S)$  &  $C_{24}(1,2,7)$, $C_{24}(5,10,11)$, \\ 

 & & & & & $C_{24}(1,7,10)$, $C_{24}(2,5,11)$. \\

(a4) & $C_{24}(1,2,9)$ & $\neq$ $C_{24}(S)$ & = $C_{24}(2,3,11)$ & $\neq$ $C_{24}(S)$ & $C_{24}(1,2,9)$, $C_{24}(3,5,10)$, \\ 

 & & & & & $C_{24}(7,9,10)$, $C_{24}(2,3,11)$. \\

(a5) & $C_{24}(2,3,5)$ &  $\neq$ $C_{24}(S)$ & = $C_{24}(2,7,9)$ & $\neq$ $C_{24}(S)$ & $C_{24}(2,3,5)$, $C_{24}(1,9,10)$, \\ 

  & & & & & $C_{24}(3,10,11)$, $C_{24}(2,7,9)$. \\

(a6) & $C_{24}(2,3,7)$ &  $\neq$ $C_{24}(S)$ & = $C_{24}(2,5,9)$ & $\neq$ $C_{24}(S)$ & $C_{24}(2,3,7)$, $C_{24}(9,10,11)$, \\ 

  & & & & & $C_{24}(1,3,10)$, $C_{24}(2,5,9)$. \\

(a7) & $C_{24}(2,3,9)$ & = $C_{24}(2,3,9)$ & = $C_{24}(2,3,9)$ & = $C_{24}(2,3,9)$ & $C_{24}(2,3,9)$, $C_{24}(3,9,10)$. \\

  & & & & &  \\

(a8) & $C_{24}(2,3,11)$ &  $\neq$ $C_{24}(S)$ & = $C_{24}(1,2,9)$ & $\neq$ $C_{24}(S)$ & $C_{24}(2,3,11)$, $C_{24}(7,9,10)$, \\ 

  & & & & & $C_{24}(3,5,10)$, $C_{24}(1,2,9)$. \\

(a9) & $C_{24}(2,5,9)$ &  $\neq$ $C_{24}(S)$ & = $C_{24}(2,3,7)$ & $\neq$ $C_{24}(S)$ & $C_{24}(2,5,9)$, $C_{24}(1,3,10)$, \\ 

  & & & & & $C_{24}(9,10,11)$, $C_{24}(2,3,7)$. \\

(a10) & $C_{24}(2,5,11)$ &  $\neq$ $C_{24}(S)$ & = $C_{24}(1,2,7)$ & $\neq$ $C_{24}(S)$ & $C_{24}(2,5,11)$, $C_{24}(1,7,10)$, \\ 

  & & & & & $C_{24}(5,10,11)$, $C_{24}(1,2,7)$. \\

(a11) & $C_{24}(2,7,9)$ &  $\neq$ $C_{24}(S)$ & = $C_{24}(2,3,5)$ & $\neq$ $C_{24}(S)$ & $C_{24}(2,7,9)$, $C_{24}(3,10,11)$, \\ 

  & & & & & $C_{24}(1,9,10)$, $C_{24}(2,3,5)$. \\

(a12) & $C_{24}(2,7,11)$ &  $\neq$ $C_{24}(S)$ & = $C_{24}(1,2,5)$ & $\neq$ $C_{24}(S)$ & $C_{24}(2,7,11)$, $C_{24}(7,10,11)$, \\ 

  & & & & & $C_{24}(1,5,10)$, $C_{24}(1,2,5)$. \\

(a13) & $C_{24}(2,9,11)$ &  $\neq$ $C_{24}(S)$ & = $C_{24}(1,2,3)$ & $\neq$ $C_{24}(S)$ & $C_{24}(2,9,11)$, $C_{24}(3,7,10)$, \\ 

  & & & & & $C_{24}(5,9,10)$, $C_{24}(1,2,3)$. \\

(b1) & $C_{24}(1,3,10)$ &  $\neq$ $C_{24}(S)$ & = $C_{24}(9,10,11)$ & $\neq$ $C_{24}(S)$ & $C_{24}(1,3,10)$, $C_{24}(2,5,9)$, \\ 

  & & & & & $C_{24}(2,3,7)$, $C_{24}(9,10,11)$. \\

(b2) & $C_{24}(1,5,10)$ & $\neq$ $C_{24}(S)$ & = $C_{24}(7,10,11)$ & $\neq$ $C_{24}(S)$ & $C_{24}(1,5,10)$, $C_{24}(1,2,5)$, \\ 

  & & & & & $C_{24}(2,7,11)$, $C_{24}(7,10,11)$. \\

(b3) & $C_{24}(1,7,10)$ &  $\neq$ $C_{24}(S)$ & = $C_{24}(5,10,11)$ & $\neq$ $C_{24}(S)$ & $C_{24}(1,7,10)$, $C_{24}(2,5,11)$, \\ 

  & & & & & $C_{24}(1,2,7)$, $C_{24}(5,10,11)$. \\

(b4) & $C_{24}(1,9,10)$ &  $\neq$ $C_{24}(S)$ & = $C_{24}(3,10,11)$ & $\neq$ $C_{24}(S)$ & $C_{24}(1,9,10)$, $C_{24}(2,3,5)$, \\ 

  & & & & & $C_{24}(2,7,9)$, $C_{24}(3,10,11)$. \\

(b5) & $C_{24}(3,5,10)$ &  $\neq$ $C_{24}(S)$ & = $C_{24}(7,9,10)$ & $\neq$ $C_{24}(S)$ & $C_{24}(3,5,10)$, $C_{24}(1,2,9)$, \\ 

  & & & & & $C_{24}(2,3,11)$, $C_{24}(7,9,10)$. \\

(b6) & $C_{24}(3,7,10)$ &  $\neq$ $C_{24}(S)$ & = $C_{24}(5,9,10)$ & $\neq$ $C_{24}(S)$ & $C_{24}(3,7,10)$, $C_{24}(2,9,11)$, \\ 

  & & & & & $C_{24}(1,2,3)$, $C_{24}(5,9,10)$. \\

(b7) & $C_{24}(3,9,10)$ & = $C_{24}(3,9,10)$ & = $C_{24}(3,9,10)$ & = $C_{24}(3,9,10)$ & $C_{24}(3,9,10)$, $C_{24}(2,3,9)$. \\

  & & & & &  \\

(b8) & $C_{24}(3,10,11)$ &  $\neq$ $C_{24}(S)$ & = $C_{24}(1,9,10)$ & $\neq$ $C_{24}(S)$ & $C_{24}(3,10,11)$, $C_{24}(2,7,9)$, \\ 

  & & & & & $C_{24}(2,3,5)$, $C_{24}(1,9,10)$. \\

(b9) & $C_{24}(5,9,10)$ & $\neq$ $C_{24}(S)$ & = $C_{24}(3,7,10)$ & $\neq$ $C_{24}(S)$ & $C_{24}(5,9,10)$, $C_{24}(1,2,3)$, \\ 

  & & & & & $C_{24}(2,9,11)$, $C_{24}(3,7,10)$. \\

(b10) & $C_{24}(5,10,11)$ & $\neq$ $C_{24}(S)$ & = $C_{24}(1,7,10)$ & $\neq$ $C_{24}(S)$ & $C_{24}(5,10,11)$, $C_{24}(1,2,7)$, \\ 

  & & & & & $C_{24}(2,5,11)$, $C_{24}(1,7,10)$. \\

(b11) & $C_{24}(7,9,10)$ & $\neq$ $C_{24}(S)$ & = $C_{24}(3,5,10)$ & $\neq$ $C_{24}(S)$ & $C_{24}(7,9,10)$, $C_{24}(2,3,11)$, \\ 

  & & & & & $C_{24}(1,2,9)$, $C_{24}(3,5,10)$. \\

(b12) & $C_{24}(7,10,11)$ & $\neq$ $C_{24}(S)$ & = $C_{24}(1,5,10)$ & $\neq$ $C_{24}(S)$ & $C_{24}(7,10,11)$, $C_{24}(2,7,11)$, \\ 

  & & & & & $C_{24}(1,2,5)$, $C_{24}(1,5,10)$. \\

(b13) & $C_{24}(9,10,11)$ & $\neq$ $C_{24}(S)$ & = $C_{24}(1,3,10)$ & $\neq$ $C_{24}(S)$ & $C_{24}(9,10,11)$, $C_{24}(2,3,7)$, \\ 

  & & & & & $C_{24}(2,5,9)$, $C_{24}(1,3,10)$. \\  \hline \hline
\end{tabular}}
\end{center}
\end{table} 

\begin{table}
	\caption{ Calculation of $\theta_{24,2,3t}(C_{24}(R))$ and $T1_{24}(C_{24}(R))$ for $t$ = 1,2,3.}
\begin{center}
\scalebox{.85}{
 \begin{tabular}{||c||*{6}{c|}|c||c||}\hline \hline 
  & & & & &  \\
S. No. & $C_{24}(R)$ &  $\theta_{24,2,3}(C_{24}(R))$ & $\theta_{24,2,6}(C_{24}(R))$ & $\theta_{24,2,9}(C_{24}(R))$ & $C_{24}(S)\in T1_{24}(C_{24}(R))$   \\ \hline \hline

  & & & & &  \\ 

(c1) & $C_{24}(1,3,4)$ &  $\neq$ $C_{24}(S)$ & = $C_{24}(4,9,11)$ & $\neq$ $C_{24}(S)$ & $C_{24}(1,3,4)$, $C_{24}(4,5,9)$, \\ 

  & & & & & $C_{24}(3,4,7)$, $C_{24}(4,9,11)$. \\

(c2) & $C_{24}(1,4,5)$ &  $\neq$ $C_{24}(S)$ & = $C_{24}(4,7,11)$ & $\neq$ $C_{24}(S)$ & $C_{24}(1,4,5)$, $C_{24}(4,7,11)$. \\  

(c3) & $C_{24}(1,4,7)$ & $\neq$ $C_{24}(S)$ & = $C_{24}(4,5,11)$ & $\neq$ $C_{24}(S)$ & $C_{24}(1,4,7)$, $C_{24}(4,5,11)$. \\  

(c4) & $C_{24}(1,4,9)$ & $\neq$ $C_{24}(S)$ & = $C_{24}(3,4,11)$ & $\neq$ $C_{24}(S)$ & $C_{24}(1,4,9)$, $C_{24}(3,4,5)$, \\ 

 & & & & & $C_{24}(4,7,9)$, $C_{24}(3,4,11)$. \\  

(c5) & $C_{24}(1,4,11)$ &  = $C_{24}(4,5,7)$ & = $C_{24}(1,4,11)$ & = $C_{24}(4,5,7)$ & $C_{24}(1,4,11)$, $C_{24}(4,5,7)$. \\

  & & & & &  \\

(c6) & $C_{24}(3,4,5)$ &  $\neq$ $C_{24}(S)$ & = $C_{24}(4,7,9)$ & $\neq$ $C_{24}(S)$ & $C_{24}(3,4,5)$, $C_{24}(1,4,9)$, \\ 

  & & & & & $C_{24}(3,4,11)$, $C_{24}(4,7,9)$. \\  

(c7) & $C_{24}(3,4,7)$ &  $\neq$ $C_{24}(S)$ & = $C_{24}(4,5,9)$ & $\neq$ $C_{24}(S)$ & $C_{24}(3,4,7)$, $C_{24}(4,9,11)$, \\ 

  & & & & & $C_{24}(1,3,4)$, $C_{24}(4,5,9)$. \\  

(c8) & $C_{24}(3,4,9)$ & = $C_{24}(3,4,9)$ & = $C_{24}(3,4,9)$ & = $C_{24}(3,4,9)$ & $C_{24}(3,4,9)$. \\  

& & & & & \\

(c9) & $C_{24}(3,4,11)$ &  $\neq$ $C_{24}(S)$ & = $C_{24}(1,4,9)$ & $\neq$ $C_{24}(S)$ & $C_{24}(3,4,11)$, $C_{24}(4,7,9)$, \\ 

  & & & & & $C_{24}(3,4,5)$, $C_{24}(1,4,9)$. \\  

(c10) & $C_{24}(4,5,7)$ & = $C_{24}(1,4,11)$ & = $C_{24}(4,5,7)$ & = $C_{24}(1,4,11)$ & $C_{24}(4,5,7)$, $C_{24}(1,4,11)$. \\  

  & & & & &  \\

(c11) & $C_{24}(4,5,9)$ & $\neq$ $C_{24}(S)$ & = $C_{24}(3,4,7)$ & $\neq$ $C_{24}(S)$ & $C_{24}(4,5,9)$, $C_{24}(1,3,4)$, \\ 

  & & & & & $C_{24}(4,9,11)$, $C_{24}(3,4,7)$. \\  

(c12) & $C_{24}(4,5,11)$ &  $\neq$ $C_{24}(S)$ & = $C_{24}(1,4,7)$ & $\neq$ $C_{24}(S)$ & $C_{24}(4,5,11)$, $C_{24}(1,4,7)$. \\ 

  & & & & & \\  

(c13) & $C_{24}(4,7,9)$ &  $\neq$ $C_{24}(S)$ & = $C_{24}(3,4,5)$ & $\neq$ $C_{24}(S)$ & $C_{24}(4,7,9)$, $C_{24}(3,4,11)$, \\ 

  & & & & & $C_{24}(1,4,9)$, $C_{24}(3,4,5)$. \\  

(c14) & $C_{24}(4,7,11)$ & $\neq$ $C_{24}(S)$ & = $C_{24}(1,4,5)$ & $\neq$ $C_{24}(S)$ & $C_{24}(4,7,11)$, $C_{24}(1,4,5)$. \\ 

& & & & & \\ 

(c15) & $C_{24}(4,9,11)$ & $\neq$ $C_{24}(S)$ & = $C_{24}(1,3,4)$ & $\neq$ $C_{24}(S)$ & $C_{24}(4,9,11)$, $C_{24}(3,4,7)$, \\ 

  & & & & & $C_{24}(4,5,9)$, $C_{24}(1,3,4)$ . \\  

(d1) & $C_{24}(1,3,6)$ & $\neq$ $C_{24}(S)$ & = $C_{24}(6,9,11)$ & $\neq$ $C_{24}(S)$ & $C_{24}(1,3,6)$, $C_{24}(5,6,9)$, \\ 

&  & & & & $C_{24}(3,6,7)$, $C_{24}(6,9,11)$. \\

(d2) & $C_{24}(1,5,6)$ & $\neq$ $C_{24}(S)$ & = $C_{24}(6,7,11)$ & $\neq$ $C_{24}(S)$ & $C_{24}(1,5,6)$, $C_{24}(6,7,11)$. \\

(d3) & $C_{24}(1,6,7)$ & $\neq$ $C_{24}(S)$ & = $C_{24}(5,6,11)$ & $\neq$ $C_{24}(S)$ & $C_{24}(1,6,7)$, $C_{24}(5,6,11)$. \\

(d4) & $C_{24}(1,6,9)$ & $\neq$ $C_{24}(S)$ & = $C_{24}(3,6,11)$ & $\neq$ $C_{24}(S)$ & $C_{24}(1,6,9)$, $C_{24}(3,5,6)$, \\ 

  & & & & & $C_{24}(6,7,9)$, $C_{24}(3,6,11)$. \\

(d5) & $C_{24}(1,6,11)$ & = $C_{24}(5,6,7)$ & = $C_{24}(1,6,11)$ & = $C_{24}(5,6,7)$ & $C_{24}(1,6,11)$, $C_{24}(5,6,7)$. \\

 & & & & & \\

(d6) & $C_{24}(3,5,6)$ & $\neq$ $C_{24}(S)$ & = $C_{24}(6,7,9)$ & $\neq$ $C_{24}(S)$ & $C_{24}(3,5,6)$, $C_{24}(1,6,9)$, \\ 

  & & & & & $C_{24}(3,6,11)$, $C_{24}(6,7,9)$. \\

(d7) & $C_{24}(3,6,7)$ & $\neq$ $C_{24}(S)$ & = $C_{24}(5,6,9)$ & $\neq$ $C_{24}(S)$ & $C_{24}(3,6,7)$, $C_{24}(6,9,11)$, \\ 

  & & & & & $C_{24}(1,3,6)$, $C_{24}(5,6,9)$. \\

(d8) & $C_{24}(3,6,9)$ & = $C_{24}(3,6,9)$ & = $C_{24}(3,6,9)$ & = $C_{24}(3,6,9)$ & $C_{24}(3,6,9)$. \\

 & & & & & \\

(d9) & $C_{24}(3,6,11)$ & $\neq$ $C_{24}(S)$ & = $C_{24}(1,6,9)$ & $\neq$ $C_{24}(S)$ & $C_{24}(3,6,11)$, $C_{24}(6,7,9)$, \\ 

  & & & & & $C_{24}(3,5,6)$, $C_{24}(1,6,9)$. \\

(d10) & $C_{24}(5,6,7)$ & = $C_{24}(1,6,11)$ & = $C_{24}(5,6,7)$ & = $C_{24}(1,6,11)$ & $C_{24}(5,6,7)$, $C_{24}(1,6,11)$. \\

 & & & & & \\

(d11) & $C_{24}(5,6,9)$ & $\neq$ $C_{24}(S)$ & = $C_{24}(3,6,7)$ & $\neq$ $C_{24}(S)$ & $C_{24}(5,6,9)$, $C_{24}(1,3,6)$, \\ 

&  & & & & $C_{24}(6,9,11)$, $C_{24}(3,6,7)$. \\

(d12) & $C_{24}(5,6,11)$ & $\neq$ $C_{24}(S)$ & = $C_{24}(1,6,7)$ & $\neq$ $C_{24}(S)$ & $C_{24}(5,6,11)$, $C_{24}(1,6,7)$. \\

(d13) & $C_{24}(6,7,9)$ & $\neq$ $C_{24}(S)$ & = $C_{24}(3,5,6)$ & $\neq$ $C_{24}(S)$ & $C_{24}(6,7,9)$, $C_{24}(3,6,11)$, \\ 

  & & & & & $C_{24}(1,6,9)$, $C_{24}(3,5,6)$. \\

(d14) & $C_{24}(6,7,11)$ & $\neq$ $C_{24}(S)$ & = $C_{24}(1,5,6)$ & $\neq$ $C_{24}(S)$ & $C_{24}(6,7,11)$, $C_{24}(1,5,6)$. \\

(d15) & $C_{24}(6,9,11)$ & $\neq$ $C_{24}(S)$ & = $C_{24}(1,3,6)$ & $\neq$ $C_{24}(S)$ & $C_{24}(6,9,11)$, $C_{24}(3,6,7)$, \\ 

  & & & & & $C_{24}(5,6,9)$, $C_{24}(1,3,6)$. \\  \hline \hline
\end{tabular}}
\end{center}
\end{table} 

\begin{table}
	\caption{ Calculation of $\theta_{24,2,3t}(C_{24}(R))$ and $T1_{24}(C_{24}(R))$ for $t$ = 1,2,3.}
\begin{center}
\scalebox{.85}{
 \begin{tabular}{||c||*{6}{c|}|c||c||}\hline \hline 
  & & & & & \\
S. No. & $C_{24}(R)$ &  $\theta_{24,2,3}(C_{24}(R))$ & $\theta_{24,2,6}(C_{24}(R))$ & $\theta_{24,2,9}(C_{24}(R))$ & $C_{24}(S)\in T1_{24}(C_{24}(R))$   \\

  & & & & &  \\ \hline \hline

(e1) & $C_{24}(1,3,8)$ & $\neq$ $C_{24}(S)$ & = $C_{24}(8,9,11)$ & $\neq$ $C_{24}(S)$ & $C_{24}(1,3,8)$,  $C_{24}(5,8,9)$, \\ 

&  & & & & $C_{24}(3,7,8)$, $C_{24}(8,9,11)$. \\

(e2) & $C_{24}(1,5,8)$ & $\neq$ $C_{24}(S)$ & = $C_{24}(7,8,11)$ & $\neq$ $C_{24}(S)$ & $C_{24}(1,5,8)$, $C_{24}(7,8,11)$. \\

 & & & & & \\

(e3) & $C_{24}(1,7,8)$ & $\neq$ $C_{24}(S)$ & = $C_{24}(5,8,11)$ & $\neq$ $C_{24}(S)$ & $C_{24}(1,7,8)$, $C_{24}(5,8,11)$. \\

 & & & & & \\

(e4) & $C_{24}(1,8,9)$ & $\neq$ $C_{24}(S)$ & = $C_{24}(3,8,11)$ & $\neq$ $C_{24}(S)$ & $C_{24}(1,8,9)$,  $C_{24}(3,5,8)$, \\ 

 & & & & & $C_{24}(7,8,9)$, $C_{24}(3,8,11)$. \\

(e5) & $C_{24}(1,8,11)$ & = $C_{24}(5,7,8)$ & = $C_{24}(1,8,11)$ & = $C_{24}(5,7,8)$ & $C_{24}(1,8,11)$, $C_{24}(5,7,8)$. \\

(e6) & $C_{24}(3,5,8)$ & $\neq$ $C_{24}(S)$ & = $C_{24}(7,8,9)$ & $\neq$ $C_{24}(S)$ & $C_{24}(3,5,8)$,  $C_{24}(1,8,9)$, \\ 

&  & & & & $C_{24}(3,8,11)$, $C_{24}(7,8,9)$. \\

(e7) & $C_{24}(3,7,8)$ & $\neq$ $C_{24}(S)$ & = $C_{24}(5,8,9)$ & $\neq$ $C_{24}(S)$ & $C_{24}(3,7,8)$,  $C_{24}(8,9,11)$, \\ 

  & & & & & $C_{24}(1,3,8)$, $C_{24}(5,8,9)$. \\

(e8) & $C_{24}(3,8,9)$ & = $C_{24}(3,8,9)$ & = $C_{24}(3,8,9)$ & = $C_{24}(3,8,9)$ & $C_{24}(3,8,9)$. \\

 & & & & & \\

(e9) & $C_{24}(3,8,11)$ & $\neq$ $C_{24}(S)$ & = $C_{24}(1,8,9)$ & $\neq$ $C_{24}(S)$ & $C_{24}(3,8,11)$,  $C_{24}(7,8,9)$, \\ 

  & & & & & $C_{24}(3,5,8)$, $C_{24}(1,8,9)$. \\

(e10) & $C_{24}(5,7,8)$ & = $C_{24}(1,8,11)$ & = $C_{24}(5,7,8)$ & = $C_{24}(1,8,11)$ & $C_{24}(5,7,8)$, $C_{24}(1,8,11)$. \\

 & & & & & \\

(e11) & $C_{24}(5,8,9)$ & $\neq$ $C_{24}(S)$ & = $C_{24}(3,7,8)$ & $\neq$ $C_{24}(S)$ & $C_{24}(5,8,9)$,  $C_{24}(1,3,8)$, \\ 

  & & & & & $C_{24}(8,9,11)$, $C_{24}(3,7,8)$. \\

(e12) & $C_{24}(5,8,11)$ & $\neq$ $C_{24}(S)$ & = $C_{24}(1,7,8)$ & $\neq$ $C_{24}(S)$ & $C_{24}(5,8,11)$,  $C_{24}(1,7,8)$. \\

 & & & & & \\

(e13) & $C_{24}(7,8,9)$ & $\neq$ $C_{24}(S)$ & = $C_{24}(3,5,8)$ & $\neq$ $C_{24}(S)$ & $C_{24}(7,8,9)$,  $C_{24}(3,8,11)$, \\ 

&  & & & & $C_{24}(1,8,9)$, $C_{24}(3,5,8)$. \\

(e14) & $C_{24}(7,8,11)$ & $\neq$ $C_{24}(S)$ & = $C_{24}(1,5,8)$ & $\neq$ $C_{24}(S)$ & $C_{24}(7,8,11)$, $C_{24}(1,5,8)$. \\

 (e15) & $C_{24}(8,9,11)$ & $\neq$ $C_{24}(S)$ & = $C_{24}(1,3,8)$ & $\neq$ $C_{24}(S)$ & $C_{24}(8,9,11)$,  $C_{24}(3,7,8)$, \\ 

  & & & & & $C_{24}(5,8,9)$, $C_{24}(1,3,8)$. \\

(f1) & $C_{24}(1,3,12)$ & $\neq$ $C_{24}(S)$ & = $C_{24}(9,11,12)$ & $\neq$ $C_{24}(S)$ & $C_{24}(1,3,12)$, $C_{24}(5,9,12)$, \\ 

  & & & & & $C_{24}(3,7,12)$, $C_{24}(9,11,12)$. \\

(f2) & $C_{24}(1,5,12)$ & $\neq$ $C_{24}(S)$ & = $C_{24}(7,11,12)$ & $\neq$ $C_{24}(S)$ & $C_{24}(1,5,12)$, $C_{24}(7,11,12)$. \\

(f3) & $C_{24}(1,7,12)$ & $\neq$ $C_{24}(S)$ & = $C_{24}(5,11,12)$ & $\neq$ $C_{24}(S)$ & $C_{24}(1,7,12)$, $C_{24}(5,11,12)$. \\

(f4) & $C_{24}(1,9,12)$ & $\neq$ $C_{24}(S)$ & = $C_{24}(3,11,12)$ & $\neq$ $C_{24}(S)$ & $C_{24}(1,9,12)$, $C_{24}(3,5,12)$, \\ 

  & & & & & $C_{24}(7,9,12)$, $C_{24}(3,11,12)$. \\

(f5) & $C_{24}(1,11,12)$ & = $C_{24}(5,7,12)$ & = $C_{24}(1,11,12)$ & = $C_{24}(5,7,12)$ & $C_{24}(1,11,12)$, $C_{24}(5,7,12)$. \\

 & & & & & \\

(f6) & $C_{24}(3,5,12)$ & $\neq$ $C_{24}(S)$ & = $C_{24}(7,9,12)$ & $\neq$ $C_{24}(S)$ & $C_{24}(3,5,12)$, $C_{24}(1,9,12)$, \\ 

  & & & & & $C_{24}(3,11,12)$, $C_{24}(7,9,12)$. \\

(f7) & $C_{24}(3,7,12)$ & $\neq$ $C_{24}(S)$ & = $C_{24}(5,9,12)$ & $\neq$ $C_{24}(S)$ & $C_{24}(3,7,12)$, $C_{24}(9,11,12)$, \\ 

&  & & & & $C_{24}(1,3,12)$, $C_{24}(5,9,12)$. \\

(f8) & $C_{24}(3,9,12)$ & = $C_{24}(3,9,12)$ & = $C_{24}(3,9,12)$ & = $C_{24}(3,9,12)$ & $C_{24}(3,9,12)$. \\

 & & & & & \\

(f9) & $C_{24}(3,11,12)$ & $\neq$ $C_{24}(S)$ & = $C_{24}(1,9,12)$ & $\neq$ $C_{24}(S)$ & $C_{24}(3,11,12)$, $C_{24}(7,9,12)$, \\ 

  & & & & & $C_{24}(3,5,12)$, $C_{24}(1,9,12)$. \\

(f10) & $C_{24}(5,7,12)$ & = $C_{24}(1,11,12)$ & = $C_{24}(5,7,12)$ & = $C_{24}(1,11,12)$ & $C_{24}(5,7,12)$, $C_{24}(1,11,12)$. \\

 & & & & & \\

(f11) & $C_{24}(5,9,12)$ & $\neq$ $C_{24}(S)$ & = $C_{24}(3,7,12)$ & $\neq$ $C_{24}(S)$ & $C_{24}(5,9,12)$, $C_{24}(1,3,12)$, \\ 

  & & & & & $C_{24}(9,11,12)$, $C_{24}(3,7,12)$. \\

(f12) & $C_{24}(5,11,12)$ & $\neq$ $C_{24}(S)$ & = $C_{24}(1,7,12)$ & $\neq$ $C_{24}(S)$ & $C_{24}(5,11,12)$, $C_{24}(1,7,12)$. \\

(f13) & $C_{24}(7,9,12)$ & $\neq$ $C_{24}(S)$ & = $C_{24}(3,5,12)$ & $\neq$ $C_{24}(S)$ & $C_{24}(7,9,12)$, $C_{24}(3,11,12)$, \\ 

  & & & & & $C_{24}(1,9,12)$, $C_{24}(3,5,12)$. \\

(f14) & $C_{24}(7,11,12)$ & $\neq$ $C_{24}(S)$ & = $C_{24}(1,5,12)$ & $\neq$ $C_{24}(S)$ & $C_{24}(7,11,12)$, $C_{24}(1,5,12)$. \\

(f15) & $C_{24}(9,11,12)$ & $\neq$ $C_{24}(S)$ & = $C_{24}(1,3,12)$ & $\neq$ $C_{24}(S)$ & $C_{24}(9,11,12)$, $C_{24}(3,7,12)$, \\ 

&  & & & & $C_{24}(5,9,12)$, $C_{24}(1,3,12)$. \\  \hline \hline
\end{tabular}}
\end{center}
\end{table} 

\begin{table}
	\caption{ Calculation of $\theta_{24,2,3t}(C_{24}(R))$ and $T1_{24}(C_{24}(R))$ for $t$ = 1,2,3.}
\begin{center}
\scalebox{.85}{
 \begin{tabular}{||c||*{6}{c|}|c||c||}\hline \hline 
  & & & & &  \\
S. No. & $C_{24}(R)$ &  $\theta_{24,2,3}(C_{24}(R))$ & $\theta_{24,2,6}(C_{24}(R))$ & $\theta_{24,2,9}(C_{24}(R))$ & $C_{24}(S)\in T1_{24}(C_{24}(R))$   \\

  & & & & &  \\ \hline \hline
	
  & & & & & \\ 

(g1) & $C_{24}(1,3,5)$ & $\neq$ $C_{24}(S)$ & = $C_{24}(7,9,11)$ & $\neq$ $C_{24}(S)$ & $C_{24}(1,3,5)$, $C_{24}(1,5,9)$, \\ 

& & & & & $C_{24}(3,7,11)$, $C_{24}(7,9,11)$. \\

(g2) & $C_{24}(1,3,7)$ & $\neq$ $C_{24}(S)$ & = $C_{24}(5,9,11)$ & $\neq$ $C_{24}(S)$ & $C_{24}(1,3,7)$, $C_{24}(5,9,11)$. \\

 & & & & & \\

(g3) & $C_{24}(1,3,9)$ & $\neq$ $C_{24}(S)$ & = $C_{24}(3,9,11)$ & $\neq$ $C_{24}(S)$ & $C_{24}(1,3,9)$, $C_{24}(3,5,9)$, \\ 

  & & & & & $C_{24}(3,7,9)$, $C_{24}(3,9,11)$. \\

(g4) & $C_{24}(1,3,11)$ & $\neq$ $C_{24}(S)$ & = $C_{24}(1,9,11)$ & $\neq$ $C_{24}(S)$ & $C_{24}(1,3,11)$, $C_{24}(5,7,9)$, \\ 

  & & & & & $C_{24}(3,5,7)$, $C_{24}(1,9,11)$. \\

(g5) & $C_{24}(3,5,7)$ & $\neq$ $C_{24}(S)$ & = $C_{24}(5,7,9)$ & $\neq$ $C_{24}(S)$ & $C_{24}(3,5,7)$, $C_{24}(1,9,11)$, \\ 

  & & & & & $C_{24}(1,3,11)$, $C_{24}(5,7,9)$. \\

(g6) & $C_{24}(3,5,9)$ & $\neq$ $C_{24}(S)$ & = $C_{24}(3,7,9)$ & $\neq$ $C_{24}(S)$ & $C_{24}(3,5,9)$, $C_{24}(1,3,9)$, \\ 

  & & & & & $C_{24}(3,9,11)$, $C_{24}(3,7,9)$. \\  \hline 
	
	& & & & & \\

(g7) & $C_{24}(3,5,11)$ & $\neq$ $C_{24}(S)$ & = $C_{24}(1,7,9)$ & $\neq$ $C_{24}(S)$ & $C_{24}(3,5,11)$, $C_{24}(1,7,9)$. \\

   & Additional & circulant graph & values of & $\theta_{24,2,t}(C_{24}(3,5,11))$ are  & given at  the bottom \\
	
& &	& &  &  of the table.  \\ \hline 

 & & & & & \\

(g8) & $C_{24}(3,7,9)$ & $\neq$ $C_{24}(S)$ & = $C_{24}(3,5,9)$ & $\neq$ $C_{24}(S)$ & $C_{24}(3,7,9)$, $C_{24}(3,9,11)$, \\ 

&  & & & & $C_{24}(1,3,9)$, $C_{24}(3,5,9)$. \\

(g9) & $C_{24}(3,7,11)$ & $\neq$ $C_{24}(S)$ & = $C_{24}(1,5,9)$ & $\neq$ $C_{24}(S)$ & $C_{24}(3,7,11)$, $C_{24}(7,9,11)$, \\ 

  & & & & & $C_{24}(1,3,5)$, $C_{24}(1,5,9)$. \\

(g10) & $C_{24}(3,9,11)$ & $\neq$ $C_{24}(S)$ & = $C_{24}(1,3,9)$ & $\neq$ $C_{24}(S)$ & $C_{24}(3,9,11)$, $C_{24}(3,7,9)$, \\ 

  & & & & & $C_{24}(3,5,9)$, $C_{24}(1,3,9)$. \\

(h1) & $C_{24}(1,3,9)$ & $\neq$ $C_{24}(S)$ & = $C_{24}(3,9,11)$ & $\neq$ $C_{24}(S)$ & $C_{24}(1,3,9)$, $C_{24}(3,5,9)$, \\ 

  & & & & & $C_{24}(3,7,9)$, $C_{24}(3,9,11)$. \\

(h2) & $C_{24}(1,5,9)$ & $\neq$ $C_{24}(S)$ & = $C_{24}(3,7,11)$ & $\neq$ $C_{24}(S)$ & $C_{24}(1,5,9)$, $C_{24}(1,3,5)$, \\ 

  & & & & & $C_{24}(7,9,11)$, $C_{24}(3,7,11)$. \\ \hline 
	
(h3) & $C_{24}(1,7,9)$ & $\neq$ $C_{24}(S)$ & = $C_{24}(3,5,11)$ & $\neq$ $C_{24}(S)$ & $C_{24}(1,7,9)$, $C_{24}(3,5,11)$. \\

   & Additional & circulant graph & values of & $\theta_{24,2,t}(C_{24}(1,7,9))$ are  & given at  the bottom \\
	
& &	& &  & of the table.  \\ \hline 

 & & & & & \\

(h4) & $C_{24}(1,9,11)$ & $\neq$ $C_{24}(S)$ & = $C_{24}(1,3,11)$ & $\neq$ $C_{24}(S)$ & $C_{24}(1,9,11)$, $C_{24}(3,5,7)$, \\ 

&  & & & & $C_{24}(5,7,9)$, $C_{24}(1,3,11)$. \\

(h5) & $C_{24}(3,5,9)$ & $\neq$ $C_{24}(S)$ & = $C_{24}(3,7,9)$ & $\neq$ $C_{24}(S)$ & $C_{24}(3,5,9)$, $C_{24}(1,3,9)$, \\ 

  & & & & & $C_{24}(3,9,11)$, $C_{24}(3,7,9)$. \\

(h6) & $C_{24}(3,7,9)$ & $\neq$ $C_{24}(S)$ & = $C_{24}(3,5,9)$ & $\neq$ $C_{24}(S)$ & $C_{24}(3,7,9)$, $C_{24}(3,9,11)$, \\ 

  & & & & & $C_{24}(1,3,9)$, $C_{24}(3,5,9)$. \\

(h7) & $C_{24}(3,9,11)$ & $\neq$ $C_{24}(S)$ & = $C_{24}(1,3,9)$ & $\neq$ $C_{24}(S)$ & $C_{24}(3,9,11)$, $C_{24}(3,7,9)$, \\ 

  & & & & & $C_{24}(3,5,9)$, $C_{24}(1,3,9)$. \\

(h8) & $C_{24}(5,7,9)$ & $\neq$ $C_{24}(S)$ & = $C_{24}(3,5,7)$ & $\neq$ $C_{24}(S)$ & $C_{24}(5,7,9)$, $C_{24}(1,3,11)$, \\ 

  & & & & & $C_{24}(1,9,11)$, $C_{24}(3,5,7)$. \\

(h9) & $C_{24}(5,9,11)$ & $\neq$ $C_{24}(S)$ & = $C_{24}(1,3,7)$ & $\neq$ $C_{24}(S)$ & $C_{24}(5,9,11)$, $C_{24}(1,3,7)$. \\

  & & & & &  \\ 

(h10) & $C_{24}(7,9,11)$ & $\neq$ $C_{24}(S)$ & = $C_{24}(1,3,5)$ & $\neq$ $C_{24}(S)$ & $C_{24}(7,9,11)$, $C_{24}(3,7,11)$, \\ 

&  & & & & $C_{24}(1,5,9)$, $C_{24}(1,3,5)$.     \\  \hline \hline
\end{tabular}}
\end{center}

\noindent
(g7) ~ $\theta_{24,2,2}(C_{24}(3,5,11))$ = $C_{24}(1,7,9)$, $\theta_{24,2,4}(C_{24}(3,5,11))$ = $C_{24}(3,5,11)$, 

\hfill $\theta_{24,2,8}(C_{24}(3,5,11))$ = $C_{24}(3,5,11)$, $\theta_{24,2,10}(C_{24}(3,5,11))$ = $C_{24}(1,7,9)$.

\vspace{.1cm}
\noindent
(h3) ~ $\theta_{24,2,2}(C_{24}(1,7,9))$ = $C_{24}(3,5,11)$, $\theta_{24,2,4}(C_{24}(1,7,9))$ = $C_{24}(1,7,9)$, ~~~

\hfill $\theta_{24,2,8}(C_{24}(1,7,9))$ = $C_{24}(1,7,9)$, $\theta_{24,2,10}(C_{24}(1,7,9))$ = $C_{24}(3,5,11)$.
\end{table} 

\begin{prm}\quad \label{p3} {\rm Show that the following pair of circulant graphs of order 24 are Type-1 isomorphic.
\begin{enumerate}
\item [\rm (1)] $C_{24}(1,3,11)$, $C_{24}(3,5,7)$;  

\item [\rm (2)] $C_{24}(1,9,11)$, $C_{24}(5,7,9)$; 

\item [\rm (3)] $C_{24}(1,3,6,11)$, $C_{24}(3,5,6,7)$;
	
\item [\rm (4)] $C_{24}(1,3,9,11)$, $C_{24}(3,5,7,9)$; 

\item [\rm (5)] $C_{24}(1,3,11,12)$, $C_{24}(3,5,7,12)$;

\item [\rm (6)] $C_{24}(1,6,9,11)$, $C_{24}(5,6,7,9)$; 

\item [\rm (7)] $C_{24}(1,9,11,12)$, $C_{24}(5,7,9,12)$; 

\item [\rm (8)] $C_{24}(1,3,6,9,11)$, $C_{24}(3,5,6,7,9)$;

\item [\rm (9)] $C_{24}(1,3,6,11,12)$, $C_{24}(3,5,6,7,12)$;

\item [\rm (10)] $C_{24}(1,3,9,11,12)$, $C_{24}(3,5,7,9,12)$;

\item [\rm (11)] $C_{24}(1,6,9,11,12)$, $C_{24}(5,6,7,9,12)$;

\item [\rm (12)] $C_{24}(1,3,6,9,11,12)$, $C_{24}(3,5,6,7,9,12)$. 
\end{enumerate} }
\end{prm}
\noindent
{\bf Solution.}\quad Given pairs of circulant graphs of order 24 are Type-1 isomorphic because of the following.
\begin{enumerate}
\item [\rm (1)] $C_{24}(7(1,3,11))$ = $C_{24}(7(1,3,11, 13,21,23))$ = $C_{24}(7,21,5, 19,,11))$ = $C_{24}(3,5,7)$. 

$\Rightarrow$  $C_{24}(1,3,11)$ and $C_{24}(3,5,7)$ are Type-1 isomorphic.
	
\item [\rm (2)] $C_{24}(7(1,9,11))$ = $C_{24}(5,7,9)$. 
$\Rightarrow$   $C_{24}(1,9,11)$ and $C_{24}(5,7,9)$  are Type-1 isomorphic.

\item [\rm (3)] $C_{24}(7(1,3,6,11))$ = $C_{24}(3,5,6,7)$. 

$\Rightarrow$  $C_{24}(1,3,6,11)$ and $C_{24}(3,5,6,7)$  are Type-1 isomorphic.
	
\item [\rm (4)] $C_{24}(5(1,3,9,11))$ = $C_{24}(3,5,7,9)$. 

$\Rightarrow$   $C_{24}(1,3,9,11)$ and $C_{24}(3,5,7,9)$  are Type-1 isomorphic.

\item [\rm (5)] $C_{24}(7(1,3,11,12))$ = $C_{24}(3,5,7,12)$. 

$\Rightarrow$  $C_{24}(1,3,11,12)$ and $C_{24}(3,5,7,12)$ are Type-1 isomorphic.

\item [\rm (6)] $C_{24}(7(1,6,9,11))$ = $C_{24}(5,6,7,9)$. 

$\Rightarrow$  $C_{24}(1,6,9,11)$ and $C_{24}(5,6,7,9)$ are Type-1 isomorphic.  

\item [\rm (7)] $C_{24}(7(1,9,11,12))$ = $C_{24}(5,7,9,12)$. 

$\Rightarrow$  $C_{24}(1,9,11,12)$ and $C_{24}(5,7,9,12)$ are Type-1 isomorphic.  

\item [\rm (8)] $C_{24}(5(1,3,6,9,11))$ = $C_{24}(3,5,6,7,9)$. 

$\Rightarrow$  $C_{24}(1,3,6,9,11)$ and $C_{24}(3,5,6,7,9)$ are Type-1 isomorphic.

\item [\rm (9)] $C_{24}(7(1,3,6,11,12))$ = $C_{24}(3,5,6,7,12)$. 

$\Rightarrow$ $C_{24}(1,3,6,11,12)$ and $C_{24}(3,5,6,7,12)$ are Type-1 isomorphic. 

\item [\rm (10)] $C_{24}(5(1,3,9,11,12))$ = $C_{24}(3,5,7,9,12)$.

$\Rightarrow$ $C_{24}(1,3,9,11,12)$ and $C_{24}(3,5,7,9,12)$ are Type-1 isomorphic.  

\item [\rm (11)] $C_{24}(7(1,6,9,11,12))$ = $C_{24}(5,6,7,9,12)$.

$\Rightarrow$ $C_{24}(1,3,9,11,12)$ and $C_{24}(3,5,7,9,12)$ are Type-1 isomorphic.  

\item [\rm (12)] $C_{24}(5(1,3,6,9,11,12))$ = $C_{24}(3,5,6,7,9,12)$.

$\Rightarrow$   $C_{24}(1,3,6,9,11,12)$ and $C_{24}(3,5,6,7,9,12)$ are Type-1 isomorphic.  \hfill $\Box$
\end{enumerate}

In the next two problems, we obtain Type-2 isomorphic circulant graphs of order 24 and show that the number of pairs of Type-2 isomorphic circulant graphs of order 24 is 32 and each pair is of Type-2 isomorphic w.r.t. $m$ = 2. At first, we show that the two pairs of circulant graphs $(a1)$ $C_{24}(1,2,11)$, $C_{24}(2,5,7)$; and $(b1)$ $C_{24}(1,10,11)$, $C_{24}(5,7,10)$ are Type-2 isomorphic w.r.t. $m$ = 2. 
 
\begin{prm} \label{p4} {\rm Show that the following pairs of circulant graphs are Type-2 isomorphic w.r.t. $m$ = 2.
\begin{enumerate}
\item [\rm (a1)]  $C_{24}(1,2,11)$ and $C_{24}(2,5,7)$; and
\item [\rm (b1)]  $C_{24}(1,10,11)$ and $C_{24}(5,7,10)$.
\end{enumerate} }
\end{prm}
\noindent
{\bf Solution.}
\begin{enumerate}
\item [\rm (a1)]   Let $R$ = $\{1,2,11\}$, $S$ = $R \cup (24-R)$, $T$ = $\{2,5,7\}$ and $n$ = 24 = $3\times 2^3$. Let $r$ = $2\in R,T$. This implies, $m$ = 2 = $\gcd(24,2)$ = $\gcd(24, r)$ is a possible value for the existence of Type-2 isomorphism of $\theta_{24,m,t}(C_{24}(R))$. We have $\theta_{24,2,t}(s)$ = $s+2jt$ for any $s\in\mathbb{Z}_{24}$ where $s$ = $2q+j$, $j,m,q\in\mathbb{N}_0$, $0 \leq j \leq 1$ and $0 \leq t \leq \frac{24}{\gcd(24, 2)}-1$ = 11. Consider,  
\\
$Ad_{24}(C_{24}(1,2,11))$ = $\{\varphi_{24,x}(C_{24}(1,2,11)): x = 1,5,7,11, 13,17,19,23\}$ 

\hspace{2.5cm}  = $\{C_{24}(x(1,2,11)): x = 1,5,7,11, 13,17,19,23\}$ 

\hspace{2.5cm} = $\{C_{24}(1,2,11), C_{24}(5,7,10)\}$  = $\{C_{24}(x(1,2,11)): x = 1,5\}$;
\\
 $\theta_{24,2,3}(C_{24}(1,2,11))$ = $\theta_{24,2,3}(C_{24}(1,2,11,  13,22,23))$ 

\hspace{2.6cm} = $C_{24}(\theta_{24,2,3}(1,2,11,  13,22,23))$  = $C_{24}(7,2,17,  19,22,5)$ 

\hspace{2.6cm} = $C_{24}(2,5,7,  17,19,22)$ = $C_{24}(2,5,7)$. 

$\Rightarrow$ $C_{24}(1,2,11)$ $\cong$ $C_{24}(2,5,7)$ and $C_{24}(2,5,7) \notin Ad_{24}(C_{24}(1,2,11))$. 

$\Rightarrow$  $C_{24}(1,2,11)$ and $C_{24}(2,5,7)$ are Type-2 isomorphic w.r.t. $m$ = 2.  

Type-2 isomorphic circulant graphs $C_{24}(1,2,11)$ and $C_{24}(2,5,7)$ w.r.t. $m$ = 2 are given in Figures 4 and 5.  In Figure 6, circulant graph $\theta_{24,2,3}(C_{24}(1,2,11))$ = $C_{24}(2,5,7)$ is given. 
\end{enumerate}
\begin{figure}[ht]
	\centerline{\hspace{.5cm} \includegraphics[width=6in]{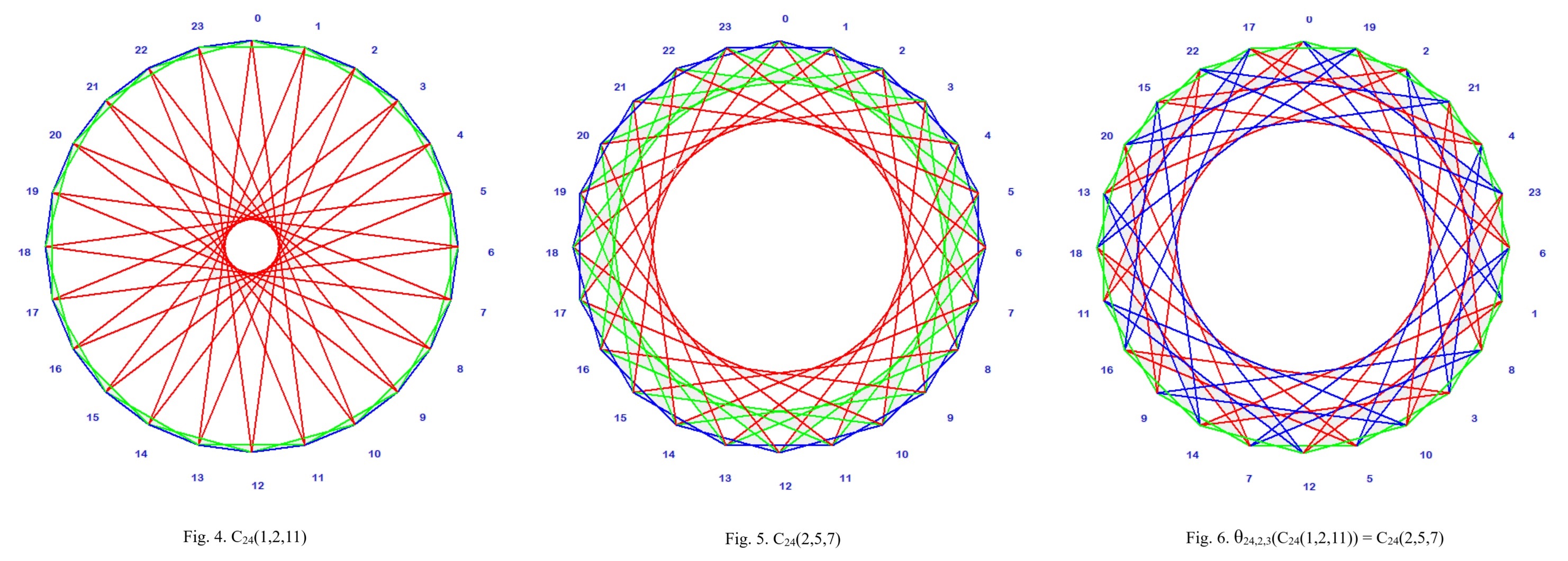}}
\end{figure}
\begin{enumerate}
\item [\rm (b1)]   Let $R$ = $\{1,10,11\}$, $S$ = $R \cup (24-R)$, $T$ = $\{5,7,10\}$ and $n$ = 24 = $3\times 2^3$. Let $r$ = $10\in R,T$. This implies, $m$ = 2 = $\gcd(24,10)$ = $\gcd(24, r)$ is a possible value for the existence of Type-2 isomorphism of $\theta_{24,m,t}(C_{24}(R))$. We have $\theta_{24,2,t}(s)$ = $s+2jt$ for any $s\in\mathbb{Z}_{24}$ where $s$ = $2q+j$, $j,m,q\in\mathbb{N}_0$, $0 \leq j \leq 1$ and $0 \leq t \leq \frac{24}{\gcd(24, 2)}-1$ = 11. Consider,  
\\
$Ad_{24}(C_{24}(1,10,11))$ = $\{C_{24}(1,10,11), C_{24}(2,5,7)\}$;
\\
 $\theta_{24,2,3}(C_{24}(1,10,11))$ = $\theta_{24,2,3}(C_{24}(1,10,11,  13,14,23))$ 

\hspace{2.8cm} = $C_{24}(\theta_{24,2,3}(1,10,11,  13,14,23))$  = $C_{24}(7,10,17,  19,14,5)$ 

\hspace{2.8cm} = $C_{24}(5,7,10,  14,17,19)$ = $C_{24}(5,7,10)$. 

$\Rightarrow$ $C_{24}(1,10,11)$ $\cong$ $C_{24}(5,7,10)$ and $C_{24}(5,7,10) \notin Ad_{24}(C_{24}(1,10,11))$. 

$\Rightarrow$ $C_{24}(1,10,11)$ and $C_{24}(5,7,10)$ are Type-2 isomorphic w.r.t. $m$ = 2.  \hfill $\Box$
\end{enumerate}

In the next problem, using remark \ref{r8} on the two pairs of Type-2 isomorphic circulant graphs $(a1)$ $C_{24}(1,2,11)$, $C_{24}(2,5,7)$; and $(b1)$ $C_{24}(1,10,11)$, $C_{24}(5,7,10)$, we find more pairs of isomorphic circulant graphs and all Type-2 isomorphic circulant graphs of order 24.

\begin{prm}\quad \label{p5} {\rm Show that each pair of circulant graphs of order 24 given below are isomorphic; classify their type of isomorphism and show that among them 32 pairs are Type-2 isomorphic circulant graphs.
\begin{enumerate}
\item [\rm (1)] $C_{24}(1,2,11)$, $C_{24}(2,5,7)$;  

\item [\rm (2)] $C_{24}(1,4,11)$, $C_{24}(4,5,7)$;  
	
\item [\rm (3)] $C_{24}(1,6,11)$, $C_{24}(5,6,7)$;  
	
\item [\rm (4)] $C_{24}(1,8,11)$, $C_{24}(5,7,8)$; 
	
\item [\rm (5)] $C_{24}(1,10,11)$, $C_{24}(5,7,10)$; 

\item [\rm (6)] $C_{24}(1,11,12)$, $C_{24}(5,7,12)$;  

\item [\rm (7)] $C_{24}(1,2,4,11)$, $C_{24}(2,4,5,7)$;  
	
\item [\rm (8)] $C_{24}(1,2,6,11)$, $C_{24}(2,5,6,7)$;
	
\item [\rm (9)] $C_{24}(1,2,8,11)$, $C_{24}(2,5,7,8)$; 
	
\item [\rm (10)] $C_{24}(1,2,10,11)$, $C_{24}(2,5,7,10)$; 
	
\item [\rm (11)] $C_{24}(1,2,11,12)$, $C_{24}(2,5,7,12)$; 
	
\item [\rm (12)] $C_{24}(1,4,6,11)$, $C_{24}(4,5,6,7)$; 
	
\item [\rm (13)] $C_{24}(1,4,8,11)$, $C_{24}(4,5,7,8)$;  

\item [\rm (14)] $C_{24}(1,4,10,11)$, $C_{24}(4,5,7,10)$;  

\item [\rm (15)] $C_{24}(1,4,11,12)$, $C_{24}(4,5,7,12)$;  

\item [\rm (16)] $C_{24}(1,6,8,11)$, $C_{24}(5,6,7,8)$;  

\item [\rm (17)] $C_{24}(1,6,10,11)$, $C_{24}(5,6,7,10)$; 
	
\item [\rm (18)] $C_{24}(1,6,11,12)$, $C_{24}(5,6,7,12)$;  

\item [\rm (19)] $C_{24}(1,8,10,11)$, $C_{24}(5,7,8,10)$;

\item [\rm (20)] $C_{24}(1,8,11,12)$, $C_{24}(5,7,8,12)$;  

\item [\rm (21)] $C_{24}(1,10,11,12)$, $C_{24}(5,7,10,12)$;
 
\item [\rm (22)] $C_{24}(1,2,4,6,11)$, $C_{24}(2,4,5,6,7)$; 

\item [\rm (23)] $C_{24}(1,2,4,8,11)$, $C_{24}(2,4,5,7,8)$;

\item [\rm (24)] $C_{24}(1,2,4,10,11)$, $C_{24}(2,4,5,7,10)$; 

\item [\rm (25)] $C_{24}(1,2,4,11,12)$, $C_{24}(2,4,5,7,12)$; 

\item [\rm (26)] $C_{24}(1,2,6,8,11)$, $C_{24}(2,5,6,7,8)$;

\item [\rm (27)] $C_{24}(1,2,6,10,11)$, $C_{24}(2,5,6,7,10)$; 

\item [\rm (28)] $C_{24}(1,2,6,11,12)$, $C_{24}(2,5,6,7,12)$;

\item [\rm (29)] $C_{24}(1,2,8,10,11)$, $C_{24}(2,5,7,8,10)$;

\item [\rm (30)] $C_{24}(1,2,8,11,12)$, $C_{24}(2,5,7,8,12)$;

\item [\rm (31)] $C_{24}(1,2,10,11,12)$, $C_{24}(2,5,7,10,12)$;

\item [\rm (32)] $C_{24}(1,4,6,8,11)$, $C_{24}(4,5,6,7,8)$;

\item [\rm (33)] $C_{24}(1,4,6,10,11)$, $C_{24}(4,5,6,7,10)$;

\item [\rm (34)] $C_{24}(1,4,6,11,12)$, $C_{24}(4,5,6,7,12)$;

\item [\rm (35)] $C_{24}(1,4,8,10,11)$, $C_{24}(4,5,7,8,10)$;

\item [\rm (36)] $C_{24}(1,4,8,11,12)$, $C_{24}(4,5,7,8,12)$;

\item [\rm (37)] $C_{24}(1,4,10,11,12)$, $C_{24}(4,5,7,10,12)$;

\item [\rm (38)] $C_{24}(1,6,8,10,11)$, $C_{24}(5,6,7,8,10)$;

\item [\rm (39)] $C_{24}(1,6,8,11,12)$, $C_{24}(5,6,7,8,12)$;

\item [\rm (40)] $C_{24}(1,6,10,11,12)$, $C_{24}(5,6,7,10,12)$;

\item [\rm (41)] $C_{24}(1,8,10,11,12)$, $C_{24}(5,7,8,10,12)$;

\item [\rm (42)] $C_{24}(1,2,4,6,8,11)$, $C_{24}(2,4,5,6,7,8)$;

\item [\rm (43)] $C_{24}(1,2,4,6,10,11)$, $C_{24}(2,4,5,6,7,10)$; 

\item [\rm (44)] $C_{24}(1,2,4,6,11,12)$, $C_{24}(2,4,5,6,7,12)$; 

\item [\rm (45)] $C_{24}(1,2,4,8,10,11)$, $C_{24}(2,4,5,7,8,10)$;

\item [\rm (46)] $C_{24}(1,2,4,8,11,12)$, $C_{24}(2,4,5,7,8,12)$;

\item [\rm (47)] $C_{24}(1,2,4,10,11,12)$, $C_{24}(2,4,5,7,10,12)$; 

\item [\rm (48)] $C_{24}(1,2,6,8,10,11)$, $C_{24}(2,5,6,7,8,10)$; 

\item [\rm (49)] $C_{24}(1,2,6,8,11,12)$, $C_{24}(2,6,5,7,8,12)$; 

\item [\rm (50)] $C_{24}(1,2,6,10,11,12)$, $C_{24}(2,5,6,7,10,12)$; 

\item [\rm (51)] $C_{24}(1,2,8,10,11,12)$, $C_{24}(2,5,8,7,10,12)$; 

\item [\rm (52)] $C_{24}(1,4,6,8,10,11)$, $C_{24}(4,5,6,7,8,10)$;

\item [\rm (53)] $C_{24}(1,4,6,8,11,12)$, $C_{24}(4,5,6,7,8,12)$;

\item [\rm (54)] $C_{24}(1,4,6,10,11,12)$, $C_{24}(4,5,6,7,10,12)$; 

\item [\rm (55)] $C_{24}(1,4,8,10,11,12)$, $C_{24}(4,5,7,8,10,12)$; 

\item [\rm (56)] $C_{24}(1,6,8,10,11,12)$, $C_{24}(5,6,7,8,10,12)$;

\item [\rm (57)] $C_{24}(1,2,4,6,8,10,11)$, $C_{24}(2,4,5,6,7,8,10)$; 

\item [\rm (58)] $C_{24}(1,2,4,6,8,11,12)$, $C_{24}(2,4,5,6,7,8,12)$; 

\item [\rm (59)] $C_{24}(1,2,4,6,10,11,12)$, $C_{24}(2,4,5,6,7,10,12)$; 

\item [\rm (60)] $C_{24}(1,2,4,8,10,11,12)$, $C_{24}(2,4,5,7,8,10,12)$; 

\item [\rm (61)] $C_{24}(1,2,6,8,10,11,12)$, $C_{24}(2,5,6,7,8,10,12)$; 

\item [\rm (62)] $C_{24}(1,4,6,8,10,11,12)$, $C_{24}(4,5,6,7,8,10,12)$; 

\item [\rm (63)] $C_{24}(1,2,4,6,8,10,11,12)$, $C_{24}(2,4,5,6,7,8,10,12)$. 
\end{enumerate}   }
\end{prm}
\noindent
{\bf Solution.}\quad Here, $n$ = 24 = $3\times 2^3$ and so the possible values of $m > 1$ $\ni$ $m$ is a divisor of $\gcd(n, r)$ = $\gcd(24, r)$, $m^3$ is a divisor of $n$ and $r\in R$ for the existence of isomorphic circulant graphs $C_{24}(R)$ of Type-2 w.r.t. $m$ is $m$ = 2 and 2 = $\gcd(24, 2)$ = $\gcd(24, 10)$ and multiple of 2 are 4 = $\gcd(24, 4)$, 6 = $\gcd(24, 6)$, 8 = $\gcd(24, 8)$ and 12 = $\gcd(24, 12)$. Thus, pairs of circulant graphs under cases (3), (4) and (19) can not be of Type-2 isomorphic. 

In problem \ref{p4}, we proved that the following two pairs of circulant graphs 
\\
$(1)$ $C_{24}(1, 2, 11)$, $C_{24}(2, 5, 7)$; and 
\\
$(5)$ $C_{24}(1, 10, 11)$, $C_{24}(5, 7, 10)$  are isomorphic of Type-2 w.r.t. $m$ = 2. 

Using remark \ref{r8} in the two pairs of Type-2 isomorphic circulant graphs, we obtain more isomorphic circulant graphs and all Type-2 isomorphic circulant graphs of order 24. We use remark \ref{r9} to establish isomorphism among circulant graphs of each pair given in the problem and remark \ref{r7} to establish Type-2 isomorphism among circulant graphs of each pair. 

\begin{enumerate}
\item [\rm (2)] $C_{24}(5(1,4,11))$ = $C_{24}(4,5,7)$. 
$\Rightarrow$  $C_{24}(1,4,11)$ and $C_{24}(4,5,7)$   are Type-1 isomorphic.
	
\item [\rm (3)] $C_{24}(5(1,6,11))$ = $C_{24}(5,6,7)$.
$\Rightarrow$  $C_{24}(1,6,11)$ and $C_{24}(5,6,7)$ are Type-1 isomorphic.
	
\item [\rm (4)] $C_{24}(5(1,8,11))$ = $C_{24}(5,7,8)$.
$\Rightarrow$  $C_{24}(1,8,11)$ and $C_{24}(5,7,8)$ are Type-1 isomorphic.
	
\item [\rm (6)] $C_{24}(5(1,11,12))$ = $C_{24}(5,7,12)$.
$\Rightarrow$  $C_{24}(1,11,12)$ and $C_{24}(5,7,12)$  are Type-1 isomorphic. 

\item [\rm (7)] $\theta_{24,2,3}(C_{24}(1,2,4,11))$ = $C_{24}(2,4,5,7)$. $\Rightarrow$  $C_{24}(1,2,4,11)$ $\cong$ $C_{24}(2,4,5,7)$. 
\\
$Ad_{24}(C_{24}(1,2,4,11))$ = $\{C_{24}(1,2,4,11), C_{24}(4,5,7,10) = C_{24}(5(1,2,4,11))\}$. 

This implies that  $C_{24}(2,4,5,7) \notin Ad_{24}(C_{24}(1,2,4,11))$. 

Hence, $C_{24}(1,2,4,11)$ and $C_{24}(2,4,5,7)$ are Type-2 isomorphic w.r.t. $m$ = 2.  
	
\item [\rm (8)] $\theta_{24,2,3}(C_{24}(1,2,6,11))$ = $C_{24}(2,5,6,7)$. $\Rightarrow$ $C_{24}(1,2,6,11)$ $\cong$ $C_{24}(2,5,6,7)$.   
\\
$Ad_{24}(C_{24}(1,2,6,11))$ = $\{C_{24}(1,2,6,11), C_{24}(5,6,7,10) = C_{24}(5(1,2,6,11))\}$. 

This implies that $C_{24}(2,5,6,7) \notin Ad_{24}(C_{24}(1,2,6,11))$. 

Hence, $C_{24}(1,2,6,11)$ and $C_{24}(2,5,6,7)$ are Type-2 isomorphic w.r.t. $m$ = 2.  
		
\item [\rm (9)] $\theta_{24,2,3}(C_{24}(1,2,8,11))$ = $C_{24}(2,5,7,8)$. $\Rightarrow$ $C_{24}(1,2,8,11)$ $\cong$ $C_{24}(2,5,7,8)$.   
\\
$Ad_{24}(C_{24}(1,2,8,11))$ = $\{C_{24}(1,2,8,11), C_{24}(5,7,8,10) = C_{24}(5(1,2,8,11))\}$. 

This implies, $C_{24}(2,5,7,8) \notin Ad_{24}(C_{24}(1,2,8,11))$. 

Hence, $C_{24}(1,2,8,11)$ and $C_{24}(2,5,7,8)$ are Type-2 isomorphic w.r.t. $m$ = 2.  
		
\item [\rm (10)] $C_{24}(5(1,2,10,11))$ = $C_{24}(2,5,7,10)$.

$\Rightarrow$   $C_{24}(1,2,10,11)$ and $C_{24}(2,5,7,10)$  are Type-1 isomorphic.
	
\item [\rm (11)] $\theta_{24,2,3}(C_{24}(1,2,11,12))$ = $C_{24}(2,5,7,12)$. 
$\Rightarrow$ $C_{24}(1,2,11,12)$ $\cong$ $C_{24}(2,5,7,12)$.    
\\
$Ad_{24}(C_{24}(1,2,11,12))$ = $\{C_{24}(1,2,11,12), C_{24}(5,7,10,12) = C_{24}(5(1,2,11,12))\}$. 

This implies, $C_{24}(2,5,7,12) \notin Ad_{24}(C_{24}(1,2,11,12))$. 

Hence, $C_{24}(1,2,11,12)$ and $C_{24}(2,5,7,12)$ are Type-2 isomorphic w.r.t. $m$ = 2.  
			
\item [\rm (12)] $C_{24}(5(1,4,6,11))$ = $C_{24}(4,5,6,7)$. 

$\Rightarrow$  $C_{24}(1,4,6,11)$ and $C_{24}(4,5,6,7)$ are Type-1 isomorphic.  

\item [\rm (13)] $C_{24}(5(1,4,8,11))$ = $C_{24}(4,5,7,8)$. 

$\Rightarrow$  $C_{24}(1,4,8,11)$ and $C_{24}(4,5,7,8)$ are Type-1 isomorphic.  

\item [\rm (14)] $ \theta_{24,2,3}(C_{24}(1,4,10,11))$ = $C_{24}(4,5,7,10)$. 
$\Rightarrow$    $C_{24}(1,4,10,11)$ $\cong$ $C_{24}(4,5,7,10)$.
\\
$Ad_{24}(C_{24}(1,4,10,11))$ = $\{C_{24}(1,4,10,11), C_{24}(2,4,5,7) = C_{24}(5(1,4,10,11))\}$. 

This implies, $C_{24}(4,5,7,10) \notin Ad_{24}(C_{24}(1,4,10,11))$. 

Hence, $C_{24}(1,4,10,11)$ and $C_{24}(4,5,7,10)$ are Type-2 isomorphic w.r.t. $m$ = 2.  
	
\item [\rm (15)] $C_{24}(5(1,4,11,12))$ = $C_{24}(4,5,7,12)$. 

$\Rightarrow$ $C_{24}(1,4,11,12)$ and $C_{24}(4,5,7,12)$  are Type-1 isomorphic.   

\item [\rm (16)] $C_{24}(5(1,6,8,11))$ = $C_{24}(5,6,7,8)$. 

$\Rightarrow$   $C_{24}(1,6,8,11)$ and $C_{24}(5,6,7,8)$ are Type-1 isomorphic. 

\item [\rm (17)] $ \theta_{24,2,3}(C_{24}(1,6,10,11))$ = $C_{24}(5,6,7,10)$. 
$\Rightarrow$  $C_{24}(1,6,10,11)$ $\cong$ $C_{24}(5,6,7,10)$.  
\\
$Ad_{24}(C_{24}(1,6,10,11))$ = $\{C_{24}(1,6,10,11), C_{24}(2,5,6,7) = C_{24}(5(1,6,10,11))\}$. 

This implies, $C_{24}(5,6,7,10) \notin Ad_{24}(C_{24}(1,6,10,11))$. 

Hence, $C_{24}(1,6,10,11)$ and $C_{24}(5,6,7,10)$ are Type-2 isomorphic w.r.t. $m$ = 2.  
		
\item [\rm (18)] $C_{24}(5(1,6,11,12))$ = $C_{24}(5,6,7,12)$. 

$\Rightarrow$ $C_{24}(1,6,11,12)$ and $C_{24}(5,6,7,12)$  are Type-1 isomorphic.   

\item [\rm (19)] $ \theta_{24,2,3}(C_{24}(1,8,10,11))$ = $C_{24}(5,7,8,10)$. 
$\Rightarrow$    $C_{24}(1,8,10,11)$ $\cong$ $C_{24}(5,7,8,10)$.
\\
$Ad_{24}(C_{24}(1,8,10,11))$ = $\{C_{24}(1,8,10,11), C_{24}(2,5,7,8) = C_{24}(5(1,8,10,11))\}$. 

This implies, $C_{24}(5,7,8,10) \notin Ad_{24}(C_{24}(1,8,10,11))$. 

Hence, $C_{24}(1,8,10,11)$ and $C_{24}(5,7,8,10)$ are Type-2 isomorphic w.r.t. $m$ = 2.  

\item [\rm (20)] $C_{24}(5(1,8,11,12))$ = $C_{24}(5,7,8,12)$. 

$\Rightarrow$ $C_{24}(1,8,11,12)$ and $C_{24}(5,7,8,12)$  are Type-1 isomorphic.   

\item [\rm (21)] $ \theta_{24,2,3}(C_{24}(1,10,11,12))$ = $C_{24}(5,7,10,12)$. 
$\Rightarrow$    $C_{24}(1,10,11,12)$ $\cong$ $C_{24}(5,7,10,12)$.
\\
$Ad_{24}(C_{24}(1,10,11,12))$ = $\{C_{24}(1,10,11,12), C_{24}(2,5,7,12) = C_{24(5}(1,10,11,12))\}$. 

This implies, $C_{24}(5,7,10,12) \notin Ad_{24}(C_{24}(1,10,11,12))$. 

Hence, $C_{24}(1,10,11,12)$ and $C_{24}(5,7,10,12)$ are Type-2 isomorphic w.r.t. $m$ = 2.  

\item [\rm (22)] $\theta_{24,2,3}(C_{24}(1,2,4,6,11))$ = $C_{24}(2,4,5,6,7)$. 
$\Rightarrow$ $C_{24}(1,2,4,6,11)$ $\cong$ $C_{24}(2,4,5,6,7)$.  
\\
$Ad_{24}(C_{24}(1,2,4,6,11))$ = $\{C_{24}(1,2,4,6,11), C_{24}(4,5,6,7,10) = C_{24}(5(1,2,4,6,11))\}$. 

This implies, $C_{24}(2,4,5,6,7) \notin Ad_{24}(C_{24}(1,2,4,6,11))$. 

Hence, $C_{24}(1,2,4,6,11)$ and $C_{24}(2,4,5,6,7)$ are Type-2 isomorphic w.r.t. $m$ = 2.  

\item [\rm (23)] $\theta_{24,2,3}(C_{24}(1,2,4,8,11))$ = $C_{24}(2,4,5,7,8)$. 
$\Rightarrow$   $C_{24}(1,2,4,8,11)$ $\cong$ $C_{24}(2,4,5,7,8)$.
\\
$Ad_{24}(C_{24}(1,2,4,8,11))$ = $\{C_{24}(1,2,4,8,11), C_{24}(4,5,8,7,10) = C_{24}(5(1,2,4,8,11))\}$. 

This implies, $C_{24}(2,4,5,7,8) \notin Ad_{24}(C_{24}(1,2,4,8,11))$. 

Hence, $C_{24}(1,2,4,8,11)$ and $C_{24}(2,4,5,7,8)$ are Type-2 isomorphic w.r.t. $m$ = 2.  

\item [\rm (24)] $C_{24}(5(1,2,4,10,11))$ = $C_{24}(2,4,5,7,10)$.

$\Rightarrow$  $C_{24}(1,2,4,10,11)$ and $C_{24}(2,4,5,7,10)$ are Type-1 isomorphic. 

\item [\rm (25)] $\theta_{24,2,3}(C_{24}(1,2,4,11,12))$ = $C_{24}(2,4,5,7,12)$. 
$\Rightarrow$   $C_{24}(1,2,4,11,12)$ $\cong$ $C_{24}(2,4,5,7,12)$.
\\
$Ad_{24}(C_{24}(1,2,4,11,12))$ = $\{C_{24}(1,2,4,11,12), C_{24}(4,5,7,10,12) = C_{24}(5(1,2,4,11,12))\}$. 

This implies, $C_{24}(2,4,5,7,12) \notin Ad_{24}(C_{24}(1,2,4,11,12))$. 

Hence, $C_{24}(1,2,4,11,12)$ and $C_{24}(2,4,5,7,12)$ are Type-2 isomorphic w.r.t. $m$ = 2.  

\item [\rm (26)] $\theta_{24,2,3}(C_{24}(1,2,6,8,11))$ = $C_{24}(2,5,6,7,8)$. 
$\Rightarrow$ $C_{24}(1,2,6,8,11)$ $\cong$ $C_{24}(2,5,6,7,8)$.  
\\
$Ad_{24}(C_{24}(1,2,6,8,11))$ = $\{C_{24}(1,2,6,8,11), C_{24}(5,6,7,8,10) = C_{24(5}(1,2,6,8,11))\}$. 

This implies, $C_{24}(2,5,6,7,8) \notin Ad_{24}(C_{24}(1,2,6,8,11))$. 

Hence, $C_{24}(1,2,6,8,11)$ and $C_{24}(2,5,6,7,8)$ are Type-2 isomorphic w.r.t. $m$ = 2.  

\item [\rm (27)] $C_{24}(5(1,2,6,10,11))$ = $C_{24}(2,5,6,7,10)$. 

$\Rightarrow$  $C_{24}(1,2,6,10,11)$ and $C_{24}(2,5,6,7,10)$ are Type-1 isomorphic. 

\item [\rm (28)] $\theta_{24,2,3}(C_{24}(1,2,6,11,12))$ = $C_{24}(2,5,6,7,12)$. 
$\Rightarrow$ $C_{24}(1,2,6,11,12)$ $\cong$ $C_{24}(2,5,6,7,12)$.  
\\
$Ad_{24}(C_{24}(1,2,6,11,12))$ = $\{C_{24}(1,2,6,11,12), C_{24}(5,6,7,10,12) = C_{24}(5(1,2,6,11,12))\}$. 

This implies, $C_{24}(2,5,6,7,12) \notin Ad_{24}(C_{24}(1,2,6,11,12))$. 

Hence, $C_{24}(1,2,6,11,12)$ and $C_{24}(2,5,6,7,12)$ are Type-2 isomorphic w.r.t. $m$ = 2.  

\item [\rm (29)] $C_{24}(5(1,2,8,10,11))$ = $C_{24}(2,5,7,8,10)$. 

$\Rightarrow$ $C_{24}(1,2,8,10,11)$ and $C_{24}(2,5,7,8,10)$ are Type-1 isomorphic.  

\item [\rm (30)] $\theta_{24,2,3}(C_{24}(1,2,8,11,12))$ = $C_{24}(2,5,7,8,12)$. 
$\Rightarrow$  $C_{24}(1,2,8,11,12)$ $\cong$ $C_{24}(2,5,7,8,12)$. 
\\
$Ad_{24}(C_{24}(1,2,8,11,12))$ = $\{C_{24}(1,2,8,11,12), C_{24}(5,7,8,10,12) = C_{24}(5(1,2,8,11,12))\}$. 

This implies, $C_{24}(2,5,7,8,12) \notin Ad_{24}(C_{24}(1,2,8,11,12))$. 

Hence, $C_{24}(1,2,8,11,12)$ and $C_{24}(2,5,7,8,12)$ are Type-2 isomorphic w.r.t. $m$ = 2.  

\item [\rm (31)] $C_{24}(5(1,2,10,11,12))$ = $C_{24}(2,5,7,10,12)$. 

$\Rightarrow$  $C_{24}(1,2,10,11,12)$ and $C_{24}(2,5,7,10,12)$ are Type-1 isomorphic. 

\item [\rm (32)] $C_{24}(5(1,4,6,8,11))$ = $C_{24}(4,5,6,7,8)$. 

$\Rightarrow$   $C_{24}(1,4,6,8,11)$ and $C_{24}(4,5,6,7,8)$ are Type-1 isomorphic. 

\item [\rm (33)] $ \theta_{24,2,3}(C_{24}(1,4,6,10,11))$ = $C_{24}(4,5,6,7,10)$. 
$\Rightarrow$    $C_{24}(1,4,6,10,11)$ $\cong$ $C_{24}(4,5,6,7,10)$.
\\
$Ad_{24}(C_{24}(1,4,6,10,11))$ = $\{C_{24}(1,4,6,10,11), C_{24}(2,4,5,6,7) = C_{24}(5(1,4,6,10,11))\}$. 

This implies, $C_{24}(4,5,6,7,10) \notin Ad_{24}(C_{24}(1,4,6,10,11))$. 

Hence, $C_{24}(1,4,6,10,11)$ and $C_{24}(4,5,6,7,10)$ are Type-2 isomorphic w.r.t. $m$ = 2.  

\item [\rm (34)] $C_{24}(5(1,4,6,11,12))$ = $C_{24}(4,5,6,7,12)$. 

$\Rightarrow$   $C_{24}(1,4,6,11,12)$ and $C_{24}(4,5,6,7,12)$ are Type-1 isomorphic. 

\item [\rm (35)] $ \theta_{24,2,3}(C_{24}(1,4,8,10,11))$ = $C_{24}(4,5,7,8,10)$. 
$\Rightarrow$    $C_{24}(1,4,8,10,11)$ $\cong$ $C_{24}(4,5,7,8,10)$.
\\
$Ad_{24}(C_{24}(1,4,8,10,11))$ = $\{C_{24}(1,4,8,10,11), C_{24}(2,4,5,7,8) = C_{24}(5(1,4,8,10,11))\}$. 

This implies, $C_{24}(4,5,7,8,10) \notin Ad_{24}(C_{24}(1,4,8,10,11))$. 

Hence, $C_{24}(1,4,8,10,11)$ and $C_{24}(4,5,7,8,10)$ are Type-2 isomorphic w.r.t. $m$ = 2.  

\item [\rm (36)] $C_{24}(5(1,4,8,11,12))$ = $C_{24}(4,5,7,8,12)$. 

$\Rightarrow$   $C_{24}(1,4,8,11,12)$ and $C_{24}(4,5,7,8,12)$ are Type-1 isomorphic. 

\item [\rm (37)] $ \theta_{24,2,3}(C_{24}(1,4,10,11,12))$ = $C_{24}(4,5,7,10,12)$. 
$\Rightarrow$    $C_{24}(1,4,10,11,12)$ $\cong$ $C_{24}(4,5,7,10,12)$
\\
$Ad_{24}(C_{24}(1,4,10,11,12))$ = $\{C_{24}(1,4,10,11,12), C_{24}(2,4,5,7,12) = C_{24}(5(1,4,10,11,12))\}$. 

This implies, $C_{24}(4,5,7,10,12) \notin Ad_{24}(C_{24}(1,4,10,11,12))$. 

Hence, $C_{24}(1,4,10,11,12)$ and $C_{24}(4,5,7,10,12)$ are Type-2 isomorphic w.r.t. $m$ = 2.  

\item [\rm (38)] $ \theta_{24,2,3}(C_{24}(1,6,8,10,11))$ = $C_{24}(5,6,7,8,10)$. 
$\Rightarrow$    $C_{24}(1,6,8,10,11)$ $\cong$ $C_{24}(5,6,7,8,10)$.
\\
$Ad_{24}(C_{24}(1,6,8,10,11))$ = $\{C_{24}(1,6,8,10,11), C_{24}(2,6,5,7,8) = C_{24}(5(1,6,8,10,11))\}$. 

This implies, $C_{24}(5,6,7,8,10) \notin Ad_{24}(C_{24}(1,6,8,10,11))$. 

Hence, $C_{24}(1,6,8,10,11)$ and $C_{24}(5,6,7,8,10)$ are Type-2 isomorphic w.r.t. $m$ = 2.  

\item [\rm (39)] $C_{24}(5(1,6,8,11,12))$ = $C_{24}(5,6,7,8,12)$. 

$\Rightarrow$   $C_{24}(1,6,8,11,12)$ and $C_{24}(5,6,7,8,12)$ are Type-1 isomorphic. 

\item [\rm (40)] $ \theta_{24,2,3}(C_{24}(1,6,10,11,12))$ = $C_{24}(5,6,7,10,12)$. 
$\Rightarrow$  $C_{24}(1,6,10,11,12)$ $\cong$ $C_{24}(5,6,7,10,12)$. 
\\
$Ad_{24}(C_{24}(1,6,10,11,12))$ = $\{C_{24}(1,6,10,11,12), C_{24}(2,5,6,7,12) = C_{24(5}(1,6,10,11,12))\}$. 

This implies, $C_{24}(5,6,7,10,12) \notin Ad_{24}(C_{24}(1,6,10,11,12))$. 

Hence, $C_{24}(1,6,10,11,12)$ and $C_{24}(5,6,7,10,12)$ are Type-2 isomorphic w.r.t. $m$ = 2.  

\item [\rm (41)] $ \theta_{24,2,3}(C_{24}(1,8,10,11,12))$ = $C_{24}(5,7,8,10,12)$. 
$\Rightarrow$    $C_{24}(1,8,10,11,12)$ $\cong$ $C_{24}(5,7,8,10,12)$.
\\
$Ad_{24}(C_{24}(1,8,10,11,12))$ = $\{C_{24}(1,8,10,11,12), C_{24}(2,5,7,8,12) = C_{24}(5(1,8,10,11,12))\}$. 

This implies, $C_{24}(5,7,8,10,12) \notin Ad_{24}(C_{24}(1,8,10,11,12))$. 

Hence, $C_{24}(1,8,10,11,12)$ and $C_{24}(5,7,8,10,12)$ are Type-2 isomorphic w.r.t. $m$ = 2.  

\item [\rm (42)] $\theta_{24,2,3}(C_{24}(1,2,4,6,8,11))$ = $C_{24}(2,4,5,6,7,8)$. 
$\Rightarrow$   $C_{24}(1,2,4,6,8,11)$ $\cong$ $C_{24}(2,4,5,6,7,8)$.
\\
$Ad_{24}(C_{24}(1,2,4,6,8,11))$ = $\{C_{24}(1,2,4,6,8,11), C_{24}(4,5,6,7,8,10) = C_{24}(5(1,2,4,6,8,11))\}$. 

This implies, $C_{24}(2,4,5,6,7,8) \notin Ad_{24}(C_{24}(1,2,4,6,8,11))$. 

Hence, $C_{24}(1,2,4,6,8,11)$ and $C_{24}(2,4,5,6,7,8)$ are Type-2 isomorphic w.r.t. $m$ = 2.  

\item [\rm (43)] $C_{24}(5(1,2,4,6,10,11))$ = $C_{24}(2,4,5,6,7,10)$.

$\Rightarrow$   $C_{24}(1,2,4,6,10,11)$ and $C_{24}(2,4,5,6,7,10)$ are Type-1 isomorphic. 

\item [\rm (44)] $\theta_{24,2,3}(C_{24}(1,2,4,6,11,12))$ = $C_{24}(2,4,5,6,7,12)$. 
$\Rightarrow$   $C_{24}(1,2,4,6,11,12)$ $\cong$ $C_{24}(2,4,5,6,7,12)$.
\\
$Ad_{24}(C_{24}(1,2,4,6,11,12))$ = $\{C_{24}(1,2,4,6,11,12)$, 

\hfill $C_{24}(4,5,6,7,10,12) = C_{24}(5(1,2,4,6,11,12))\}$. 

This implies, $C_{24}(2,4,5,6,7,12) \notin Ad_{24}(C_{24}(1,2,4,6,11,12))$. 

Hence, $C_{24}(1,2,4,6,11,12)$ and $C_{24}(2,4,5,6,7,12)$ are Type-2 isomorphic w.r.t. $m$ = 2.  

\item [\rm (45)] $C_{24}(5(1,2,4,8,10,11))$ = $C_{24}(2,4,5,7,8,10)$. 

$\Rightarrow$    $C_{24}(1,2,4,8,10,11)$ and $C_{24}(2,4,5,7,8,10)$ are Type-1 isomorphic.

\item [\rm (46)] $\theta_{24,2,3}(C_{24}(1,2,4,8,11,12))$ = $C_{24}(2,4,5,7,8,12)$. 
$\Rightarrow$   $C_{24}(1,2,4,8,11,12)$ $\cong$ $C_{24}(2,4,5,7,8,12)$.
\\
$Ad_{24}(C_{24}(1,2,4,8,11,12))$ = $\{C_{24}(1,2,4,8,11,12)$, 

\hfill $C_{24}(4,5,7,8,10,12) = C_{24}(5(1,2,4,8,11,12))\}$. 

This implies, $C_{24}(2,4,5,7,8,12) \notin Ad_{24}(C_{24}(1,2,4,8,11,12))$. 

Hence, $C_{24}(1,2,4,8,11,12)$ and $C_{24}(2,4,5,7,8,12)$ are Type-2 isomorphic w.r.t. $m$ = 2.  

\item [\rm (47)] $C_{24}(5(1,2,4,10,11,12))$ = $C_{24}(2,4,5,7,10,12)$. 

$\Rightarrow$ $C_{24}(1,2,4,10,11,12)$ and $C_{24}(2,4,5,7,10,12)$ are Type-1 isomorphic.   

\item [\rm (48)] $C_{24}(5(1,2,6,8,10,11))$ = $C_{24}(2,5,6,7,8,10)$. 

$\Rightarrow$  $C_{24}(1,2,6,8,10,11)$ and $C_{24}(2,5,6,7,8,10)$  are Type-1 isomorphic.  

\item [\rm (49)] $\theta_{24,2,3}(C_{24}(1,2,6,8,11,12))$ = $C_{24}(2,6,5,7,8,12)$. 
$\Rightarrow$  $C_{24}(1,2,6,8,11,12)$ $\cong$ $C_{24}(2,6,5,7,8,12)$. 
\\
$Ad_{24}(C_{24}(1,2,6,8,11,12))$ = $\{C_{24}(1,2,6,8,11,12)$, 

\hfill $C_{24}(5,6,7,8,10,12) = C_{24}(5(1,2,6,8,11,12))\}$. 

This implies, $C_{24}(2,5,6,7,8,12) \notin Ad_{24}(C_{24}(1,2,6,8,11,12))$. 

Hence, $C_{24}(1,2,6,8,11,12)$ and $C_{24}(2,5,6,7,8,12)$ are Type-2 isomorphic w.r.t. $m$ = 2.  

\item [\rm (50)] $C_{24}(5(1,2,6,10,11,12))$ = $C_{24}(2,5,6,7,10,12)$. 

$\Rightarrow$    $C_{24}(1,2,6,10,11,12)$ and $C_{24}(2,5,6,7,10,12)$ are Type-1 isomorphic.

\item [\rm (51)] $C_{24}(5(1,2,8,10,11,12))$ = $C_{24}(2,5,8,7,10,12)$. 

$\Rightarrow$   $C_{24}(1,2,8,10,11,12)$ and $C_{24}(2,5,8,7,10,12)$  are Type-1 isomorphic.

\item [\rm (52)] $ \theta_{24,2,3}(C_{24}(1,4,6,8,10,11))$ = $C_{24}(4,5,6,7,8,10)$. 
$\Rightarrow$    $C_{24}(1,4,6,8,10,11)$ $\cong$ $C_{24}(4,5,6,7,8,10)$.
\\
$Ad_{24}(C_{24}(1,4,6,8,10,11))$ = $\{C_{24}(1,4,6,8,10,11)$, 

\hfill $C_{24}(2,4,5,6,7,8) = C_{24}(5(1,4,6,8,10,11))\}$. 

This implies, $C_{24}(4,5,6,7,8,10) \notin Ad_{24}(C_{24}(1,4,6,8,10,11))$. 

Hence, $C_{24}(1,4,6,8,10,11)$ and $C_{24}(4,5,6,7,8,10)$ are Type-2 isomorphic w.r.t. $m$ = 2.  

\item [\rm (53)]  $C_{24}(5(1,4,6,8,11,12))$ = $C_{24}(4,5,6,7,8,12)$.

$\Rightarrow$   $C_{24}(1,4,6,8,11,12)$ and $C_{24}(4,5,6,7,8,12)$  are Type-1 isomorphic.

\item [\rm (54)] $ \theta_{24,2,3}(C_{24}(1,4,6,10,11,12))$ = $C_{24}(4,5,6,7,10,12)$.

$\Rightarrow$    $C_{24}(1,4,6,10,11,12)$ $\cong$ $C_{24}(4,5,6,7,10,12)$.
\\
$Ad_{24}(C_{24}(1,4,6,10,11,12))$ = $\{C_{24}(1,4,6,10,11,12)$, 

\hfill $C_{24}(2,4,5,6,7,12) = C_{24}(5(1,4,6,10,11,12))\}$.

 This implies, $C_{24}(4,5,6,7,10,12) \notin Ad_{24}(C_{24}(1,4,6,10,11,12))$. 

Hence, $C_{24}(1,4,6,10,11,12)$ and $C_{24}(4,5,6,7,10,12)$ are Type-2 isomorphic w.r.t. $m$ = 2.  

\item [\rm (55)] $ \theta_{24,2,3}(C_{24}(1,4,8,10,11,12))$ = $C_{24}(4,5,7,8,10,12)$.
 
$\Rightarrow$  $C_{24}(1,4,8,10,11,12)$ $\cong$ $C_{24}(4,5,7,8,10,12)$.  
\\
$Ad_{24}(C_{24}(1,4,8,10,11,12))$ = $\{C_{24}(1,4,8,10,11,12)$, 

\hfill $C_{24}(2,4,5,7,8,12) = C_{24}(5(1,4,8,10,11,12))\}$. 

This implies, $C_{24}(4,5,7,8,10,12) \notin Ad_{24}(C_{24}(1,4,8,10,11,12))$. 

Hence, $C_{24}(1,4,8,10,11,12)$ and $C_{24}(4,5,7,8,10,12)$ are Type-2 isomorphic w.r.t. $m$ = 2.  

\item [\rm (56)] $ \theta_{24,2,3}(C_{24}(1,6,8,10,11,12))$ = $C_{24}(5,6,7,8,10,12)$.
 
$\Rightarrow$   $C_{24}(1,6,8,10,11,12)$ $\cong$ $C_{24}(5,6,7,8,10,12)$.
\\
$Ad_{24}(C_{24}(1,6,8,10,11,12))$ = $\{C_{24}(1,6,8,10,11,12)$, 

\hfill $C_{24}(2,6,5,7,8,12) = C_{24}(5(1,6,8,10,11,12))\}$. 

This implies, $C_{24}(5,6,7,8,10,12) \notin Ad_{24}(C_{24}(1,6,8,10,11,12))$. 

Hence, $C_{24}(1,6,8,10,11,12)$ and $C_{24}(5,6,7,8,10,12)$ are Type-2 isomorphic w.r.t. $m$ = 2.  

\item [\rm (57)] $C_{24}(5(1,2,4,6,8,10,11))$ = $C_{24}(5(1,2,4,6,8,10,11, 13,14,16,18,20,22,23))$ 

\hfill = $C_{24}(5,10,20,6,16,2,7, 17,22,8,18,4,14,19)$ = $C_{24}(2,4,5,6,7,8,10)$. 

$\Rightarrow$   $C_{24}(1,2,4,6,8,10,11)$ and $C_{24}(2,4,5,6,7,8,10)$ are Type-1 isomorphic.
 
\item [\rm (58)] $\theta_{24,2,3}(C_{24}(1,2,4,6,8,11,12))$ = $\theta_{24,2,3}(C_{24}(1,2,4,6,8,11, 12, 13,16,18,20,22,23))$ 

\hfill = $C_{24}(7,2,4,6,8,17, 12, 19,16,18,20,22,5)$ = $C_{24}(2,4,5,6,7,8,12)$.

$\Rightarrow$ $C_{24}(1,2,4,6,8,11,12)$ $\cong$ $C_{24}(2,4,5,6,7,8,12)$. 
\\
$Ad_{24}(C_{24}(1,2,4,6,8,11,12))$ = $\{C_{24}(1,2,4,6,8,11,12)$, 

\hfill $C_{24}(4,5,6,7,8,10,12) = C_{24}(5(1,2,4,6,8,11,12))\}$. 

This implies, $C_{24}(2,4,5,6,7,8,12) \notin Ad_{24}(C_{24}(1,2,4,6,8,11,12))$. 

Hence, $C_{24}(1,2,4,6,8,11,12)$ and $C_{24}(2,4,5,6,7,8,12)$ are Type-2 isomorphic w.r.t. $m$ = 2.  

\item [\rm (59)] $C_{24}(5(1,2,4,6,10,11,12))$ = $C_{24}(5(1,2,4,6,10,11, 12, 13,14,18,20,22,23))$ 

\hfill = $C_{24}(5,10,20,6,2,7, 12, 17,22,18,4,14,19))$ = $C_{24}(2,4,5,6,7,10,12)$.

$\Rightarrow$   $C_{24}(1,2,4,6,10,11,12)$ and $C_{24}(2,4,5,6,7,10,12)$ are Type-1 isomorphic.

\item [\rm (60)] $C_{24}(5(1,2,4,8,10,11,12))$ = $C_{24}(5(1,2,4,8,10,11, 12, 13,14,16,20,22,23))$ 

\hfill = $C_{24}(5,10,20,16,2,7, 12, 17,22,8,4,14,19))$ = $C_{24}(2,4,5,7,8,10,12)$.

$\Rightarrow$   $C_{24}(1,2,4,8,10,11,12)$ and $C_{24}(2,4,5,7,8,10,12)$ are Type-1 isomorphic.

\item [\rm (61)] $C_{24}(5(1,2,6,8,10,11,12))$ = $C_{24}(5(1,2,6,8,10,11, 12, 13,14,16,18,22,23))$ 

\hfill = $C_{24}(5,10,6,16,2,7, 12, 17,22,8,18,14,19))$ = $C_{24}(2,5,6,7,8,10,12)$.

$\Rightarrow$   $C_{24}(1,2,6,8,10,11,12)$ and $C_{24}(2,5,6,7,8,10,12)$ are Type-1 isomorphic.

\item [\rm (62)] $\theta_{24,2,3}(C_{24}(1,4,6,8,10,11,12))$ = $\theta_{24,2,3}(C_{24}(1,4,6,8,10,11, 12, 13,14,16,18,20,23))$ 

\hfill = $C_{24}(7,4,6,8,10,17, 12, 19,14,16,18,20,5)$ = $C_{24}(4,5,6,7,8,10,12)$.

$\Rightarrow$ $C_{24}(1,4,6,8,10,11,12)$ $\cong$ $C_{24}(4,5,6,7,8,10,12)$. 
\\
$Ad_{24}(C_{24}(1,4,6,8,10,11,12))$ = $\{C_{24}(1,4,6,8,10,11,12)$, 

\hfill $C_{24}(2,4,5,6,7,8,12) = C_{24}(5(1,4,6,8,10,11,12))\}$. 
\\
This implies, $C_{24}(4,5,6,7,8,10,12) \notin Ad_{24}(C_{24}(1,4,6,8,10,11,12))$. 
\\
Hence, $C_{24}(1,4,6,8,10,11,12)$ and $C_{24}(4,5,6,7,8,10,12)$ are Type-2 isomorphic w.r.t. $m$ = 2.  

\item [\rm (63)] $C_{24}(5(1,2,4,6,8,10,11,12))$ = $C_{24}(5(1,2,4,6,8,10,11, 12, 13,14,16,18,20,22,23))$ 

\hfill = $C_{24}(5,10,20,6,16,2,7, 12, 17,22,8,18,4,14,19))$ = $C_{24}(2,4,5,6,7,8,10,12)$. 

$\Rightarrow$   $C_{24}(1,2,4,6,8,10,11,12)$ and $C_{24}(2,4,5,6,7,8,10,12)$ are Type-1 isomorphic.  
\end{enumerate}

Thus, out of the 63 pairs of isomorphic circulant graphs of order 24, exactly 32 pairs are Type-2 isomorphic and all these 32 pairs of circulant graphs are Type-2 isomorphic w.r.t. $m$ = 2.  All the 32 pairs of isomorphic circulant graphs of Type-2 w.r.t. $m$ = 2 are given below. 
\begin{enumerate}
\item [\rm (1)] $C_{24}(1,2,11)$, $C_{24}(2,5,7)$;  

\item [\rm (5)] $C_{24}(1,10,11)$, $C_{24}(5,7,10)$; 

\item [\rm (7)] $C_{24}(1,2,4,11)$, $C_{24}(2,4,5,7)$;  
	
\item [\rm (8)] $C_{24}(1,2,6,11)$, $C_{24}(2,5,6,7)$;
	
\item [\rm (9)] $C_{24}(1,2,8,11)$, $C_{24}(2,5,7,8)$; 
	
\item [\rm (11)] $C_{24}(1,2,11,12)$, $C_{24}(2,5,7,12)$; 
	
\item [\rm (14)] $C_{24}(1,4,10,11)$, $C_{24}(4,5,7,10)$; 
	
\item [\rm (17)] $C_{24}(1,6,10,11)$, $C_{24}(5,6,7,10)$; 
	
\item [\rm (19)] $C_{24}(1,8,10,11)$, $C_{24}(5,7,8,10)$;

\item [\rm (21)] $C_{24}(1,10,11,12)$, $C_{24}(5,7,10,12)$;

\item [\rm (22)] $C_{24}(1,2,4,6,11)$, $C_{24}(2,4,5,6,7)$; 

\item [\rm (23)] $C_{24}(1,2,4,8,11)$, $C_{24}(2,4,5,7,8)$;

\item [\rm (25)] $C_{24}(1,2,4,11,12)$, $C_{24}(2,4,5,7,12)$; 

\item [\rm (26)] $C_{24}(1,2,6,8,11)$, $C_{24}(2,5,6,7,8)$;

\item [\rm (28)] $C_{24}(1,2,6,11,12)$, $C_{24}(2,5,6,7,12)$;

\item [\rm (30)] $C_{24}(1,2,8,11,12)$, $C_{24}(2,5,7,8,12)$;

\item [\rm (33)] $C_{24}(1,4,6,10,11)$, $C_{24}(4,5,6,7,10)$;

\item [\rm (35)] $C_{24}(1,4,8,10,11)$, $C_{24}(4,5,7,8,10)$;

\item [\rm (37)] $C_{24}(1,4,10,11,12)$, $C_{24}(4,5,7,10,12)$;

\item [\rm (38)] $C_{24}(1,6,8,10,11)$, $C_{24}(5,6,7,8,10)$;

\item [\rm (40)] $C_{24}(1,6,10,11,12)$, $C_{24}(5,6,7,10,12)$;

\item [\rm (41)] $C_{24}(1,8,10,11,12)$, $C_{24}(5,7,8,10,12)$;

\item [\rm (42)] $C_{24}(1,2,4,6,8,11)$, $C_{24}(2,4,5,6,7,8)$;

\item [\rm (44)] $C_{24}(1,2,4,6,11,12)$, $C_{24}(2,4,5,6,7,12)$; 

\item [\rm (46)] $C_{24}(1,2,4,8,11,12)$, $C_{24}(2,4,5,7,8,12)$;

\item [\rm (49)] $C_{24}(1,2,6,8,11,12)$, $C_{24}(2,6,5,7,8,12)$; 

\item [\rm (52)] $C_{24}(1,4,6,8,10,11)$, $C_{24}(4,5,6,7,8,10)$;

\item [\rm (54)] $C_{24}(1,4,6,10,11,12)$, $C_{24}(4,5,6,7,10,12)$; 

\item [\rm (55)] $C_{24}(1,4,8,10,11,12)$, $C_{24}(4,5,7,8,10,12)$; 

\item [\rm (56)] $C_{24}(1,6,8,10,11,12)$, $C_{24}(5,6,7,8,10,12)$;

\item [\rm (58)] $C_{24}(1,2,4,6,8,11,12)$, $C_{24}(2,4,5,6,7,8,12)$; 

\item [\rm (62)] $C_{24}(1,4,6,8,10,11,12)$, $C_{24}(4,5,6,7,8,10,12)$. \hfill $\Box$
\end{enumerate}

Given a circulant graph $C_n(R)$ having isomorphic circulant graphs of Type-2 w.r.t. $m$, remark \ref{r8} helps us to obtain more  isomorphic graphs which covers Type-2 w.r.t. $m$ as well as some Adam's isomorphic graphs of $C_n(R)$. In this section, we obtain all the 12 triples of Type-2 isomorphic circulant graphs of order 27. We also find isomorphic circulant graphs $C_{48}(r_1,r_2,r_3)$ and $C_{96}(r_1,r_2,r_3)$ of Type-2 w.r.t. $m$ = 2 and $C_{54}(r_1,r_2,r_3,r_4)$ of Type-2 w.r.t. $m$ = 3 and give a few open problems on these type of isomorphic circulant graphs.

\subsection{On Type-2 isomorphic circulant graphs $C_{27}(R)$}

In this section, we find all isomorphic circulant graphs of Type-2 of order 27. In problem 3.16 in \cite{v2-1}, we proved that circulant graphs $C_{27}(1,3,8,10)$, $C_{27}(2,3,7,11)$ and $C_{27}(3,4,5,13)$ are isomorphic of Type-2 w.r.t. $m$ = 3. In the next problem, using remark \ref{r8} in the above triple Type-2 isomorphic circulant graphs, we find all Type-2 isomorphic circulant graphs of order 27.

 \begin{prm}\quad \label{p6} {\rm Show that there are 12 triples of Type-2 isomorphic circulant graphs of order 27 and they are all isomorphic of Type-2 w.r.t. $m$ = 3. }
 \end{prm}
 \noindent
 {\bf Solution.}\quad Here, $n$ = 27 = $3^3$. Therefore, possible value of $m > 1$ $\ni$ $m$ divides $\gcd(n, r)$ and $m^3$ divides $n$ is $m$ = 3. Also, $m$ = 3 = $\gcd(n, r)$ = $\gcd(27, r)$ = $\gcd(27, 3)$ = $\gcd(27, 6)$ = $\gcd(27, 12)$. We start with finding isomorphic circulant graphs of Type-2 w.r.t. $m$ = 3. In problem 3.16 in \cite{v2-1}, it is proved that circulant graphs 
 
 (1) $C_{27}(1,3,8,10)$, $C_{27}(2,3,7,11)$ and $C_{27}(3,4,5,13)$ are isomorphic of Type-2 w.r.t. $m$ = 3. 
 
 Using remark \ref{r8} in this triple of isomorphic circulant graphs of Type-2 w.r.t. $m$ = 3, we get the following triples of isomorphic circulant graphs of the form $C_{27}(R)$ and each of these triples are either Type-1 isomorphic circulant graphs or Type-2 isomorphic circulant graphs w.r.t. $m$ = 3. 	
 
 \begin{enumerate} 	
	
\item [\rm (2)]	$C_{27}(1,6,8,10)$, $C_{27}(2,6,7,11)$, $C_{27}(4,5,6,13)$; 
			
\item [\rm (3)]	$C_{27}(1,8,9,10)$, $C_{27}(2,7,9,11)$, $C_{27}(4,5,9,13)$; 
			
\item [\rm (4)]	$C_{27}(1,8,10,12)$, $C_{27}(2,7,11,12)$, $C_{27}(4,5,12,13)$;
		
\item [\rm (5)]	$C_{27}(1,3,6,8,10)$, $C_{27}(2,3,6,7,11)$, $C_{27}(3,4,5,6,13)$; 
		
\item [\rm (6)]	$C_{27}(1,3,8,9,10)$, $C_{27}(2,3,7,9,11)$, $C_{27}(3,4,5,9,13)$; 
		
\item [\rm (7)]	$C_{27}(1,3,8,10,12)$, $C_{27}(2,3,7,11,12)$, $C_{27}(3,4,5,12,13)$;
		
\item [\rm (8)]	$C_{27}(1,6,8,9,10)$, $C_{27}(2,6,7,9,11)$, $C_{27}(4,5,6,9,13)$;
		
\item [\rm (9)]	$C_{27}(1,6,8,10,12)$, $C_{27}(2,6,7,11,12)$, $C_{27}(4,5,6,12,13)$; 
		
\item [\rm (10)]	$C_{27}(1,8,9,10,12)$, $C_{27}(2,7,9,11,12)$, $C_{27}(4,5,9,12,13)$; 
		
\item [\rm (11)] $C_{27}(1,3,6,8,9,10)$, $C_{27}(2,3,6,7,9,11)$, $C_{27}(3,4,5,6,9,13)$;
		
\item [\rm (12)] $C_{27}(1,3,6,8,10,12)$, $C_{27}(2,3,6,7,11,12)$, $C_{27}(3,4,5,6,12,13)$;
		
\item [\rm (13)] $C_{27}(1,3,8,9,10,12)$, $C_{27}(2,3,7,9,11,12)$, $C_{27}(3,4,5,9,12,13)$; 
		
\item [\rm (14)] $C_{27}(1,6,8,9,10,12)$, $C_{27}(2,6,7,9,11,12)$, $C_{27}(4,5,6,9,12,13)$;

\item [\rm (15)] $C_{27}(1,3,6,8,9,10,12)$, $C_{27}(2,3,6,7,9,11,12)$, $C_{27}(3,4,5,6,9,12,13)$. 
\end{enumerate}	

And among them the following three triples are of Type-1 isomorphic since
\\
~(3)  	$C_{27}(2(1,8,9,10))$ = $C_{27}(2(1,8,9,10, 17,18,19,26))$ = $C_{27}(2,16,18,20, 7,9,11,25)$ = $C_{27}(2,7,9,11)$, 

~	$C_{27}(4(1,8,9,10))$ = $C_{27}(4(1,8,9,10, 17,18,19,26))$ = $C_{27}(4,5,9,13, 14,18,22,23)$ = $C_{27}(4,5,9,13)$;
\\
(12)  $C_{27}(2(1,3,6,8,10,12))$ = $C_{27}(2(1,3,6,8,10,12, 15,17,19,21,24,26))$ 

\hfill = $C_{27}(2,6,12,16,20,24, 3,7,11,15,21,25)$ = $C_{27}(2,3,6,7,11,12)$ and

  ~~ $C_{27}(4(1,3,6,8,10,12))$ = $C_{27}(4(1,3,6,8,10,12, 15,17,19,21,24,26))$ 

\hfill = $C_{27}(4,12,24,5,13,21, 6,14,22,3,15,23)$ = $C_{27}(3,4,5,6,12,13)$;	
\\	
(15) $C_{27}(2(1,3,6,8,9,10,12))$ = $C_{27}(2(1,3,6,8,9,10,12, 15,17,18,19,21,24,26))$ 

\hfill = $C_{27}(2,6,12,16,18,20,24, 3,7,9,11,15,21,25)$ = $C_{27}(2,3,6,7,9,11,12)$ and

 ~~ $C_{27}(4(1,3,6,8,9 10,12))$ = $C_{27}(4(1,3,6,8,9,10,12, 15,17,18,19,21,24,26))$ 

\hfill = $C_{27}(4,12,24,5,9,13,21, 6,14,18,22,3,15,23)$ = $C_{27}(3,4,5,6,9,12,13)$.	

And all others are of Type-2 isomorphic w.r.t. $m$ = 3 since 
	
\begin{enumerate}
\item [\rm (2)]	$\theta_{27,3,1}(C_{27}(1,6,8,10))$ = $C_{27}(4,5,6,13)$ and $\theta_{27,3,2}(C_{27}(1,6,8,10))$ = $C_{27}(2,6,7,11)$.
		
 $\Rightarrow$ $C_{27}(1,6,8,10)$ $\cong$ $C_{27}(4,5,6,13)$ $\cong$ $C_{27}(2,6,7,11)$. 
\\		
Also,	$Ad_{27}(C_{27}(1,6,8,10))$ = $\{C_{27}(1,6,8,10), C_{27}(2,7,11,12), C_{27}(3,4,5,13)\}$. 
		
 $\Rightarrow$  $C_{27}(2,6,7,11),C_{27}(4,5,6,13)\notin Ad_{27}(C_{27}(1,6,8,10))$.
	 
$\Rightarrow$ $C_{27}(1,6,8,10)$, $C_{27}(4,5,6,13)$,  $C_{27}(2,6,7,11)$ are isomorphic of Type-2 w.r.t. $m$ = 3.
		
\item [\rm (4)]	$\theta_{27,3,1}(C_{27}(1,8,10,12))$ = $C_{27}(4,5,12,13)$ and $\theta_{27,3,2}(C_{27}(1,8,10,12))$ = $C_{27}(2,7,11,12)$.
		
$\Rightarrow$ $C_{27}(1,8,10,12)$ $\cong$ $C_{27}(4,5,12,13)$ $\cong$ $C_{27}(2,7,11,12)$. 
\\		
Also,	$Ad_{27}(C_{27}(1,8,10,12))$ = $\{C_{27}(1,8,10,12), C_{27}(2,3,7,11), C_{27}(4,5,6,13)\}$. 
		
$\Rightarrow$  $C_{27}(2,7,11,12),C_{27}(4,5,12,13)\notin Ad_{27}(C_{27}(1,8,10,12))$.
		
$\Rightarrow$ $C_{27}(1,8,10,12)$, $C_{27}(4,5,12,13)$,  $C_{27}(2,7,11,12)$ are isomorphic of Type-2 w.r.t. $m$ = 3.
		
\item [\rm (5)]	$\theta_{27,3,1}(C_{27}(1,3,6,8,10))$ = $C_{27}(3,4,5,6,13)$ and $\theta_{27,3,2}(C_{27}(1,3,6,8,10))$ = $C_{27}(2,3,6,7,11)$.

$\Rightarrow$ $C_{27}(1,3,6,8,10)$ $\cong$ $C_{27}(3,4,5,6,13)$ $\cong$ $C_{27}(2,3,6,7,11)$. 
\\
Also,	$Ad_{27}(C_{27}(1,3,6,8,10))$ = $\{C_{27}(1,3,6,8,10), C_{27}(2,6,7,11,12), C_{27}(3,4,5,12,13)\}$. 

$\Rightarrow$  $C_{27}(2,3,6,7,11),C_{27}(3,4,5,6,13)\notin Ad_{27}(C_{27}(1,3,6,8,10))$.

$\Rightarrow$ $C_{27}(1,3,6,8,10)$, $C_{27}(3,4,5,6,13)$,  $C_{27}(2,3,6,7,11)$ are Type-2 isomorphic w.r.t. $m$ = 3.

\item [\rm (6)]	$\theta_{27,3,1}(C_{27}(1,3,8,9,10))$ = $C_{27}(3,4,5,9,13)$ and $\theta_{27,3,2}(C_{27}(1,3,8,9,10))$ = $C_{27}(2,3,7,9,11)$.

$\Rightarrow$ $C_{27}(1,3,8,9,10)$ $\cong$ $C_{27}(3,4,5,9,13)$ $\cong$ $C_{27}(2,3,7,9,11)$. 
\\
Also,	$Ad_{27}(C_{27}(1,3,8,9,10))$ = $\{C_{27}(1,3,8,9,10), C_{27}(2,6,7,9,11), C_{27}(4,5,9,12,13)\}$. 

$\Rightarrow$  $C_{27}(2,3,7,9,11),C_{27}(3,4,5,9,13)\notin Ad_{27}(C_{27}(1,3,8,9,10))$.

$\Rightarrow$ $C_{27}(1,3,8,9,10)$, $C_{27}(3,4,5,9,13)$,  $C_{27}(2,3,7,9,11)$ are Type-2 isomorphic w.r.t. $m$ = 3.

\item [\rm (7)]	$\theta_{27,3,1}(C_{27}(1,3,8,10,12))$ = $C_{27}(3,4,5,12,13)$ $\&$ $\theta_{27,3,2}(C_{27}(1,3,8,10,12))$ = $C_{27}(2,3,7,11,12)$.

$\Rightarrow$ $C_{27}(1,3,8,10,12)$ $\cong$ $C_{27}(3,4,5,12,13)$ $\cong$ $C_{27}(2,3,7,11,12)$. 
\\
Also,	$Ad_{27}(C_{27}(1,3,8,10,12))$ = $\{C_{27}(1,3,8,10,12),  C_{27}(2,3,6,7,11), C_{27}(4,5,6,12,13)\}$. 

$\Rightarrow$  $C_{27}(2,3,7,11,12), C_{27}(3,4,5,12,13)\notin Ad_{27}(C_{27}(1,3,8,10,12))$.
\\
$\Rightarrow$ $C_{27}(1,3,8,10,12)$, $C_{27}(3,4,5,12,13)$,  $C_{27}(2,3,7,11,12)$ are Type-2 isomorphic w.r.t. $m$ = 3.

\item [\rm (8)]	$\theta_{27,3,1}(C_{27}(1,6,8,9,10))$ = $C_{27}(4,5,6,9,13)$ and $\theta_{27,3,2}(C_{27}(1,6,8,9,10))$ = $C_{27}(2,6,7,9,11)$.

$\Rightarrow$ $C_{27}(1,6,8,9,10)$ $\cong$ $C_{27}(4,5,6,9,13)$ $\cong$ $C_{27}(2,6,7,9,11)$. 
\\
Also,	$Ad_{27}(C_{27}(1,6,8,9,10))$ = $\{C_{27}(1,6,8,9,10), C_{27}(2,7,9,11,12), C_{27}(3,4,5,9,13)\}$. 

$\Rightarrow$  $C_{27}(2,6,7,9,11),C_{27}(4,5,6,9,13)\notin Ad_{27}(C_{27}(1,6,8,9,10))$.

$\Rightarrow$ $C_{27}(1,6,8,9,10)$, $C_{27}(4,5,6,9,13)$,  $C_{27}(2,6,7,9,11)$ are Type-2 isomorphic w.r.t. $m$ = 3.

\item [\rm (9)]	$\theta_{27,3,1}(C_{27}(1,6,8,10,12))$ = $C_{27}(4,5,6,12,13)$ $\&$ $\theta_{27,3,2}(C_{27}(1,6,8,10,12))$ = $C_{27}(2,6,7,11,12)$.

$\Rightarrow$ $C_{27}(1,6,8,10,12)$ $\cong$ $C_{27}(4,5,6,12,13)$ $\cong$ $C_{27}(2,6,7,11,12)$. 
\\
Also,	$Ad_{27}(C_{27}(1,6,8,10,12))$ = $\{C_{27}(1,6,8,10,12), C_{27}(2,3,7,11,12), C_{27}(3,4,5,6,13)\}$. 

$\Rightarrow$  $C_{27}(2,6,7,11,12), C_{27}(4,5,6,12,13)\notin Ad_{27}(C_{27}(1,6,8,10,12))$.
\\
$\Rightarrow$ $C_{27}(1,6,8,10,12)$, $C_{27}(4,5,6,12,13)$,  $C_{27}(2,6,7,11,12)$ are Type-2 isomorphic w.r.t. $m$ = 3.

\item [\rm (10)] $\theta_{27,3,1}(C_{27}(1,8,9,10,12))$ = $C_{27}(4,5,9,12,13)$ $\&$ $\theta_{27,3,2}(C_{27}(1,8,9,10,12))$ = $C_{27}(2,7,9,11,12)$.

$\Rightarrow$ $C_{27}(1,8,9,10,12)$ $\cong$ $C_{27}(4,5,9,12,13)$ $\cong$ $C_{27}(2,7,9,11,12)$. 
\\
Also,	$Ad_{27}(C_{27}(1,8,9,10,12))$ = $\{C_{27}(1,8,9,10,12), C_{27}(2,3,7,9,11), C_{27}(4,5,6,9,13)\}$. 

$\Rightarrow$  $C_{27}(2,7,9,11,12), C_{27}(4,5,9,12,13)\notin Ad_{27}(C_{27}(1,8,9,10,12))$.
\\
$\Rightarrow$ $C_{27}(1,8,9,10,12)$, $C_{27}(4,5,9,12,13)$,  $C_{27}(2,7,9,11,12)$ are Type-2 isomorphic w.r.t. $m$ = 3.
		
\item [\rm (11)] $\theta_{27,3,1}(C_{27}(1,3,6,8,9,10))$ = $C_{27}(3,4,5,6,9,13)$ and
\\
  $\theta_{27,3,2}(C_{27}(1,3,6,8,9,10))$ = $C_{27}(2,3,6,7,9,11)$.

$\Rightarrow$ $C_{27}(1,3,6,8,9,10)$ $\cong$ $C_{27}(3,4,5,6,9,13)$ $\cong$ $C_{27}(2,3,6,7,9,11)$. Also,
\\
$Ad_{27}(C_{27}(1,3,6,8,9,10))$ = $\{C_{27}(1,3,6,8,9,10), C_{27}(2,6,7,9,11,12), C_{27}(3,4,5,9,12,13)\}$. 

$\Rightarrow$  $C_{27}(2,3,6,7,9,11),C_{27}(3,4,5,6,9,13)\notin Ad_{27}(C_{27}(1,3,6,8,9,10))$.

$\Rightarrow$ $C_{27}(1,3,6,8,9,10)$, $C_{27}(3,4,5,6,9,13)$,  $C_{27}(2,3,6,7,9,11)$ are Type-2 isomorphic w.r.t. $m$ = 3.
	
\item [\rm (13)] $\theta_{27,3,1}(C_{27}(1,3,8,9,10,12))$ = $C_{27}(3,4,5,9,12,13)$ and
\\
  $\theta_{27,3,2}(C_{27}(1,3,8,9,10,12))$ = $C_{27}(2,3,7,9,11,12)$.

$\Rightarrow$ $C_{27}(1,3,8,9,10,12)$ $\cong$ $C_{27}(3,4,5,9,12,13)$ $\cong$ $C_{27}(2,3,7,9,11,12)$. Also,
\\
$Ad_{27}(C_{27}(1,3,8,9,10,12))$ = $\{C_{27}(1,3,8,9,10,12), C_{27}(2,3,6,7,9,11), C_{27}(4,5,6,9,12,13)\}$. 

$\Rightarrow$  $C_{27}(2,3,7,9,11,12), C_{27}(3,4,5,9,12,13)\notin Ad_{27}(C_{27}(1,3,8,9,10,12))$.

This implies that $C_{27}(1,3,8,9,10,12)$, $C_{27}(3,4,5,9,12,13)$,  $C_{27}(2,3,7,9,11,12)$ are Type-2 isomorphic w.r.t. $m$ = 3.
		
\item [\rm (14)] $\theta_{27,3,1}(C_{27}(1,6,8,9,10,12))$ = $C_{27}(4,5,6,9,12,13)$ and
\\
 $\theta_{27,3,2}(C_{27}(1,6,8,9,10,12))$ = $C_{27}(2,6,7,9,11,12)$.

$\Rightarrow$ $C_{27}(1,6,8,9,10,12)$ $\cong$ $C_{27}(4,5,6,9,12,13)$ $\cong$ $C_{27}(2,6,7,9,11,12)$. Also, 
\\
$Ad_{27}(C_{27}(1,6,8,9,10,12))$ = $\{C_{27}(1,6,8,9,10,12), C_{27}(2,3,7,9,11,12), C_{27}(3,4,5,6,9,13)\}$. 

$\Rightarrow$  $C_{27}(2,6,7,9,11,12), C_{27}(4,5,6,9,12,13) \notin Ad_{27}(C_{27}(1,6,8,9,10,12))$.

$\Rightarrow$ $C_{27}(1,6,8,9,10,12)$, $C_{27}(4,5,6,9,12,13)$,  $C_{27}(2,6,7,9,11,12)$ are Type-2 isomorphic w.r.t. $m$ = 3.
\end{enumerate}	
  
  Thus, there are 12 triples of Type-2 isomorphic circulant graphs of order 27 and all these triples of isomorphic circulant graphs are Type-2 isomorphic w.r.t. m = 3. Type-2 isomorphic circulant graphs $C_{27}(1,6,8,10)$, $C_{27}(2,6,7,11)$ and $C_{27}(4,5,6,13)$ w.r.t. $m$ = 3 are given in Figures 7, 8 and 9 and in Figures 10 and 11 circulant graphs $\theta_{27,3,1}(C_{27}(1,6,8,10))$ = $C_{27}(4,5,6,13)$ and $\theta_{27,3,2}(C_{27}(1,6,8,10))$ = $C_{27}(2,6,7,11)$ are given. \hfill $\Box$   
  \begin{figure}[ht]
  	\centerline{\includegraphics[width=6.2in]{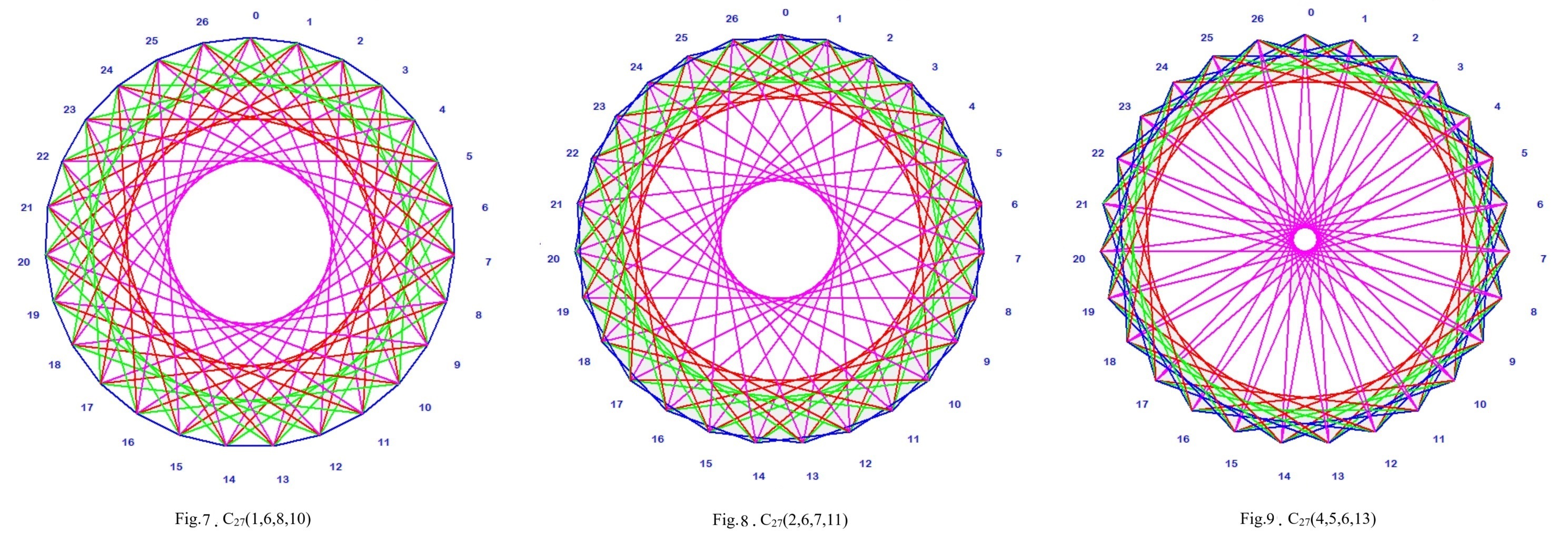}}
  \end{figure}
\begin{figure}[ht]
	\hspace{5cm}\includegraphics[width=4.4in]{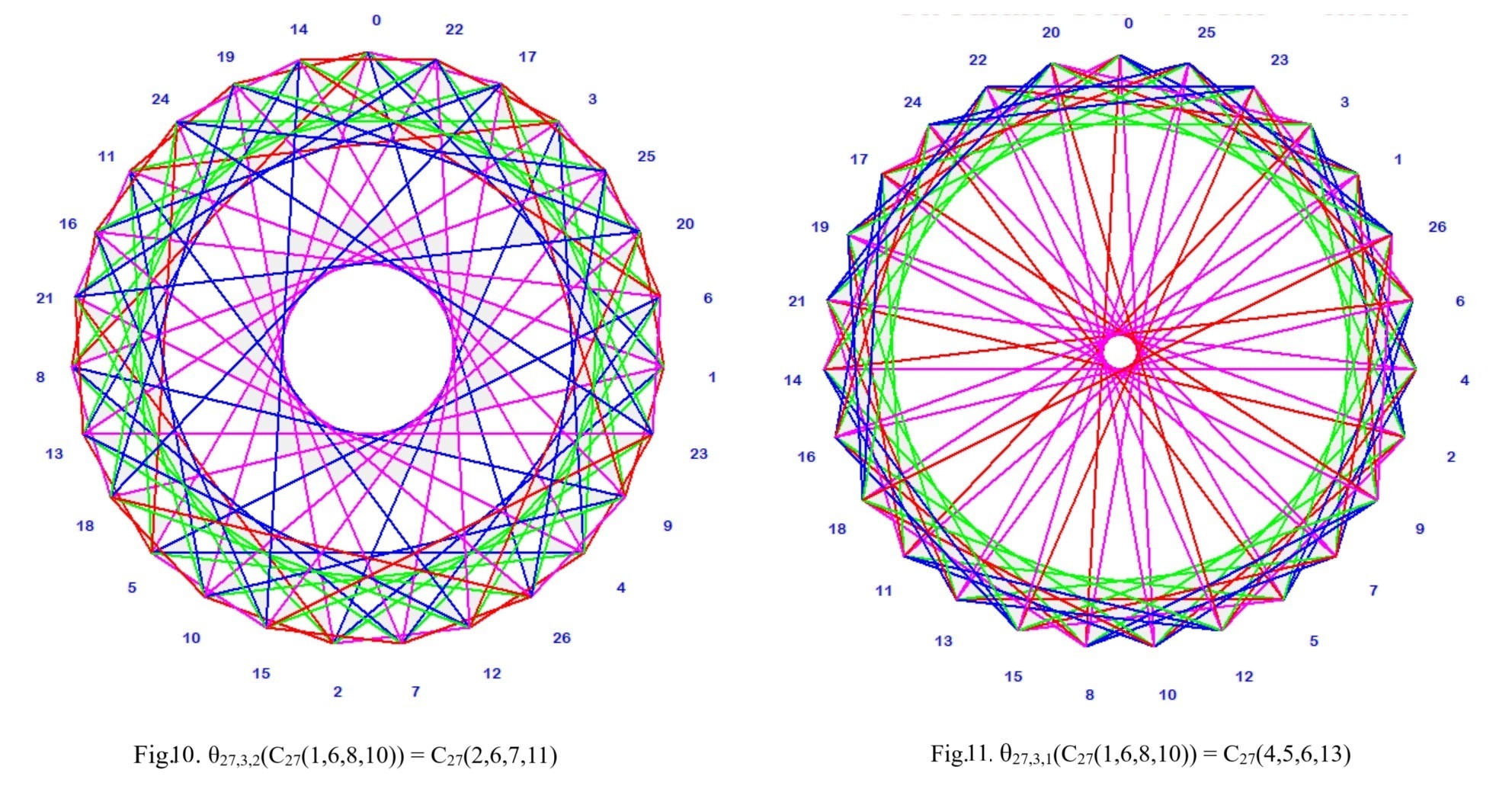}
\end{figure}

\section{Conclusion}

 The author feels that this paper provides a lot of scope for further research to obtain more families of Type-2 isomorphic circulant graphs.
 
\vspace{.1cm}
\noindent
\textbf{Declaration of competing interest}\quad 
The author declares that he has no conflict of interest.


\end{document}